# Smarandache semirings, semifields and semivector spaces

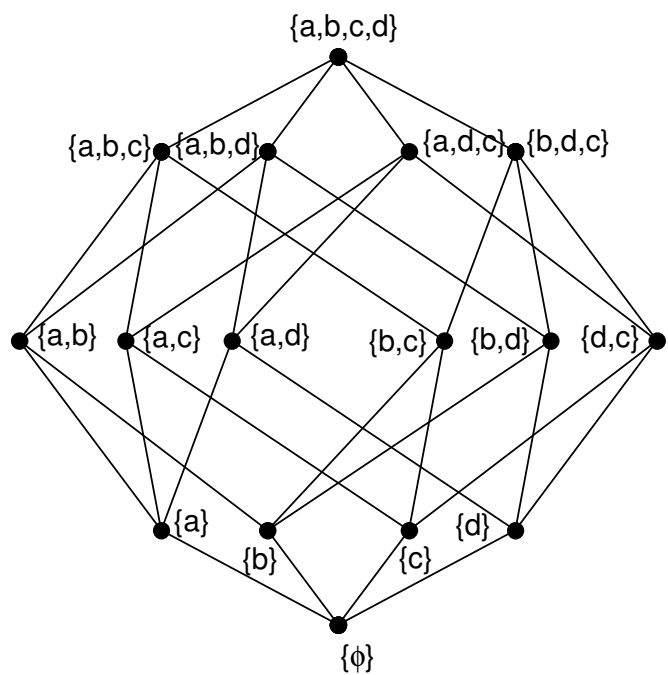

**W. B. VASANTHA KANDASAMY**

# Smarandache Semirings, Semifields and Semivector spaces

W. B. Vasantha Kandasamy
Department of Mathematics
Indian Institute of Technology, Madras
Chennai – 600036, India
Email: vasantha@iitm.ac.in
web: http://mat.iitm.ac.in/~wbv

*The picture on the cover is a Boolean algebra constructed using the power set P(X) where X = {a, b, c} which is a finite Smarandache semiring of order 16.*

# CONTENTS

# PREFACE

Smarandache notions, which can be undoubtedly characterized as revolutionary mathematics, has the capacity of being utilized to analyse, study and introduce, naturally, the concepts of several structures by means of extension or identification as a substructure. Several researchers around the world working on Smarandache notions have systematically carried out this study. This is the first book on the Smarandache algebraic structures that have two binary operations.

Semirings are algebraic structures with two binary operations enjoying several properties and it is the most generalized structure — for all rings and fields are semirings. The study of this concept is very meagre except for a very few research papers. Now, when we define Smarandache semirings (S-semiring), we make the richer structure of semifield to be contained in a S-semiring; and this S-semiring is of the first level. To have the second level of S-semirings, we need a still richer structure, viz. field to be a subset in a S-semiring. This is achieved by defining a new notion called the Smarandache mixed direct product. Likewise we also define the Smarandache semifields of level II. This study makes one relate, compare and contrasts weaker and stronger structures of the same set.

The motivation for writing this book is two-fold. First, it has been our aim to give an insight into the Smarandache semirings, semifields and semivector spaces. Secondly, in order to make an organized study possible, we have also introduced all the concepts about semirings, semifields and semivector spaces; since, to the best of our knowledge we do not have books, which solely deals with these concepts. This book introduces several new concepts about Smarandache semirings, semifields and semivector spaces. We assume at the outset that the reader has a strong background in algebra that will enable one to follow and understand the book completely.

This book consists of seven chapters. The first chapter introduces the basic concepts which are very essential to make the book self-contained. The second chapter is solely devoted to the introduction of semirings and its properties. The notions about semifields and semivector spaces are introduced in the third chapter. Chapter four, which is one of the major parts of this book contains a complete systematic introduction of all concepts together with a sequential analysis of these concepts. Examples are provided abundantly to make the abstract definitions and results easy and explicit to the reader.

Further, we have also given several problems as exercises to the student/ researcher, since it is felt that tackling these research problems is one of the ways to get deeply involved in the study of Smarandache semirings, semifields and semivector spaces. The fifth chapter introduces Smarandache semifields and elaborates some of its properties. The concept of Smarandache semivector spaces are introduced and analysed in chapter six. The final chapter includes 25 research problems and they will certainly be a boon to any researcher. It is also noteworthy to mention that at the end of every chapter we have provided a bibliographical list for supplementary reading, since referring to and knowing these concepts will equip and enrich the researcher's knowledge. The book also contains a comprehensive index.

Finally, following the suggestions and motivations of Dr. Minh Perez we have introduced the Smarandache anti semiring, anti semifield and anti semivector space. On his suggestion I have at each stage introduced II level of Smarandache semirings and semifields. Overall in this book we have totally defined 65 concepts related to the Smarandache notions in semirings and its generalizations.

I deeply acknowledge my children Meena and Kama whose joyful persuasion and support encouraged me to write this book.

# CHAPTER ONE
# PRELIMINARY NOTIONS

This chapter gives some basic notions and concepts used in this book to make this book self-contained. The serious study of semirings is very recent and to the best of my knowledge we do not have many books on semirings or semifields or semivector spaces. The purpose of this book is two-fold, firstly to introduce the concepts of semirings, semifields and semivector spaces (which we will shortly say as semirings and its generalizations), which are not found in the form of text. Secondly, to define Smarandache semirings, semifields and semivector spaces and study these newly introduced concepts.

In this chapter we recall some basic properties of semigroups, groups, lattices, Smarandache semigroups, fields, vector spaces, group rings and semigroup rings. We assume at the outset that the reader has a good background in algebra.

## 1.1 Semigroups, Groups and Smarandache Semigroups

In this section we just recall the definition of these concepts and give a brief discussion about these properties.

**DEFINITION 1.1.1**: *Let S be a non-empty set, S is said to be a semigroup if on S is defined a binary operation '$*$' such that*

1. *For all a, b $\in$ S we have a $*$ b $\in$ S (closure).*
2. *For all a, b, c $\in$ S we have (a $*$ b) $*$ c = a $*$ (b $*$ c) (associative law), We denote by (S, $*$) the semigroup.*

**DEFINITION 1.1.2**: *If in a semigroup (S, $*$), we have a $*$ b = b $*$ a for all a, b $\in$ S we say S is a commutative semigroup.*

*If the number of elements in the semigroup S is finite we say S is a finite semigroup or a semigroup of finite order, otherwise S is of infinite order. If the semigroup S contains an element e such that e $*$ a = a $*$ e = a for all a $\in$ S we say S is a semigroup with identity e or a monoid. An element x $\in$ S, S a monoid is said to be invertible or has an inverse in S if there exist a y $\in$ S such that xy = yx = e.*

**DEFINITION 1.1.3**: *Let (S, $*$) be a semigroup. A non-empty subset H of S is said to be a subsemigroup of S if H itself is a semigroup under the operations of S.*

**DEFINITION 1.1.4**: *Let (S, $*$) be a semigroup, a non-empty subset I of S is said to be a right ideal of S if I is a subsemigroup of S and for all s $\in$ S and i $\in$ I we have is $\in$ I.*

*Similarly one can define left ideal in a semigroup. We say I is an ideal of a semigroup if I is simultaneously a left and a right ideal of S.*

**DEFINITION 1.1.5**: *Let (S, ∗) and (S₁, o) be two semigroups. We say a map ϕ from (S, ∗) → (S₁, o) is a semigroup homomorphism if ϕ(s₁ ∗ s₂) = ϕ(s₁) o ϕ(s₂) for all s₁, s₂ ∈ S.*

*Example 1.1.1*: $Z_9$ = {0, 1, 2, … , 8} is a commutative semigroup of order nine under multiplication modulo 9 with unit.

*Example 1.1.2*: S = {0, 2, 4, 6, 8, 10} is a semigroup of finite order, under multiplication modulo 12. S has no unit but S is commutative.

*Example 1.1.3*: Z be the set of integers. Z under usual multiplication is a semigroup with unit of infinite order.

*Example 1.1.4*: 2Z = {0, ±2, ±4, … , ±2n …} is an infinite semigroup under multiplication which is commutative but has no unit.

*Example 1.1.5*: Let $S_{2\times 2} = \left\{ \begin{pmatrix} a & b \\ c & d \end{pmatrix} \bigg/ a,b,c,d \in Z_4 \right\}$. $S_{2\times 2}$ is a finite non-commutative semigroup under matrix multiplication modulo 4, with unit $I_{2\times 2} = \begin{pmatrix} 1 & 0 \\ 0 & 1 \end{pmatrix}$.

*Example 1.1.6*: Let $M_{2\times 2} = \left\{ \begin{pmatrix} a & b \\ c & d \end{pmatrix} \bigg/ a,b,c,d \in Q, \text{the field of rationals} \right\}$. $M_{2\times 2}$ is a non-commutative semigroup of infinite order under matrix multiplication with unit $I_{2\times 2} = \begin{pmatrix} 1 & 0 \\ 0 & 1 \end{pmatrix}$.

*Example 1.1.7*: Let Z be the semigroup under multiplication pZ = {0, ± p, ± 2p, …} is an ideal of Z, p any positive integer.

*Example 1.1.8*: Let $Z_{14}$ = {0, 1, 2, … , 13} be the semigroup under multiplication. Clearly I = {0, 7} and J = {0, 2, 4, 6, 8, 10, 12} are ideals of $Z_{14}$.

*Example 1.1.9*: Let X = {1, 2, 3, … , n} where n is a finite integer. Let S (n) denote the set of all maps from the set X to itself. Clearly S (n) is a semigroup under the composition of mappings. S(n) is a non-commutative semigroup with $n^n$ elements in it; in fact S(n) is a monoid as the identity map is the identity element under composition of mappings.

*Example 1.1.10*: Let S(3) be the semigroup of order 27, (which is for n = 3 described in example 1.1.9.) It is left for the reader to find two sided ideals of S(3).

**Notation**: Throughout this book S(n) will denote the semigroup of mappings of any set X with cardinality of X equal to n. Order of S(n) is denoted by o (S(n)) or |S(n)| and S(n) has $n^n$ elements in it.

Now we just recall the definition of group and its properties.

**DEFINITION 1.1.6**: *A non-empty set of elements G is said to from a group if in G there is defined a binary operation, called the product and denoted by '•' such that*

1. *a, b ∈ G implies a • b ∈ G (Closure property)*
2. *a, b, c ∈ G implies a • (b • c) = (a • b) • c (associative law)*
3. *There exists an element e ∈ G such that a • e = e • a = a for all a ∈ G (the existence of identity element in G).*
4. *For every a ∈ G there exists an element $a^{-1}$ ∈ G such that $a • a^{-1} = a^{-1} • a = e$ (the existence of inverse in G).*

*A group G is abelian or commutative if for every a, b ∈ G a • b = b • a. A group, which is not abelian, is called non-abelian. The number of distinct elements in G is called the order of G; denoted by o (G) = |G|. If o (G) is finite we say G is of finite order otherwise G is said to be of infinite order.*

**DEFINITION 1.1.7**: *Let (G, o) and ($G_1$, ∗) be two groups. A map φ: G to $G_1$ is said to be a group homomorphism if φ (a • b) = φ(a) ∗ φ(b) for all a, b ∈ G.*

**DEFINITION 1.1.8**: *Let (G, ∗) be a group. A non-empty subset H of G is said to be a subgroup of G if (H, ∗) is a group, that is H itself is a group.*

For more about groups refer. (I. N. Herstein and M. Hall).

Throughout this book by $S_n$ we denote the set of all one to one mappings of the set X = {$x_1$, … , $x_n$} to itself. The set $S_n$ together with the composition of mappings as an operation forms a non-commutative group. This group will be addressed in this book as symmetric group of degree n or permutation group on n elements. The order of $S_n$ is finite, only when n is finite. Further $S_n$ has a subgroup of order n!/2 , which we denote by $A_n$ called the alternating group of $S_n$ and S = $Z_p$ \ {0} when p is a prime under the operations of usual multiplication modulo p is a commutative group of order p-1.

Now we just recall the definition of Smarandache semigroup and give some examples. As this notion is very new we may recall some of the important properties about them.

**DEFINITION 1.1.9**: *The Smarandache semigroup (S-semigroup) is defined to be a semigroup A such that a proper subset of A is a group. (with respect to the same induced operation).*

**DEFINITION 1.1.10**: *Let S be a S semigroup. If every proper subset of A in S, which is a group is commutative then we say the S-semigroup S to be a Smarandache commutative semigroup and if S is a commutative semigroup and is a S-semigroup then obviously S is a Smarandache commutative semigroup.*

*Let S be a S-semigroup, o(S) = number of elements in S that is the order of S, if o(S) is finite we say S is a finite S-semigroup otherwise S is an infinite S-semigroup.*

***Example 1.1.11***: Let $Z_{12} = \{0, 1, 2, \ldots, 11\}$ be the modulo integers under multiplication mod 12. $Z_{12}$ is a S-semigroup for the sets $A_1 = \{1, 5\}$, $A_2 = \{9, 3\}$, $A_3 = \{4, 8\}$ and $A_4 = \{1, 5, 7, 11\}$ are subgroups under multiplication modulo 12.

***Example 1.1.12***: Let $S(5)$ be the symmetric semigroup. $S_5$ the symmetric group of degree 5 is a proper subset of $S(5)$ which is a group. Hence $S(5)$ is a Smarandache semigroup.

***Example 1.1.13***: Let $M_{n \times n} = \{(a_{ij}) \,/\, a_{ij} \in Z\}$ be the set of all $n \times n$ matrices; under matrix multiplication $M_{n \times n}$ is a semigroup. But $M_{n \times n}$ is a S-semigroup if we take $P_{n \times n}$ the set of all is a non-singular matrixes of $M_{n \times n}$, it is a group under matrix multiplication.

***Example 1.1.14***: $Z_p = \{0, 1, 2, \ldots, p-1\}$ is a semigroup under multiplication modulo p. The set $A = \{1, p-1\}$ is a subgroup of $Z_p$. Hence $Z_p$ for all primes p is a S-semigroup.

For more about S-semigroups one can refer [9,10,11,18].

**PROBLEMS**:

1. For the semigroup $S_{3 \times 3} = \{(a_{ij}) \,/\, a_{ij} \in Z_2 = \{0, 1\}\}$; (the set of all 3×3 matirixes with entries from $Z_2$) under multiplication.

    i. Find the number of elements in $S_{3 \times 3}$.
    ii. Find all the ideals of $S_{3 \times 3}$.
    iii. Find only the right ideals of $S_{3 \times 3}$.
    iv. Find all subsemigroups of $S_{3 \times 3}$.

2. Let $S(21)$ be the set of all mappings of a set $X = \{1, 2, \ldots, 21\}$ with 21 elements to itself. $S(21) = S(X)$ is a semigroup under composition of mappings.

    i. Find all subsemigroups of $S(X)$ which are not ideals.
    ii. Find all left ideals of $S(X)$.
    iii. How many two sided ideals does $S(X)$ contain?

3. Find all the ideals of $Z_{28} = \{0, 1, 2, \ldots, 27\}$, the semigroup under multiplication modulo 28.

4. Construct a homomorphism between the semigroups. $S_{3 \times 3}$ given in problem 1 and $Z_{28}$ given in problem 3. Find the kernel of this homomorphism. ($\phi: S_{3 \times 3} \to Z_{28}$) where ker $\phi = \{x \in S_{3 \times 3} \,/\, \phi(x) = 1\}$.

5. Does there exist an isomorphism between the semigroups $S(4)$ and $Z_{256}$? Justify your answer with reasons.

6. Find all the right ideals of $S(5)$. Can $S(5)$ have ideals of order 120?

7. Find all the subgroups of $S_4$.

8. Does there exist an isomorphism between the groups $G = \langle g/g^6 = 1 \rangle$ and $S_3$?

9. Can you construct a group homomorphism between $g = \langle g/g^6 = 1 \rangle$ and $S_4$? Prove or disprove.

10. Does there exist a group homomorphism between $G = \langle g/g^{11} = 1 \rangle$ and the symmetric group $S_5$?

11. Can we have a group homomorphism between $G = \langle g/g^p = 1 \rangle$ and the symmetric group $S_q$ (where p and q are two distinct primes)?

12. Find a group homomorphism between $D_{2n}$ and $S_n$. ($D_{2n}$ is called the dihedral group of order 2n given by the following relation, $D_{2n} = \{a, b/ a^2 = b^n = 1;$ bab $= a\}$.

13. Give an example of a S-semigroup of order 7 (other than $Z_7$).

14. Does a S-semigroup of order 2 exist? Justify!

15. Find a S-semigroup of order 16.

16. Find all subgroups of the S-semigroup, $Z_{124} = \{0, 1, 2, \ldots, 123\}$ under multiplication modulo 124.

17. Can $Z_{25} = \{0, 1, 2, \ldots, 24\}$ have a subset of order 6 which is a group ($Z_{25}$ is a semigroup under usual multiplication modulo 25)?

18. Let $Z_{121} = \{0, 1, 2, \ldots, 120\}$ be the semigroup under multiplication modulo 121. Can $Z_{121}$ have subgroups of even order? If so find all of them.

## 1.2 Lattices

In this section we just recall the basic results about lattices used in this book.

**DEFINITION 1.2.1**: *Let A and B be non-empty sets. A relation R from A to B is a subset of A × B. Relations from A to B are called relations on A, for short, if (a, b) ∈ R then we write aRb and say that 'a is in relation R to b'. Also if a is not in relation R to b, we write aRb.*

*A relation R on a non-empty set A may have some of the following properties:*

> *R is reflexive if for all a in A we have aRa.*
> *R is symmetric if for a and b in A: aRb implies bRa.*
> *R is anti symmetric if for all a and b in A; aRb and bRa imply a = b.*
> *R is transitive if for a, b, c in A; aRb and bRc imply aRc.*

*A relation R on A is an equivalence relation if R is reflexive, symmetric and transitive. In the case [a] = {b ∈ A| aRb}, is called the equivalence class of a for any a ∈ A.*

**DEFINITION 1.2.2**: *A relation R on a set A is called a partial order (relation) if R is reflexive, anti symmetric and transitive. In this case (A, R) is a partially ordered set or poset.*

We denote the partial order relation by ≤ or ⊆.

**DEFINITION 1.2.3**: *A partial order relation ≤ on A is called a total order if for each a, b ∈ A, either a ≤ b or b ≤ a. {A, ≤} is called a chain or a totally ordered set.*

*Example 1.2.1*: Let A = {1, 2, 3, 4, 7}, (A, ≤) is a total order. Here '≤' is the usual "less than or equal to" relation.

*Example 1.2.2*: Let X = {1, 2, 3}, the power set of X is denoted by P(X) = {ϕ, X, {a}, {b}, {c}, {a, b}, {c, b}, {a, c}}. P(X) under the relation ' ⊆ ' "inclusion" as subsets or containment relation is a partial order on P(X).

It is important or interesting to note that finite partially ordered sets can be represented by Hasse Diagrams. Hasse diagram of the poset A given in example 1.2.1:

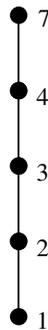

**Figure 1.2.1**

Hasse diagram of the poset P(X) described in example 1.2.2 is as follows:

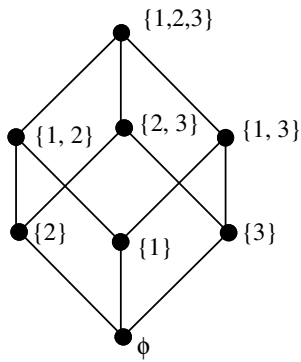

**Figure 1.2.2**

**DEFINITION 1.2.4**: *Let (A, ≤) be a poset and B ⊆ A.*

i)     *a ∈ A is called an upper bound of B if and only if for all b ∈ B, b ≤ a.*
ii)    *a ∈ A is called a lower bound of B if and only if for all b ∈ B; a ≤ b.*
iii)   *The greatest amongst the lower bounds, whenever it exists is called the infimum of B, and is denoted by inf B.*
iv)   *The least upper bound of B whenever it exists is called the supremum of B and is denoted by sup B.*

Now with these notions and notations we define a semilattice.

**DEFINITION 1.2.5**: *A poset (L, ≤) is called a semilattice order if for every pair of elements x, y in L the sup (x, y) exists (or equivalently we can say inf (x, y) exist).*

**DEFINITION 1.2.6**: *A poset (L, ≤) is called a lattice ordered if for every pair of elements x, y in L the sup (x, y) and inf (x, y) exists.*

It is left for the reader to verify the following result.

**Result:**
1. Every ordered set is lattice ordered.
2. In a lattice ordered set (L, ≤) the following statements are equivalent for all x and y in L.
    a. $x \leq y$
    b. Sup $(x, y) = y$
    c. Inf $(x, y) = x$.

Now as this text uses also the algebraic operations on a lattice we define an algebraic lattice.

**DEFINITION 1.2.7**: *An algebraic lattice (L, ∩, ∪) is a non-empty set L with two binary operations ∪ (join) and ∩ (meet) (also called union or sum and intersection or product respectively) which satisfy the following conditions for all x, y, z ∈ L.*

$L_1.$ $x \cap y = y \cap x,$                 $x \cup y = y \cup x$
$L_2.$ $x \cap (y \cap z) = (x \cap y) \cap z,$     $x \cup (y \cup z) = (x \cup y) \cup z$
$L_3.$ $x \cap (x \cup y) = x,$            $x \cup (x \cap y) = x.$

*Two applications of $L_3$ namely $x \cap x = x \cap (x \cup (x \cap x)) = x$ lead to the additional condition $L_4.$ $x \cap x = x$, $x \cup x = x$. $L_1$ is the commutative law, $L_2$ is the associative law, $L_3$ is the absorption law and $L_4$ is the idempotent law.*

The connection between lattice ordered sets and algebraic lattices is as follows:

**Result:**
1. Let (L, ≤) be a lattice ordered set. If we define $x \cap y = $ inf $(x, y)$ and $x \cup y = $ sup $(x, y)$ then (L, ∪, ∩) is an algebraic lattice.

2. Let (L, ∪, ∩) be an algebraic lattice. If we define x ≤ y if and only if x ∩ y = x (or x ≤ y if and only if x ∪ y = y) then (L, ≤) is a lattice ordered set.

This result is left as an exercise for the reader to verify.

Thus it can be verified that the above result yields a one to one relationship between algebraic lattices and lattice ordered sets. Therefore we shall use the term lattice for both concepts. |L| = o(L) denotes the order (that is cardinality) of the lattice L.

***Example 1.2.3***: $L_4^5$ be a lattice given by the following Hasse diagram:

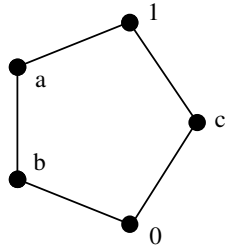

**Figure 1.2.3**

This lattice will be called as the pentagon lattice in this book.

***Example 1.2.4***: Let $L_3^5$ be the lattice given by the following Hasse diagram:

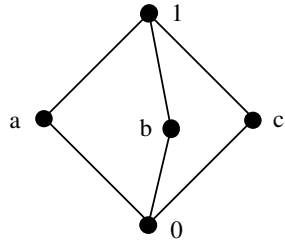

**Figure 1.2.4**

This lattice will be addressed as diamond lattice in this book.

**DEFINITION 1.2.8**: *Let (L, ≤) be a lattice. If '≤' is a total order on L and L is lattice order we call L a chain lattice. Thus we see in a chain lattice L we have for every pair a, b ∈ L we have either a ≤ b or b ≤ a.*

Chain lattices will play a major role in this book.

***Example 1.2.5***: Let L be [a, b] any closed interval on the real line, [a, b] under the total order is a chain lattice.

***Example 1.2.6***: [0,∞) is also a chain lattice of infinite order. Left for the reader to verify.

*Example 1.2.7*: [-∞, 1] is a chain lattice of infinite cardinality.

*Example 1.2.8*: Take [0, 1] = L the two element set. L is the only 2 element lattice and it is a chain lattice having the following Hasse diagram and will be denoted by $C_2$.

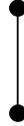

**Figure 1.2.5**

**DEFINITION 1.2.9**: *A non-empty subset S of a lattice L is called a sublattice of L if S is a lattice with respect to the restriction of ∩ and ∪ of L onto S.*

**DEFINITION 1.2.10**: *Let L and M be any two lattices. A mapping f: L → M is called a*

1. *Join homomorphism if x ∪ y = z ⇒ f(x) ∪ f(y) = f(z)*
2. *Meet homomorphism if x ∩ y = z ⇒ f(x) ∩ f(y) = f(z)*
3. *Order homomorphism if x ≤ y imply f(x) ≤ f(y) for all x, y ∈ L.*

*f is a lattice homomorphism if it is both a join and a meet homomorphism. Monomorphism, epimorphism, isomorphism of lattices are defined as in the case of other algebraic structures.*

**DEFINITION 1.2.11**: *A lattice L is called modular if for all x, y, z ∈ L, x ≤ z imply x ∪ (y ∩ z) = (x ∪ y) ∩ z.*

**DEFINITION 1.2.12**: *A lattice L is called distributive if either of the following conditions hold good for all x, y, z in L. x ∪ (y ∩ z) = (x ∪ y) ∩ (x ∪ z) or x ∩ (y ∪ z) = (x ∩ y) ∪ (x ∩ z) called the distributivity equations.*

It is left for reader to verify the following result:

**Result**: A lattice L is distributive if and only if for all x, y, z ∈ L. (x ∩ y) ∪ (y ∩ z) ∪ (z ∩ x) = (x ∪ y) ∩ (y ∪ z) ∩ (z ∪ x).

*Example 1.2.9*: The following lattice L given by the Hasse diagram is distributive.

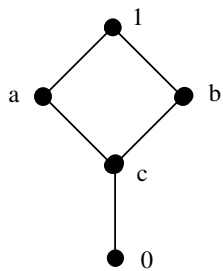

**Figure 1.2.6**

*Example 1.2.10*:

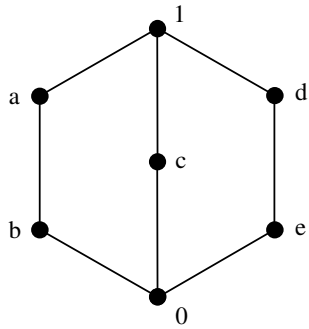

**Figure 1.2.7**

This lattice is non-distributive left for the reader to verify.

*Example 1.2.11*: Prove P(X) the power set of X where X = (1, 2) is a lattice with 4 elements given by the following Hasse diagram:

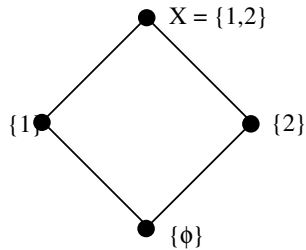

**Figure 1.2.8**

*Example 1.2.12*: The lattice

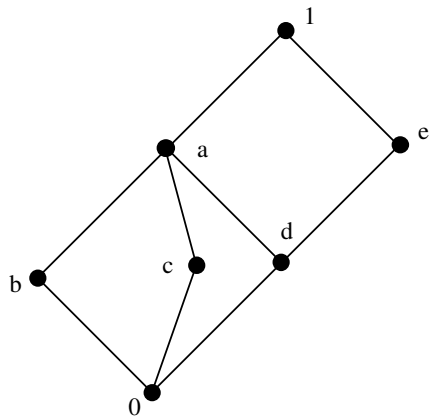

**Figure 1.2.9**

is modular and not distributive. Left for the reader to verify.

**DEFINITION 1.2.13**: *A lattice L with 0 and 1 is called complemented if for each $x \in L$ there is atleast one element y such that $x \cap y = 0$ and $x \cup y = 1$, y is called a complement of x.*

*Example 1.2.13*: The lattice with the following Hasse diagram:

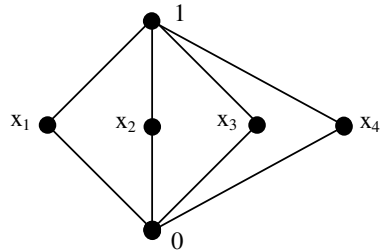

**Figure 1.2.10**

is such that $x_i \cap x_j = 0$, $i \neq j$; $x_i \cup x_j = 1$, $i \neq j$, each $x_i$ has a complement $x_j$, $i \neq j$.

**Result**: If L is a distributive lattice then each $x \in L$ has atmost one complement which is denoted by x'. This is left for the reader to verify.

**DEFINITION 1.2.14**: *A complemented distributive lattice is called a Boolean algebra (or a Boolean lattice). Distributivity in a Boolean algebra guarantees the uniqueness of complements.*

**DEFINITION 1.2.15**: *Let $B_1$ and $B_2$ be two Boolean algebras. The mapping $\phi: B_1 \to B_2$ is called a Boolean algebra homomorphism if $\phi$ is a lattice homomorphism and for all $x \in B_1$, we have $\phi(x') = (\phi(x))'$.*

*Example 1.2.14*: Let $X = \{x_1, x_2, x_3, x_4\}$. P(X) = power set of X, is a Boolean algebra with 16 elements in it. This is left for the reader to verify.

**PROBLEMS**:

1. Prove the diamond lattice is non-distributive but modular.

2. Prove the pentagon lattice is non-distributive and non-modular.

3. Prove the lattice with Hasse diagram is non-modular.

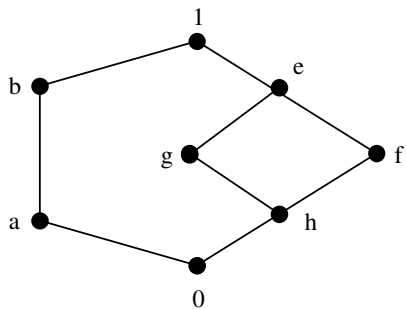

**Figure 1.2.11**

4. Find all sublattices of the lattice given in Problem 3. Does this lattice contain the pentagon lattice as a sublattice?

5. Prove all chain lattices are distributive.

6. Prove all lattices got from the power set of a set is distributive.

7. Prove a lattice L is distributive if and only if for all x, y, z ∈ L, $x \cap y = x \cap z$ and $x \cup y = x \cup z$ imply $y = z$.

8. Prove for any set X with n elements P(X), the power set of X is a Boolean algebra with $2^n$ elements in it.

9. Is the lattice with the following Hasse diagram, distributive? complemented? modular?

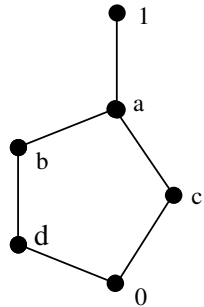

**Figure 1.2.12**

10. Prove for any lattice L without using the principle of duality the following conditions are equivalent.

    1. ∀ a, b, c ∈ L, $(a \cup b) \cap c = (a \cap c) \cup (b \cap c)$
    2. ∀ a, b, c ∈ L, $(a \cap b) \cup c = (a \cup c) \cap (b \cup c)$.

11. Prove if a, b, c are elements of a modular lattice L with the property $(a \cup b) \cap c = 0$ then $a \cap (b \cup c) = a \cap b$.

12. Prove in any lattice we have $[(x \cap y) \cup (x \cap z)] \cap [(x \cap y) \cup (y \cap z)] = x \cap y$ for all x, y, z ∈ L.

13. Prove a lattice L is modular if and only if for all x, y, z ∈ L, $x \cup (y \cap (x \cup z)) = (x \cup y) \cap (x \cup z)$.

14. Prove b has 2 complements a and c in the pentagon lattice given by the diagram

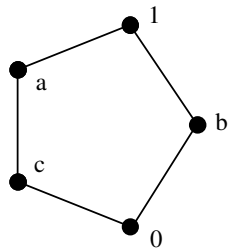

**Figure 1.2.13**

15. Give two examples of lattices of order 8 and 16, which are not Boolean algebras.

16. How many Boolean algebras are there with four elements 0, 1, a and b?

## 1.3 Rings and Fields

In this section we mainly introduce the concept of ring and field. This is done for two reasons, one to enable one to compare a field and a semifield. Second, to study group rings and semigroup rings. We do not give all the properties about field but what is essential alone is given, as the book assumes that the reader must have a good background of algebra.

**DEFINITION 1.3.1**: *A ring is a set R together with two binary operations + and •  called addition and multiplication, such that*

1. *(R, +) is an abelian group*
2. *The product r • s of any two elements r, s $\in$ R is in R and multiplication is associative.*
3. *For all r, s, t $\in$ R, r • (s + t) = r • s + r • t and (r + s) • t = r • t + s • t (distributive law).*

*We denote the ring by (R, +, •) or simply by R. In general the neutral element in (R,+) will always be denoted by 0, the additive inverse of r $\in$ R is –r. Instead of r • s we will denote it by rs. Clearly these rings are "associative rings". Let R be a ring, R is said to be commutative, if a • b = b • a for all a, b $\in$ R. If there is an element 1 $\in$ R such that 1 • r = r • 1 = r for all r $\in$ R, then 1 is called the identity (or unit) element.*

*If r • s = 0 implies r = 0 or s = 0 for all r, s $\in$ R then R is called integral. A commutative integral ring with identity is called an integral domain.*

*If R \ {0} is a group then R is called a skew field or a division ring. If more over, R is commutative we speak of a field. The characteristic of R is the smallest natural number k with kr = r + … + r (k times) equal to zero for all r $\in$ R. We then write k = characteristic R. If no such k exists we put characteristic R = 0.*

*Example 1.3.1*: Q be the set of all rationals. (Q, +, •) is a field of characteristic 0.

*Example 1.3.2*: Let Z be the set of integers. (Z, +, •) is a ring which is in fact an integral domain.

*Example 1.3.3*: Let $M_{n \times n}$ be the collection of all n × n matrices with entries from Q. $M_{n \times n}$ with matrix addition and matrix multiplication is a ring which is non-commutative and this ring has zero divisors that is, $M_{n \times n}$ is not a skew field or a division ring.

*Example 1.3.4*: Let $M'_{n \times n}$ denote the set of all non-singular matrices with entries from Q that is given in example 1.3.3. Clearly $M'_{n \times n}$ is a division ring or a skew field.

***Example 1.3.5***: Let R be the set of reals, R is a field of characteristic 0.

***Example 1.3.6***: Let $Z_{28} = \{0, 1, 2, \ldots, 27\}$. $Z_{28}$ with usual addition and multiplication modulo 28 is a ring. Clearly $Z_{28}$ is a commutative ring with $7.4 \equiv 0 \pmod{28}$ that is $Z_{28}$ has zero divisors.

***Example 1.3.7***: Let $Z_{23} = \{0, 1, 2, \ldots, 22\}$ be the ring of integers modulo 23. $Z_{23}$ is a field of characteristic 23.

**DEFINITION 1.3.2**: *Let F be a field. A proper subset S of F is said to be a subfield if S itself is a field under the operations of F.*

**DEFINITION 1.3.3**: *Let F be a field. If F has no proper subfields then F is said to be a prime field.*

***Example 1.3.8***: $Z_p = \{0, 1, 2, \ldots, p-1\}$ where p is a prime, is a prime field of characteristic p.

***Example 1.3.9***: Let Q be the field of rationals. Q has no proper subfield. Q is the prime field of characteristic 0; all prime fields of characteristic 0 are isomorphic to Q.

***Example 1.3.10***: Let R be the field of reals. R has the subset $Q \subset R$ and Q is a field; so R is not a prime field and characteristic $R = 0$.

**DEFINITION 1.3.4**: *Let R be any ring. A proper subset S of R is said to be a subring of R if S is a ring under the operations of R.*

***Example 1.3.11***: Let $Z_{20} = \{0, 1, 2, \ldots, 19\}$ is the ring of integers modulo 20. Clearly $A = \{0, 10\}$ is a subring of $Z_{20}$.

***Example 1.3.12***: Let Z be the ring of integers $5Z \subset Z$ is the subring of Z.

***Example 1.3.13***: Let R be a commutative ring and R[x] be the polynomial ring. $R \subset R[x]$ is a subring of R[x]. In fact R[x] is an integral domain if and only if R is an integral domain (left as an exercise for the reader to verify).

**DEFINITION 1.3.5**: *Let R and S be any two rings. A map $\phi: R \to S$ is said to be a ring homomorphism*
*if $\phi(a + b) = \phi(a) + \phi(b)$ and $\phi(ab) = \phi(a) \bullet \phi(b)$ for all $a, b \in R$.*

**DEFINITION 1.3.6**: *Let R be a ring. I a non-empty subset of R is called right (left) ideal of R if*

1. *I is a subring.*
2. *For $r \in R$ and $i \in I$, $ir \in I$ (or $ri \in I$).*

*If I is simultaneously both a right and a left ideal of R we say I is an ideal of R. Thus ideals are subrings but all subrings are not ideals.*

*Example 1.3.14*: Let Z be the ring of integers. $pZ = \{0, \pm p, \pm 2p, \ldots\}$ for any $p \in Z$ is an ideal of Z.

*Example 1.3.15*: Let $Z_{12} = \{0, 1, 2, \ldots, 11\}$ be the ring of integers modulo 12. $I = \{0, 6\}$ is an ideal of $Z_{12}$, $P = \{0, 3, 6, 9\}$ is also an ideal of $Z_{12}$ $I_2 = \{0, 2, 4, 6, 8, 10\}$ is an ideal of $Z_{12}$.

*Example 1.3.16*: Let R[x] be a polynomial ring. $p(x) = p_0 + p_1 x + \ldots + p_n x^n$ be a polynomial of degree n ($p_n \neq 0$). Clearly p(x) generates an ideal. We leave it for the reader to check this fact. We denote the ideal generated by p(x) by $\langle p(x) \rangle$.

**DEFINITION 1.3.7**: *Let $\phi: R \to R'$ be a ring homomorphism the kernel of $\phi$ denoted by ker $\phi = \{x \in R \,/\, \phi(x) = 0\}$ is an ideal of R.*

**DEFINITION 1.3.8**: *Let R be any ring, I an ideal of R. The set $R/I = \{a + I \,/\, a \in R\}$ is defined as the quotient ring. For this quotient ring, I serves as the additive identity.*

The reader is requested to prove R / I is a ring.

**PROBLEMS**:

1. Let F be a field. Prove F has no ideals.

2. Find all ideals of the ring $Z_{24}$.

3. Prove $Z_{29}$ has no ideals.

4. Let $M_{2\times 2} = \left\{ \begin{pmatrix} a & b \\ c & d \end{pmatrix} \Big/ a, b, c, d \in Z_3 = \{0, 1, 2\} \right\}$, $M_{2\times 2}$ is a ring under usual matrix addition and matrix multiplication.
   - i. Find one right ideal of $M_{2\times 2}$.
   - ii. Find one left ideal of $M_{2\times 2}$.
   - iii. Find an ideal of $M_{2\times 2}$.

5. In problem 4 find a subring of $M_{2\times 2}$, which is not an ideal of $M_{2\times 2}$.

6. Let $Z_7[x]$ be the polynomial ring. Suppose $p(x) = x^2 + 3$. Find the ideal I generated by p(x).

7. Let $Z_{12} = \{0, 1, 2, \ldots, 11\}$ be the ring of modulo integers 12. Let the ideal $I = \{0, 2, 4, 6, 8, 10\}$. Find the quotient ring $Z_{12} / I$. Is $Z_{12} / I$ a field?

8. Let $Z_7[x]$ be the polynomial ring over $Z_7$. $I = \langle x^3 + 1 \rangle$ be the ideal generated by the polynomial $p(x) = x^3 + 1$. Find $Z_7[x] / \langle x^3 + 1 \rangle$. When will $Z_7[x] / \langle x^3 + 1 \rangle$ be a field?

9. Find $\dfrac{Z[x]}{\langle x^2 + 1 \rangle}$.

10. Find all principal ideals in $Z_3^5[x]$ = {all polynomials of degree ≤ less than or equal to 5}. (Hint: We say any ideal is principal if it is generated by a single element).

11. Construct a prime field with 53 elements.

12. Prove $Z_2[x] / \langle x^2 + x +1 \rangle$ is a non-prime field with 4 elements in it.

13. Let Z be the ring of integers, prove nZ for some positive integer n is a principal ideal of Z.

14. Can Z have ideals, which are not principal?

15. Can $Z_n$ (n any positive integer) have ideals, which are not principal ideals of $Z_n$?

16. Let $Z_{24}$ = {0, 1, 2, … , 23} be the ring of integers modulo 24. Find an ideal I in $Z_{24}$ so that the quotient ring $Z_{24}$ / I has the least number of elements in it.

## 1.4 Vector spaces

In this section we introduce the concept of vector spaces mainly to compare and contrast with semivector spaces built over semifields. We just recall the most important definitions and properties about vector spaces.

**DEFINITION 1.4.1**: *A vector space (or a linear space) consists of the following*

1. *a field F of scalars.*
2. *a set V of objects called vectors.*
3. *a rule (or operation) called vector addition, which associates with each pair of vectors $\alpha, \beta$ in V a vector $\alpha + \beta$ in V in such a way that*

   i. *addition is commutative, $\alpha + \beta = \beta + \alpha$.*
   ii. *addition is associative; $\alpha + (\beta + \gamma) = (\alpha + \beta) + \gamma$.*
   iii. *there is a unique vector 0 in V, called the zero vector, such that $\alpha + 0 = \alpha$ for all $\alpha \in V$.*
   iv. *for each vector $\alpha$ in V there is a unique vector $-\alpha$ in V such that $\alpha + (-\alpha) = 0$.*
   v. *a rule (or operation) called scalar multiplication which associates with each scalar c in F and a vector $\alpha$ in V a vector $c\alpha$ in V called the product of c and $\alpha$ in V such that*

      a. *$1.\alpha = \alpha$ for every $\alpha \in V$.*
      b. *$(c_1, c_2) \alpha = c_1 (c_2 \alpha)$.*
      c. *$c (\alpha + \beta) = c\alpha + c\beta$.*
      d. *$(c_1 + c_2)\alpha = c_1\alpha + c_2\alpha$ for $c_1, c_2, c \in F$ and $\alpha, \beta \in V$.*

It is important to state that vector space is a composite object consisting of a field F, a set of ' vectors' and two operations with certain special properties. The same set of vectors may be part of a number of distinct vector spaces. When there is no chance of confusion, we may simply refer to the vector space as V. We shall say 'V is a vector space over the field F'.

*Example 1.4.1*: Let R[x] be the polynomial ring where R is the field of reals. R[x] is a vector space over R.

*Example 1.4.2*: Let Q be the field of rationals and R the field of reals. R is a vector space over Q.

It is important and interesting to note that Q is not a vector space over R in the example 1.4.2.

*Example 1.4.3*: Let F be any field $V = F \times F = \{(a, b) / a, b \in F\}$. It is left for the reader to verify V is a vector space over F.

*Example 1.4.4*: Let $V = \{M_{n \times m}\} = \{(a_{ij}) / a_{ij} \in Q\}$. V is the set of all $n \times m$ matrices with entries from Q. It is easily verified that V is a vector space over Q.

**DEFINITION 1.4.2**: *Let V be a vector space over the field F. Let $\beta$ be a vector in V, $\beta$ is said to be a linear combination of vectors $\alpha_1, \ldots, \alpha_n$ in V provided there exists scalars $c_1, c_2, \ldots, c_n$ in F such that $\beta = c_1\alpha_1 + \ldots + c_n\alpha_n = \sum_{i=1}^{n} c_i \alpha_i$.*

**DEFINITION 1.4.3**: *Let V be a vector space over the field F. A subspace of V is a subset W of V which is itself a vector space over F with the operations of vector addition and scalar multiplication on V.*

**DEFINITION 1.4.4**: *Let S be a set of vectors in a vector space V. The subspace spanned by S is defined to be the intersection W of all subspaces of V which contain S. When S is a finite set of vectors say $S = \{\alpha_1, \ldots, \alpha_n\}$ we shall simply call W the subspace spanned by the vectors $\alpha_1, \alpha_2, \ldots, \alpha_n$.*

**DEFINITION 1.4.5**: *Let V be a vector space over F. A subset S of V is said to be linearly dependent (or simply dependent) if there exist distinct vectors $\alpha_1, \alpha_2, \ldots, \alpha_n$ in S and scalars $c_1, \ldots, c_n$ in F not all of which are 0 such that $\alpha_1 c_1 + \ldots + \alpha_n c_n = 0$.*

*A set that is not linearly dependent is called linearly independent. If the set S contains only finitely many vectors $\alpha_1, \alpha_2, \ldots, \alpha_n$ we sometimes say that $\alpha_1, \ldots, \alpha_n$ are dependent (or independent) instead of saying S is dependent (or independent).*

**DEFINITION 1.4.6**: *Let V be a vector space over the field F. A basis for V is a linearly independent set of vectors in V, which spans the space V. The space V is finite dimensional if it has a finite basis which spans V, otherwise we say V is infinite dimensional.*

***Example 1.4.5***: Let $V = F \times F \times F = \{(x_1, x_2, x_3) / x_1, x_2, x_3 \in F\}$ where F is a field. V is a vector space over F. The set $\beta = \{(1, 0, 0), (0, 1, 0), (0, 0, 1)\}$ is a basis for V.

It is left for the reader to verify that $\beta$ spans $V = F \times F \times F = F^3$; we can say $F^3$ is a vector space over F of dimension three.

***Example 1.4.6***: Let F be a field and $F^n = F \times \ldots \times F$ (n times), $F^n$ is a vector space over F. A set of basis for $F^n$ over F is $\beta = \{(1, 0, 0, \ldots, 0), (0, 1, 0, \ldots, 0), (0, 0, 1, 0, \ldots, 0), \ldots, (0, 0, \ldots, 0, 1)\}$. It can be shown, $F^n$ is spanned by $\beta$ the dimension of $F^n$ is n. We call this particular basis as the standard basis of $F^n$.

***Example 1.4.7***: Let $F_n[x]$ be a vector space over the field F; where $F_n[x]$ contains all polynomials of degree less or equal to n. Now $\beta = \{1, x, x^2, \ldots, x^n\}$ is a basis of $F_n[x]$. The dimension of $F_n[x]$ is $n + 1$.

***Example 1.4.8***: Let $F[x]$ be the polynomial ring which is a vector space over F. Now the set $\{1, x, x^2, \ldots, x^n, \ldots\}$ is a basis of $F[x]$. The dimension of the vector space $F[x]$ over F is infinite.

***Remark***: A vector space V over F can have many basis but for that vector space the number of elements in each of the basis is the same; which is the dimension of V.

**DEFINITION 1.4.7**: *Let V and W be two vector spaces defined over the same field F. A linear transformation T: V $\rightarrow$ W is a function from V to W such that $T(c\alpha + \beta) = cT(\alpha) + T(\beta)$ for all $\alpha, \beta \in V$ and for all scalars $c \in F$.*

***Remark***: The linear transformation leaves the scalars of the field invariant. Linear transformation cannot be defined if we take vector spaces over different fields. If W = V then the linear transformation from V to V is called the linear operator.

**DEFINITION 1.4.8**: *Let L (V, W) denote the collection of all linear transformation of the vector space V to W, V and W vector spaces defined over the field F. L(V, W) is a vector space over F.*

***Example 1.4.9***: Let $R^3$ be a vector space defined over the reals R. $T(x_1, x_2, x_3) = (3x_1, x_1, -x_2, 2x_1 + x_2 + x_3)$ is a linear operator on $R^3$.

***Example 1.4.10***: Let $R^3$ and $R^2$ be vector spaces defined over R. T is a linear transformation from $R^3$ into $R^2$ given by $T(x_1, x_2, x_3) = (x_1 + x_2, 2x_3 - x_1)$. It is left for the reader to verify T is a linear transformation.

***Example 1.4.11***: Let $V = F \times F \times F$ be a vector space over F. Check whether the 3 sets are 3 distinct sets of basis for V.

1. $\{(1, 5, 0), (0, 7, 1), (3, 8, 8)\}$.
2. $\{(4, 2, 0), (2, 0, 4), (0, 4, 2)\}$.
3. $\{(-7, 2, 1), (0, -3, 5), (7, 0, -1)\}$.

**PROBLEMS**:

1.  Let $M_{3\times 5} = \{(a_{ij}) \mid a_{ij} \in Q\}$ denote the set of all 3×5 matrices with entries from Q the rational field.

    i.   Prove $M_{3\times 5}$ is a vector space over Q the rationals.
    ii.  Find a basis of $M_{3\times 5}$.
    iii. What is the dimension of $M_{3\times 5}$?

2.  Prove L(V, W) is a vector space over F if V and W are vector spaces over F.

3.  Let V be a vector space of dimension 3 over a field F and W be a vector space of dimension 5 defined over F. Find the dimension of L(V, W) over F.

4.  Suppose V is a vector space of dimension n over F. If $B = \{v_1, \dots, v_n\}$ and $B' = \{w_1, \dots, w_n\}$ are two distinct basis of V. Find a method by which one basis can be represented in terms of the other (The Change of Basis rule).

5.  Show the spaces $M_{n\times n} = \{(a_{ij}) \mid a_{ij} \in Q\}$ the set of n×n matrices with entries from Q is isomorphic with L(V, V) = {set of all linear operators from V to V}. (V is a n-dimensional vector space over Q).

6.  Prove we can always get a matrix associated with any linear operator from a finite dimensional vector space V to V.

7.  Let T be the linear operator on $R^4$.
    $T(x_1, x_2, x_3, x_4) = (x_1 + 3x_3 - x_4, x_3 + 3x_4 - x_2, 5x_2 - x_4, x_1 + x_2 + x_3 + x_4)$.

    i.   What is the matrix of T in the standard basis for $R^4$?
    ii.  What is the matrix of T in the basis $\{\alpha_1, \alpha_2, \alpha_3, \alpha_4\}$ where
        $\alpha_1 = (1, 1, 1, 0)$
        $\alpha_2 = (0, 0, 3, 4)$
        $\alpha_3 = (0, 5, 0, 2)$
        $\alpha_4 = (1, 0, 0, 1)$.

## 1.5 Group rings and semigroup rings

In this section we introduce the notion of group rings and semigroup rings; the main motivation for introducing these concepts is that in this book we will define analogously group semirings and semigroup semirings where the rings are replaced by semirings. Several new properties not existing is the case of group rings is found in the case of group semirings. Throughout this section by the ring R we mean either R is a field or R is a commutative ring with 1. G can be any group but we assume the operation on the group G is only multiplication. S is a semigroup under multiplication.

**DEFINITION 1.5.1**: *Let R be a ring and G a group the group ring RG of the group G over the ring R consists of all finite formal sums of the form $\sum_{i} \alpha_i g_i$ (i runs over finite number) where $\alpha_i \in R$ and $g_i \in G$ satisfying the following conditions:*

i. $\sum_{i=1}^{n} \alpha_i g_i = \sum_{i=1}^{n} \beta_i g_i \Leftrightarrow \alpha_i = \beta_i$ *for i = 1, 2, ... ,n.*

ii. $\left(\sum_{i=1}^{n} \alpha_i g_i\right) + \left(\sum_{i=1}^{n} \beta_i g_i\right) = \sum_{i=1}^{n} (\alpha_i + \beta_i) g_i$ .

iii. $\left(\sum_{i=1}^{n} \alpha_i g_i\right)\left(\sum_{i=1}^{n} \beta_i h_i\right) = \sum_{i=1}^{n} \gamma_k m_k$ *where $g_i h_j = m_k$ and $\gamma_k = \sum \alpha_i \beta_j$*

iv. $r_i g_i = g_i r_i$ *for all $g_i \in G$ and $r_i \in R$.*

v. $r \sum_{i=1}^{n} r_i g_i = \sum_{i=1}^{n} (r r_i) g_i$ *for $r \in R$ and $\sum r_i g_i \in RG$.*

*RG is an associative ring with $0 \in R$ as its additive identity. Since $I \in R$ we have $G = 1 \bullet G \subseteq RG$ and $R \bullet e = R \subseteq RG$, where e is the identity of the group G. If we replace the group G by the semigroup S with identity we get the semigroup ring RS in the place of the group ring RG.*

***Example 1.5.1***: Q be the field of rationals and $G = \langle g\,/\,g^2 = 1\rangle$ be the cyclic group. The group ring QG = {a + bg | a, b ∈ Q and g ∈ G} is a commutative ring of characteristic 0.

***Example 1.5.2***: Let $Z_8$ = {0, 1, 2, ... , 7} be the ring of integers modulo 8. $S_3$ be the symmetric group of degree 3. $Z_8 S_3$ is the group ring of $S_3$ over $Z_8$. $Z_8 S_3$ is a non-commutative ring of characteristic 8.

***Example 1.5.3***: Let R be the real field, $S_n$ the symmetric group of degree n. The group ring $RS_n$ is a non-commutative ring of characteristic 0. This is not a skew field for $RS_n$ has zero divisors.

***Example 1.5.4***: Let $Z_5$ = {0, 1, 2, 3, 4} be the prime field of characteristic 5. $G = \langle g\,|\,g^{12} = 1\rangle$ be the cyclic group of order 12. The group ring $Z_5 G$ is a commutative ring with characteristic 5 and has zero divisors.

**PROBLEMS**:

1. Let Q be the field of rationals. $G = S_8$ be the symmetric group of degree 8. Find in the group ring QG = QS, a right ideal and an ideal.

2. S(6) be the symmetric semigroup. Let $Z_6$ = {0, 1, 2, ... , 5} be the ring of integers modulo 6. Find in the semigroup ring $Z_6 S(6)$.

    i. Ideals.
    ii. Right ideals.
    iii. Zero divisors.

iv.   Subrings which are not ideals.

3. Let $Z_5S_3$ be the group ring of the group $S_3$ over the prime field $Z_5$. $Z_4S_7$ the group ring of the symmetric group $S_7$ over $Z_4$.

   i.    Construct a ring homomorphism $\phi$ from $Z_5S_3$ to $Z_4S_7$.
   ii.   Find ker $\phi$.
   iii.  Find the quotient ring $Z_5S_3$ / ker $\phi$.

4. Let $Z_2S_3$ and $Z_2S(3)$ be the group ring and the semigroup ring. Can we construct a homomorphism from $Z_2S_3$ to $Z_2S(3)$?

5. Find all zero divisors in $Z_3S(3)$, the semigroup ring of the semigroup $S(3)$ over the prime field $Z_3$.

6. Find all zero divisors in $Z_3S_3$, the group ring of the group $S_3$ over the ring $Z_3$.

7. $Z_3S(3)$ or $Z_3S_3$ which has more number of zero divisors? (Hint: We know $S_3 \subset S(3)$ use this to prove the result).

## Supplementary Reading

# CHAPTER TWO
# SEMIRINGS AND ITS PROPERTIES

The study of the concept of semiring is very meagre. In my opinion I have not come across a textbook that covers completely all the properties of semirings. Hence this complete chapter is devoted to introduction of semirings, polynomial semirings and many new properties about it, analogous to rings. This chapter also gives for the sake of completeness the definition of several types of special semirings like ∗-semirings, congruence simple semirings and so on. We do not intend to give all definition or all properties instead we expect the reader to refer those papers which are enlisted in the supplementary reading at the end of this chapter. This chapter has six sections. In sections one and two we define semirings and give examples and prove some basic properties. Section three shows how lattices are used to construct semirings. Polynomial semirings are introduced in section four. Section five defines and recalls the definitions of group semirings and semigroup semirings. The final section mainly recalls some of the special types of semirings like c-semirings, ∗-semirings, inductive ∗-semirings etc.

## 2. 1 Definition and examples of semirings

In this section we introduce the concept of semirings and give some examples. This is mainly carried out because we do not have many textbook for semirings except in the book 'Handbook of Algebra' Vol. I, by Udo, which carries a section on semirings and semifields.

**DEFINITION (LOUIS DALE)**: *Let S be a non-empty set on which is defined two binary operations addition '+' and multiplication '•' satisfying the following conditions:*

1. *(S, +) is a commutative monoid.*
2. *(S, •) is a semigroup.*
3. *(a + b) • c = a • c + b • c and a • (b + c) = a • b + a • c for all a, b, c in S.*

*That is multiplication '•' distributes over the operation addition '+'. (S, +, •) is a semiring.*

**DEFINITION (LOUIS DALE)**: *The semiring (S, +, •) is said to be a commutative semiring if the semigroup (S, •) is a commutative semigroup. If (S, •) is not a commutative semigroup we say S is a non-commutative semiring.*

**DEFINITION (LOUIS DALE)**: *If in the semiring (S, +, •), (S, •) is a monoid that is there exists $1 \in S$ such that $a \cdot 1 = 1 \cdot a = a$ for all $a \in S$. We say the semiring is a semiring with unit.*

Throughout this book $Z^+$ will denote the set of all positive integers and $Z^o = Z^+ \cup \{0\}$ will denote the set of all positive integers with zero. Similarly $Q^o = Q^+ \cup \{0\}$ will denote the set of all positive rationals with zero and $R^o = R^+ \cup \{0\}$ denotes the set of all positive reals with zero.

**DEFINITION (LOUIS DALE)**: *Let (S, +, •) be a semiring. We say the semiring is of characteristic m if ms = s + … + s (m times) equal to zero for all s ∈ S. If no such m exists we say the characteristic of the semiring S is 0 and denote it as characteristic S = 0. In case S has characteristic m then we say characteristic S = m.*

Here it is interesting to note that certain semirings have no characteristic associated with it. This is the main deviation from the nature of rings. We say the semiring S is finite if the number of elements in S is finite and is denoted by |S| or o(S). If the number of elements in S is not finite we say S is of infinite cardinality. Now we give some examples of semirings.

*Example 2.1.1*: Let $Z^o = Z^+ \cup \{0\}$. $(Z^o, +, •)$ is a semiring of infinite cardinality and the characteristic $Z^o$ is 0. Further $Z^o$ is a commutative semiring with unit.

*Example 2.1.2*: Let $Q^o = Q^+ \cup \{0\}$. $(Q^o, +, o)$ is also a commutative semiring with unit of infinite cardinality and characteristic $Q^o$ is 0.

*Example 2.1.3*: Let

$$M_{2\times 2} = \left\{ \begin{pmatrix} a & b \\ c & d \end{pmatrix} / a, b, c, d \in Z^o \right\}$$

= set of all 2 × 2 matrices with entries from $Z^o$. Clearly $(M_{2\times 2}, +, o)$ is a semiring under matrix addition and matrix multiplication. $M_{2\times 2}$ is a non-commutative semiring of characteristic zero with unit element $\begin{pmatrix} 1 & 0 \\ 0 & 1 \end{pmatrix}$ and is of infinite cardinality.

*Example 2.1.4*: Let S be the chain lattice given by the following Hasse diagram:

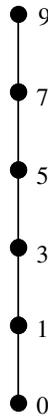

**Figure 2.1.1**

S is a semiring with inf and sup as binary operations on it. This is a commutative semiring of finite cardinality or order. This semiring has no characteristic associated with it.

*Example 2.1.5*: Consider the following lattice given by the Hasse diagram:

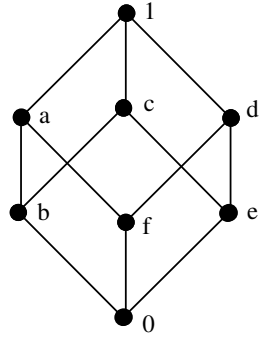

**Figure 2.1.2**

It can be verified that this lattice is also a semiring which is commutative with finite cardinality and has no characteristic associated with it.

**DEFINITION 2.1.1**: *Let $S_1$ and $S_2$ be two semirings. The direct product of $S_1 \times S_2 = \{(s_1, s_2) / s_1 \in S_1$ and $s_2 \in S_2\}$ is also a semiring with component-wise operation. Similarly if $S_1, S_2, …, S_n$ be n semirings. The direct product of these semirings denoted by $S_1 \times S_2 \times S_3 \times … \times S_n = \{(s_1, s_2, …, s_n) / s_i \in S_1; i = 1, 2, …, n\}$ is a semiring also known as the direct product of semirings.*

*Example 2.1.6*: Let $Z^o$ be the semiring $Z^o \times Z^o \times Z^o = \{(a, b, c) / a, b, c \in Z^o\}$ is a semiring. This enjoys vividly different properties from $Z^o$.

*Example 2.1.7*: Let $S_1$ be the two-element chain lattice and $S_2$ be the lattice given by the following Hasse diagram:

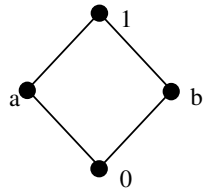

**Figure 2.1.3**

$S_1 \times S_2 = \{(0, 0), (0, a), (0, b), (0, 1), (1, 0), (1, a), (1, b), (1, 1)\}$ is a semiring under the operations of the lattices, $S_1 \times S_2$ has the following Hasse diagram:

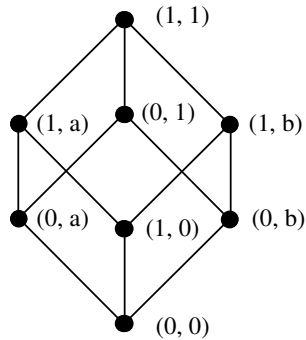

**Figure 2.1.4**

$S = S_1 \times S_2$ is a semiring (left for the reader as an exercise to verify).

*Example 2.1.8*: Let S be the lattice with the following Hasse diagram

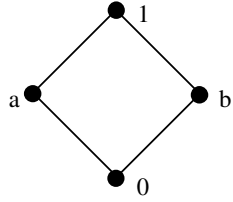

**Figure 2.1.5**

S is a commutative semiring with unit and is of finite cardinality. Let

$$S_{2 \times 2} = \left\{ \begin{pmatrix} x_1 & x_2 \\ x_3 & x_4 \end{pmatrix} \Big/ x_1, x_2, x_3, x_4 \in S \right\}$$

= set of all $2 \times 2$ matrices with entries from the semiring $S = \{0, 1, a, b\}$. Let $A, B \in S_{2 \times 2}$, where $A = \begin{pmatrix} a_1 & a_2 \\ a_3 & a_4 \end{pmatrix}$ and $B = \begin{pmatrix} b_1 & b_2 \\ b_3 & b_4 \end{pmatrix}$. Define '+' on $S_{2 \times 2}$ as $A + B = \begin{pmatrix} a_1 & a_2 \\ a_3 & a_4 \end{pmatrix} + \begin{pmatrix} b_1 & b_2 \\ b_3 & b_4 \end{pmatrix} = \begin{pmatrix} a_1 \cup b_1 & a_2 \cup b_2 \\ a_3 \cup b_3 & a_4 \cup b_4 \end{pmatrix}$. Clearly $(S_{2 \times 2}, '+')$ is a commutative monoid. The matrix $\begin{pmatrix} 0 & 0 \\ 0 & 0 \end{pmatrix}$ acts as the additive identity. For $A, B \in S_{2 \times 2}$ define • as $A \bullet B =$

$\begin{pmatrix} a_1 & a_2 \\ a_3 & a_4 \end{pmatrix} \bullet \begin{pmatrix} b_1 & b_2 \\ b_3 & b_4 \end{pmatrix} = \begin{pmatrix} (a_1 \cap b_1) \cup (a_2 \cap b_3) & (a_1 \cap b_2) \cup (a_2 \cap b_4) \\ (a_3 \cap b_1) \cup (a_4 \cap b_3) & (a_3 \cap b_2) \cup (a_4 \cap b_4) \end{pmatrix}$.

Clearly $(S_{2 \times 2}, \bullet)$ is a semigroup. It is left for the reader to verify $(S_{2 \times 2}, +, \bullet)$ is a semiring. This semiring is of finite cardinality but is non-commutative for if $A = \begin{pmatrix} a & b \\ 0 & 1 \end{pmatrix}$ and $B = \begin{pmatrix} 1 & a \\ b & b \end{pmatrix}$. $A \bullet B = \begin{pmatrix} 1 & 1 \\ b & b \end{pmatrix}$. Now $B \bullet A = \begin{pmatrix} a & 1 \\ 0 & b \end{pmatrix}$. So $A \bullet B \neq B \bullet A$ for $A, B \in S_{2 \times 2}$.

Thus we have seen semirings of characteristic 0 which is both commutative and non-commutative having infinite cardinality. We are yet to find semiring of characteristic n, n ≠ 0 and we have seen finite cardinality semirings with no characteristic associated with it both commutative and non-commutative.

**PROBLEMS**:

1. Give an example of a semiring of order 27.

2. How many elements does the semiring $M_{3\times 3} = \{(a_{ij})/ a_{ij} \in L$ where the lattice L is given by the following Hasse diagram}?

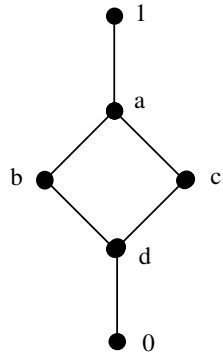

**Figure 2.1.6**

3. Prove $M_{3\times 3}$ given in Problem 2 is a non-commutative semiring.

4. Find the unit element of $M_{3\times 3}$.

5. Let $X = \{x_1, x_2, x_3, x_4, x_5, x_6\}$ and P(X) be the power set of X. Prove P(X) is a semiring.

6. Can non-commutative semirings of finite order n exist for every integer n?

7. Can any finite semiring have characteristic p? (p any finite integer)

8. Give an example of a finite semiring which has no characteristic associated with it.

9. Give an example of an infinite non-commutative semiring of characteristic zero.

## 2.2 Semirings and its properties

In this section we introduce properties like subsemirings, ideals in semirings, zero divisors, idempotents and units in semirings.

**DEFINITION 2.2.1**: *Let S be a semiring. P a subset of S. P is said to be a subsemiring of S if P itself is a semiring.*

*Example 2.2.1*: Let $Z^o$ be the semiring. $2Z^o = \{0, 2, 4, \ldots\}$ is a subsemiring of $Z^o$.

*Example 2.2.2*: Let $Z^o[x]$ be the polynomial semiring. $Z^o \subseteq Z^o[x]$ is a subsemiring of $Z^o[x]$.

*Example 2.2.3*: Let S be a distributive lattice which is a semiring given by the following Hasse diagram:

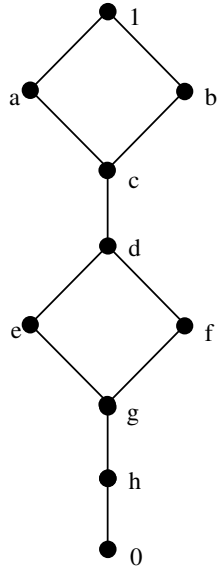

**Figure 2.2.1**

L = {1, a, c, d, e, g, h, 0} is a subsemiring of S.

**DEFINITION 2.2.2**: *Let S be a semiring. I be a non-empty subset of S. I is a right (left) ideal of S if*

1. *I is a subsemiring.*
2. *For all i ∈ I and s ∈ S we have is ∈ I (si ∈ I).*

**DEFINITION 2.2.3**: *Let S be a semiring. A non-empty subset I of S is said to be an ideal of S if I is simultaneously a right and left ideal of S.*

*Example 2.2.4*: Let $Z^o$ be the semiring. $nZ^o$ for any integer. n is an ideal of $Z^o$.

*Example 2.2.5*: Let S be a semiring given by the following lattice whose Hasse diagram :

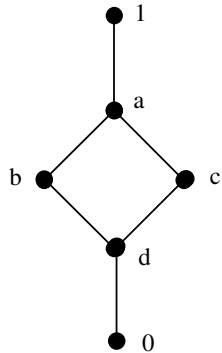

**Figure 2.2.2**

Clearly I = {0, d, c, a, b} is an ideal of S.

**DEFINITION 2.2.4**: *Let S be a semiring. $x \in S \setminus \{0\}$ is said to be a zero divisor of S if there exists $y \neq 0$ in S such that $x \bullet y = 0$.*

**DEFINITION 2.2.5**: *Let S be a non-commutative semiring with unit. $x \in S$ is said to be have right (left) unit if their exists a $y \in S$ such that $xy = 1$ (or $yx = 1$).*

**DEFINITION 2.2.6**: *Let S be a semiring. $x \in S$ is an idempotent, if $x \bullet x = x^2 = x$.*

*Example 2.2.6*: Let S be a semiring given by the following Hasse diagram a, b ∈ S is such that $a \bullet b = 0$ i.e. a is a zero divisor and $a^2 = a$ and $b^2 = b$ so the semiring has idempotents.

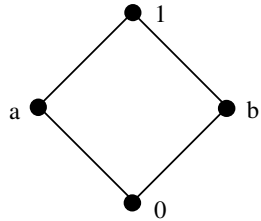

**Figure 2.2.3**

*Example 2.2.7*: Let S be a semiring given by $S = Z^o \times Z^o \times Z^o$. S has zero divisors given by a = (8 0, 2) and b = (0, 6, 0) and $a \bullet b = (0, 0, 0)$.

*Example 2.2.8*: Let $S = Q^o \times Q^o \times Q^o \times Q^o$ be the semiring under component wise addition and multiplication. For S is a semiring with unit (1, 1, 1, 1). Let a = (3, ¼, 5/3, 7/2) ∈ S. The inverse of a is b = (1/3, 4, 3/5, 2/7) ∈ S is such that ab = (1, 1, 1, 1). Thus S has units.

**DEFINITION 2.2.7**: *Let S and S' be two semirings. A mapping $\phi : S \to S'$ is called the semiring homomorphism if $\phi(a + b) = \phi(a) + \phi(b)$ and $\phi(a \bullet b) = \phi(a) \bullet \phi(b)$ for all a, b ∈ S.*

*If $\phi$ is a one to one and onto map we say $\phi$ to be a semiring isomorphism.*

**DEFINITION (LOUIS DALE)**: *Let S be a semiring. We say S is a strict semiring if $a + b = 0$ implies $a = 0$ and $b = 0$.*

Chris Monico calls this concept as zero sum free.

*Example 2.2.9*: $Z^o$ the semiring is a strict semiring.

*Example 2.2.10*: The semiring $Q^o$ is a strict semiring.

*Example 2.2.11*: The chain lattice L given by the following Hasse diagram is a strict semiring.

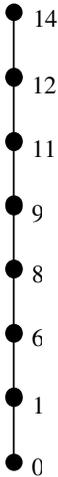

**Figure 2.2.4**

**DEFINITION 2.2.8**: *Let S be a semiring with unit 1. We say an element x is invertible or has an inverse if there exists a y ∈ S such that xy = yx = 1.*

*Example 2.2.12*: Let $Q^o$ be the semiring of rationals. Every element in $Q^o \setminus \{0\}$ is a unit in $Q^o$. All elements are invertible.

*Example 2.2.13*: Let $M_{2\times 2} = \left\{ \begin{pmatrix} a & b \\ c & d \end{pmatrix} \bigg/ a,b,c,d \in Q^o \right\}$ be the semiring under matrix addition and matrix multiplication. $M_{2\times 2}$ is a non-commutative semiring with unit and has non-trivial units in it.

**PROBLEMS**:

1. Find 3 subsemirings of the semiring S where S is a distributive lattice given by the following Hasse diagram.

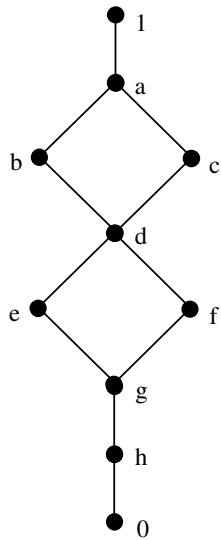

**Figure 2.2.5**

2. How many subsemirings does $Z^o$ contain?

3. Is every subsemiring in $Z^o$ an ideal of $Z^o$?

4. Find all ideals of the semiring given in Problem 1.

5. Find at least two right ideals of the semiring $M_{5\times 5} = \{(a_{ij})/ a_{ij} \in Z^o\}$.

6. Can the semiring $Q^o$ have ideals?

7. Can the semiring $M_{2\times 2} = \left\{ \begin{pmatrix} a & b \\ c & d \end{pmatrix} / a,b,c,d \in R^o \right\}$ have ideals? Justify your answer.

8. Let $L_{2\times 2} = \left\{ \begin{pmatrix} a & b \\ c & d \end{pmatrix} / a,b,c,d \in L \right\}$ (the set of all $2 \times 2$ matrices with entries from the lattice L given by the following Hasse diagram).

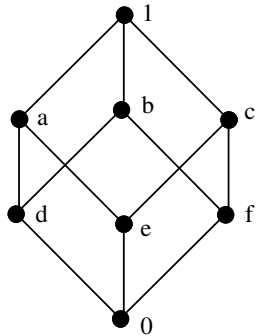

**Figure 2.2.6**

  a. Find the right ideals of $L_{2\times 2}$.
  b. Find the ideals of $L_{2\times 2}$.
  c. Find the number of elements in $L_{2\times 2}$.

9. Find non-trivial idempotents in the semiring $L_{2\times 2}$ given in problem 8.

10. Can we find non-commutative semirings of order n for every positive integer n ≥ 2? (Hint: We have commutative semirings of order n for every positive integer n, n ≥ 2 as we have chain lattices of all orders.)

11. Can $Z^o$ have non-trivial idempotents? zerodivisors? units?

12. Is $R^o$ a strict semiring?

13. Give an example of a semiring which is not a strict semiring.

14. Let $M_{3\times3} = \{(a_{ij})/ a_{ij} \in Z^o\}$ = {set of all $3 \times 3$ matrices with entries from $Z^o$} be the semiring under matrix addition and matrix multiplication. Can $M_{3\times3}$ have units? zero divisors? idempotents?

15. Find a semiring homomorphism between $C_6$ and the semiring $L_5$ where the Hasse diagram for $C_6$ and $L_5$ are given below:

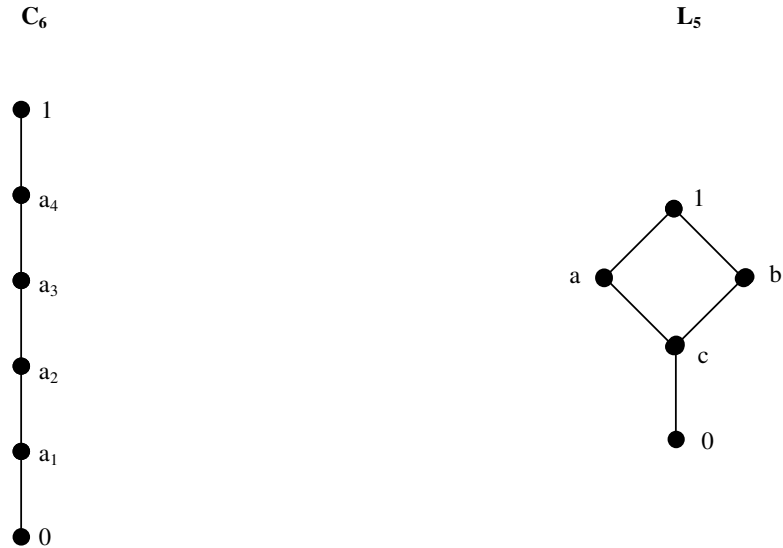

**Figure 2.2.7**

## 2.3 Semirings using distributive lattices

In this section we prove all distributive lattices are semirings and obtain some special and unique properties enjoyed by these class of semirings. Further this is the only class of semirings known to us which is of finite order and has no characteristic.

Throughout this book by $C_n$ we will denote a chain lattice with n elements given by the following Hasse diagram:

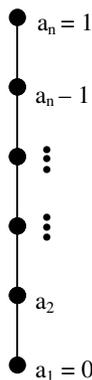

**Figure 2.3.1**

That is we assume $C_n$ is a chain lattive which has 0 as the least element and 1 as the greatest element.

*Example 2.3.1*: Let $C_5$ be the chain lattice given by the following Hasse diagram

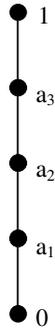

**Figure 2.3.2**

It is a semigroup under union '$\cup$' or max and intersection '$\cap$' or min as binary operations. Since $C_5$ is a distributive lattice $C_5$ is a semiring.

**THEOREM 2.3.1**: *Let $(L, \cup, \cap)$ be a distributive lattice. L is a semiring.*

*Proof*: L is a semiring as under '$\cup$', L is a commutative semigroup with '0' as identity and L under '$\cap$' is a commutative semigroup; and $\cup$ and $\cap$ distribute over each other. Hence $(L, \cup, \cap)$ is a semiring.

**COROLLARY**: All chain lattices are semirings.

*Remark*: All lattices are not semirings. For the lattice L given by the following Hasse diagram is not a semiring.

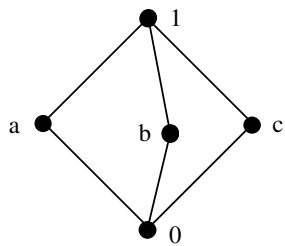

**Figure 2.3.3**

Since distributive law is not true, we see L is not a semiring.

**THEOREM 2.3.2**: *Let B be a Boolean algebra. B is a semiring*.

*Proof*: B is a distributive lattice so by theorem 2.3.1, B is a semiring.

**THEOREM 2.3.3**: *A semiring which is a distributive lattice has no characteristic associated with it.*

*Proof*: Since if L is a distributive lattice which is a semiring, we see na = a + … + a (n times) will not be zero but equal to a as a + a = a for all a ∈ L. Hence the claim.

**THEOREM 2.3.4**: *The direct product of distributive lattices is a semiring.*

*Proof*: Since direct product of distributive lattices is distributive, we see direct product of distributive lattices is a semiring.

*Example 2.3.2*: Let $S = L_1 \times L_2 \times L_3$ be the direct product of the lattices given by the following Hasse diagram. S is a semiring.

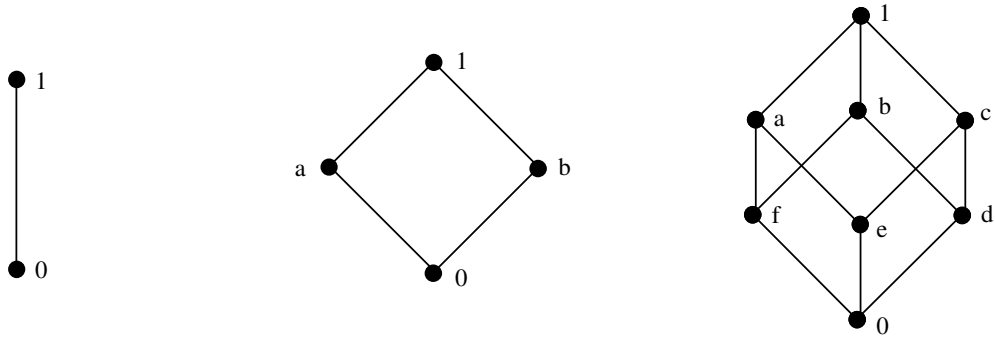

**Figure 2.3.4**

Using these distributive lattices L we can construct finite non-commutative semirings by constructing square matrices with entries from L. Let L be a finite distributive lattice with 0 and 1. Denote the elements in L by $L = \{0, x_1, x_2, …, x_n, 1\}$. $M_{n \times n} = \{(a_{ij}) / a_{ij} \in L\} = \{$set of all n × n matrices with entries from L$\}$. The two operations sum and product denoted by + and • is defined as follows:

Let $a = (a_{ij})$, $b = (b_{ij})$, $a + b = (a_{ij} \cup b_{ij})$ where $\cup$ is the '∪' of the lattice L.

$$a \bullet b = \begin{pmatrix} a_{11} & \cdots & a_{1n} \\ a_{21} & \cdots & a_{2n} \\ \vdots & & \vdots \\ a_{n1} & \cdots & a_{nn} \end{pmatrix} \bullet \begin{pmatrix} b_{11} & \cdots & b_{1n} \\ b_{21} & \cdots & b_{2n} \\ \vdots & & \vdots \\ b_{n1} & \cdots & b_{nn} \end{pmatrix} = \begin{pmatrix} c_{11} & \cdots & c_{1n} \\ c_{21} & \cdots & c_{2n} \\ \vdots & & \vdots \\ c_{n1} & \cdots & c_{nn} \end{pmatrix}$$

where $c_{11} = [(a_{11} \cap b_{11}) \cup (a_{12} \cap b_{21}) \cup … \cup (a_{1n} \cap b_{n1})]$. Similarly $c_{12}$ etc. are calculated.

Since '∪' and '∩' distribute over each other as L is a distributive lattice, we see $\{M_{n \times n}, +, \bullet\}$ is a finite non-commutative semiring. This semiring has the zero matrix to be the additive identity and $I_{n \times n}$ the matrix with the diagonal elements 1 and all other entries as zero as the identity element. Thus only by this method we are in a position to get finite non-commutative semirings which are also strict semirings but has zero divisors.

*Example 2.3.3*: Let L be a lattice given by the following Hasse diagram:

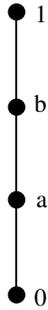

**Figure 2.3.5**

$M_{3\times 3} = \{(a_{ij})/ a_{ij} \in L\}$ be the collection of all $3 \times 3$ matrices with entries from L. Clearly $M_{3\times 3}$ is a finite non-commutative strict semiring with unit and with zero divisors and idempotents under the matrix operations described in example 2.3.2.

Now we see all Boolean algebras are semirings, for a Boolean algebra is a distrubutive complemented lattice. Here we will be using the following definitions and results about Boolean algebras.

**DEFINITION (GRATZER G.)**: *Let L be a lattice with zero. $a \in L$ is called an atom if for all $b \in L$, $0 < b \leq a \Rightarrow b = a$.*

**Result**: Let B be a finite Boolean algebra and let A denote the set of all atoms in B. Then B is isomorphic to P(A) the power set of A. i.e. $(B, \cup, \cap) \cong (P(A), \cup, \cap)$. The proof of this result is left to the reader. But we illustrate this by two examples.

*Example 2.3.4*: Let B be the Boolean algebra given by the following Hasse diagram

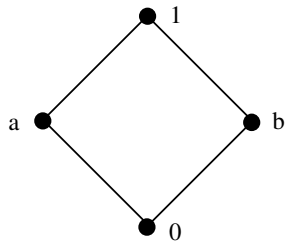

**Figure 2.3.6**

The atoms of B are a and b. So the set of atoms A = {a, b}. Now P(A) = {φ, {a}, {b}, {a,b}}. Now P(A) is a distributive complemented lattice under the operation union '∪' and intersection '∩' of subsets and complementation. (P(A), ∪, ∩) is a lattice having the following Hasse diagram:

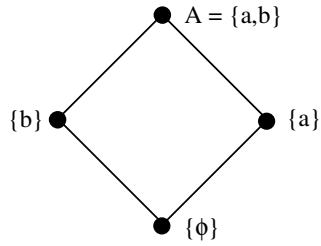

**Figure 2.3.7**

Clearly B ≅ {P(A), ∪, ∩} as lattices by the map ψ: B → {P(A), ∪, ∩} where ψ(0) = φ, ψ(a) = {a}, ψ(b) = {b} and ψ(1) = A = {a, b}.

Hence the result mentioned is true for this example. Now consider the example:

*Example 2.3.5*: Let (B, ∪, ∩) be the Boolean algebra given by the following Hasse diagram:

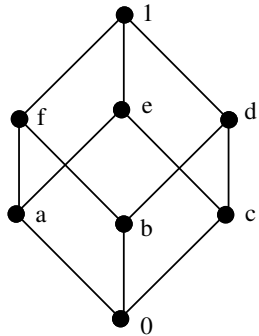

**Figure 2.3.8**

The set of atoms of B is A = {a, b, c}. P(A) = {φ, {a}, {b}, {c}, {a, b}, {a, c}, {b, c}, {a, b, c} = A}. (P(A), ∪, ∩) is a lattice having the following Hasse diagram:

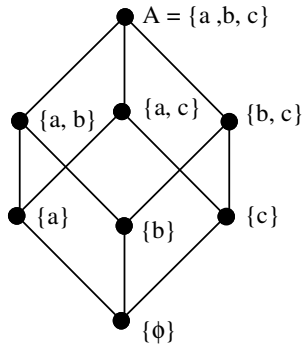

**Figure 2.3.9**

It is left for the reader to realize the lattices {P(A), ∪, ∩} and {B, ∪, ∩} are isomorphic.

Thus the semirings which are got by taking finite Boolean algebras or direct product of finite Boolean algebras are finite commutative semirings with unit and zero divisors in them.

**PROBLEMS**:

1. Find the number of elements in $M_{5\times 5} = \{(a_{ij})/ a_{ij} \in C_2\}$. (where $C_2$ is the chain lattice with 2 elements having the Hasse diagram).

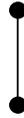

**Figure 2.3.10**

$M_{5\times 5}$ the set of all $5 \times 5$ matrices with entries from $C_2$.

   1. Find an idempotent matrix in $M_{5\times 5}$.
   2. Prove $M_{5\times 5}$ is non-commutative.
   3. Prove $M_{5x5}$ has no characteristic.
   4. Find zero divisor in $M_{5\times 5}$.
   5. Does $M_{5\times 5}$ have units?

2. Let $B_1$ and $B_2$ be two distributive lattices having the following Hasse diagram:

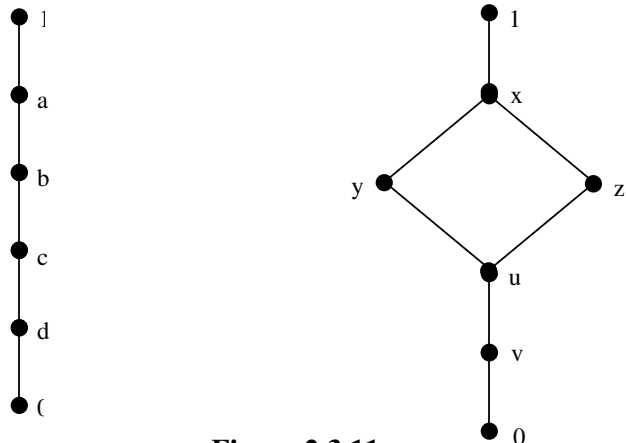

**Figure 2.3.11**

Show $B_1 \times B_2$ is a semiring and find the number of elements in $B_1 \times B_2$. Can $B_1 \times B_2$ have non-trivial idempotents? units? zero divisors? Is $B_1 \times B_2$ a commutative semiring?

3. Is the lattice given by the following Hasse diagram a semiring? Justify your answer.

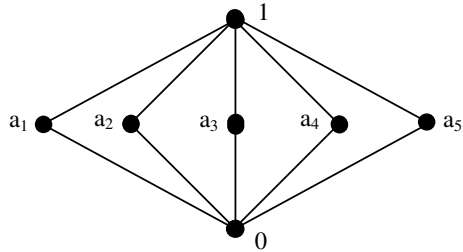

**Figure 2.3.12**

4. Find the number of elements in $M_{3\times 3} = \{(a_{ij})/a_{ij} \in L\}$ where L is the distributive lattice given by the following Hasse diagram:

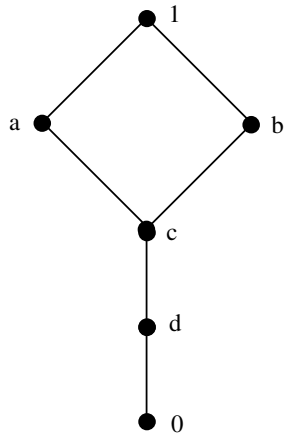

**Figure 2.3.13**

Show explicitly $M_{3\times 3}$ is non-commutative. Find a zero divisor and an idempotent in $M_{3\times 3}$.

5. Does the semiring $\{P(A), \cup, \cap\}$ where $A = \{x_1, x_2, x_3, x_4, x_6\}$ have zero divisors? Prove every element is an idempotent in $P(A)$. Is $P(A)$ a commutative semiring?

## 2.4 Polynomial Semirings

In this section we introduce polynomial semirings. Polynomial semirings have been studied by [3]. We see like polynomial rings we can discuss about the nature of the polynomial semirings over $Q^o$, $Z^o$ or $R^o$. When semirings are distributive lattices we have several nice properties about them.

**DEFINITION (LOUIS DALE):** *Let S be a commutative semiring with unit and x an indeterminate. The polynomial semiring S[x] consists of all symbols of the form $s_o +$* $s_1x + \ldots + s_nx^n$ *where n can be any non-negative integer and where the coefficients $s_0$, $s_1, \ldots, s_n$ are all in S. In order to make a semiring out of S[x] we must be able to recognize when two elements in it are equal, we must be able to add and multiply elements in S[x] so that the axioms of a semiring hold true for S[x]. If $p(x) = a_0 + a_1x + \ldots + a_mx^m$ and $q(x) = b_0 + b_1x + \ldots + b_nx^n$ are in S[x] then p(x) = q(x) if and only if for every integer $i \geq 0$; $a_i = b_i$.*

*Two polynomials are declared equal if and only if their corresponding coefficients are equal. Now we say $p(x) + q(x) = c_0 + c_ix \ldots + c_tx^t$, where for each i, $c_i = a_i + b_i$. $0 = 0 + 0x + \ldots + 0x^n$ acts as the additive identity. So S[x] is a commutative monoid under +. It is left for the reader to verify S[x] under multiplication of two polynomials is a commutative semigroup. So (S[x], +, •) is a commutative semiring with unit 1.*

**Example 2.4.1**: Let $Z^o$ be the semiring. $Z^o[x]$ is the polynomial semiring. $Z^o[x]$ is a strict polynomial semiring. $Z^o[x]$ has no zero divisors.

**Examples 2.4.2**: Let S be the distributive lattice given by the following Hasse diagram:

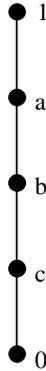

**Figure 2.4.1**

S[x] is a polynomial semiring which is commutative and has no zero divisors.

**Example 2.4.3**: Let S be the semiring given by the following distributive lattice which has the Hasse diagram as given below:

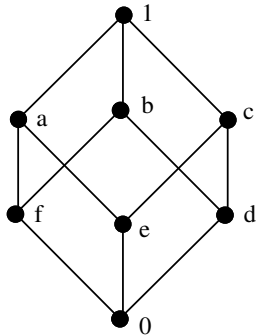

**Figure 2.4.2**

S[x] is a commutative polynomial ring having zero divisors. We are mainly interested in introducing the polynomial semiring for the construction of semivector spaces. Apart from this we are not going to give in this book any special properties enjoyed by the polynomial semirings or about the polynomials in them. For more about polynomial semirings one can refer Louis Dale. We can as in the case of rings study related properties about polynomial semirings.

**PROBLEMS**:

1. Let $Z^o[x]$ be the polynomial semiring. Find two ideals in $Z^o[x]$.

2. Can $Q^o[x]$ have ideals? If so find any two ideals in $Q^o[x]$.

3. Is S[x] where S is the semiring (i.e. the distributive lattice given by the following Hasse diagram).

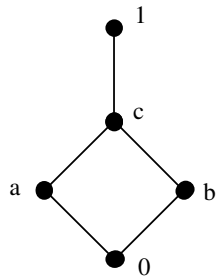

**Figure 2.4.3**

   a. Find ideals of S[x].
   b. Find subsemirings of S[x] which are not ideals of S[x].
   c. Find some zero divisors in S[x]\S.
   d. Can S[x]\S have idempotents?

4. Can $Z^o[x]$ have subsemirings which are not ideals?

5. Can $R^o[x]$ have ideals? Find subsemirings of $R^o[x]$ which are not ideals. Is $Q^o[x]$ a subsemiring which is not an ideal of $R^o[x]$?

6. Can $Z^o[x]$ have a polynomial p(x) such that p(a) = 0 for some $a \in Z^0 \setminus \{0\}$?

7. Define concepts like reducibility of polynomials in $Q^o[x]$ or in general in a polynomial semiring.

## 2.5 Group Semirings and Semigroup Semirings

This section is completely devoted to the introduction and study of group semirings and semigroup semirings which are not found in literature except one or two papers that have been done by the author in this field. The construction of these will give finite non-commutative semrings other than the ones discussed in section 2.3. Further this method will also give infinite non-commutative semirings of characteristic zero.

**DEFINITION 2.5.1**: *Let S be a strict commutative semiring with unit. G any group under multiplication. SG the group semiring of the group G over the semiring S is defined analogous to group rings and contains finite formal sums of the form $\sum_{i=1}^{n} s_i g_i$ where i runs over a finite number n with $s_i \in S$ and $g_i \in G$ satisfying the following conditions:*

1. $\sum_i s_i g_i = \sum_i t_i g_i \Leftrightarrow s_i = t_i$, $s_i, t_i \in S$ and $g_i \in G$ for all i.

2. $\left(\sum_i s_i g_i + \sum_i t_i g_i\right) = \sum_i (s_i + t_i) g_i$ for all $g_i \in G$ and $s_i, t_i \in S$.

3. $\left(\sum_i s_i g_i\right)\left(\sum_j t_j g_j\right) = \sum m_k p_k$ where $m_k = \Sigma s_i t_j$ and $p_k = g_i h_j$, $g_i, h_j \in G$ and $s_i, t_j \in S$.

4. $s_i g_i = g_i s_i$ for $s_i \in S$ and $g_i \in G$.

5. $s\left(\sum_i s_i g_i\right) = \sum_i s s_i g_i$ for $g_i \in G$ and $s, s_i \in S$.

6. As $1 \in G$ and $1 \in S$ we have $1 \bullet G = G \subseteq SG$ and $S \bullet 1 = S \subseteq SG$.

*The group semiring SG will be a commutative semiring if G is a commutative group. If G is non-commutative, the group semiring SG will be a non-commutative semiring.*

Let S be a semiring of infinite order and characteristic 0. Then SG will be a infinite semiring of characteristic 0, whatever be the group G. If S is a finite semiring say a finite distributive lattice then the group semiring SG will be a finite semiring provided G is a finite group, otherwise SG will be an infinite semiring but in both the cases for the semiring SG the characteristic remains undefined as we know for distributive lattices the characteristic does not exist.

Further we will show only in cases when we use the semiring as distributive lattices the group semiring enjoys a special property namely that the group elements of the group G can be represented in terms of the group semiring elements i.e. we can express $g \in G$ as $g = \alpha \bullet \beta$ where $\alpha, \beta \in SG$ and $\alpha, \beta \notin G$. This is a very unique and a different property enjoyed by the group semirings SG when S is a distributive non-chain lattice.

**DEFINITION 2.5.2**: *Let S be a strict commutative semiring with unit and L be a semigroup under multiplication with unit. The semigroup semiring SL is defined analogous to group semiring defined in this section; i.e. if we replace the group G by the semigroup L we get the semigroup semiring.*

*Example 2.5.1*: Let $Z^o$ be the semiring. $S_3$ be the symmetric group of degree 3. $Z^o S_3$ is the group semiring which is an infinite non-commutative semiring of characteristic zero. It can be checked that $Z^o S_3$ is also a strict semiring.

*Example 2.5.2*: Let $Q^o$ be the semiring. $G = \langle g / g^n = 1 \rangle$ be the cyclic group of order n. Clearly $Q^o G$ is the group semiring of infinite order which is commutative and has characteristic 0.

*Example 2.5.3*: Let the semiring S be the chain lattice, $C_2$. $S_4$ be the symmetric group of degree 4. The group semiring $SS_4$ is a non-commutative finite semiring with no characteristic associated with it.

*Example 2.5.4*: Let S be the semiring $C_3$. $G = \langle g \rangle$ be the infinite cyclic group. SG is an infinite commutative semiring and has no characteristic associated with it.

It is important to note that semirings are the only known structures with two binary operations for which the three distinct possibilities can occur.

1. The semirings can have zero characteristic; for example $Q^o$, $Z^o$, $R^o$.
2. Semirings have no finite or infinite characteristic defined. Example all distributive lattices.
3. Semirings with finite characteristic n where n is a positive integer.

It is left as an open research problem to construct semirings with finite characteristic n, n a positive integer.

*Example 2.5.5*: Let $Z^o$ be the semiring. $S(5)$ be the semigroup of all mappings of a five element set to itself under composition of mappings. $Z^o S(5)$ is the semigroup semiring which is non-commutative and has characteristic zero.

*Example 2.5.6*: Let $C_5$ be the chain lattice which is a semiring denoted by S. $S(3)$ the semigroup of mappings. The semigroup semiring $SS(3)$ is a finite non-commutative semiring with no characteristic associated with it.

*Example 2.5.7*: Let X be a finite set. P(X) the power set of X. L = P(X) is a semigroup under the operation '∪'. Let $Q^o$ be the semiring. $Q^o L$ is a semigroup semiring which is commutative and its characteristic is zero.

We do not deal and develop the properties of group semirings or semigroup semirings in this section. As our main motivation is only the study of Smarandache semirings and for this study these concepts will help us in giving concrete examples. We give only an example of how group elements are representable in terms of elements from group semiring and give a result regarding this.

*Example 2.5.8*: Let S be the semiring given by the following Hasse diagram:

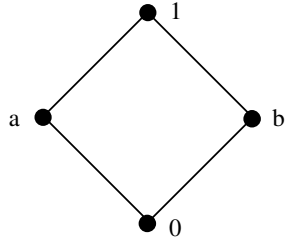

**Figure 2.5.1**

Let $G = \langle g/g^5 = 1\rangle$ be a cyclic group of degree 5. SG the group semiring of G over S. Let $\alpha = ag + bg^3$ and $\beta = ag^3 + bg$. Now $\alpha, \beta \in SG\backslash G$; but $\alpha \bullet \beta = (ag + bg^3)(ag^3 + bg) = ag^4 + bg^4 = (a + b)g^4$ (as $a + b = 1$) $= g^4 \in G$.

Thus we see $\alpha, \beta \notin G$ but $\alpha \bullet \beta \in G$. This is solely a new property enjoyed by group semirings when the semirings used are Boolean algebras or order greater than 2.

**THEOREM 2.5.1**: *Let $G = \langle g/g^n = 1\rangle$ be a finite cyclic group of order n. S be any Boolean algebra of order strictly greater than two. Then we have every element $g^i \in G$ can be represented as $\alpha \bullet \beta$ where $\alpha, \beta \in SG\backslash G$.*

*Proof*: Let $g^i \in G$ with $k + r = i$. Take $\alpha = (a_i g^k + a_i' g^r)$ and $\beta = (a_i' g^k + a_i g^r)$ where $a_i$ and $a_i'$ are atoms of S such that $a_i a_i' = 0$ and $a_i + a_i' = 1$. Now $\alpha \bullet \beta = (a_i g^k + a_i' g^r)(a_i' g^k + a_i g^r) = (a_i + a_i')g^{k+r} = 1.g^i = g^i \in G$.

This concept of group semirings and semigroup semirings will be used when we study Smarandache semirings.

**PROBLEMS**:

1. Find ideals and subsemirings in $Z^o S_3$.

2. Can $Q^o S_3$ have zero divisors? Justify.

3. Does $R^o G$ for any group G have idempotents?

4. Is $R^o G$ a strict semiring?

5. Let S be the semiring got from the distributive lattice with the following Hasse diagram

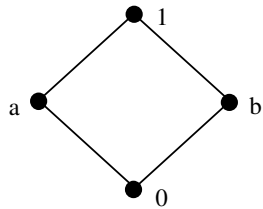

**Figure 2.5.2**

Can the group semiring $SS_4$ have zero divisors? idempotents? units?

6. Find right ideals of $SS_4$ given in problem 5.

7. Find two sided ideals in $Q^oS_3$.

8. Find subsemirings in $SS_4$ which are not ideals of $SS_4$.

9. Let $Z^oS(5)$ be the semigroup semiring.

   a. Find ideals in $Z^oS(5)$
   b. Find left ideals in $Z^oS(5)$
   c. Give a subsemiring in $Z^oS(5)$ which is not an ideal of $Z^oS(5)$
   d. Can $Z^oS(5)$ have zero divisors? units? idempotents?

10. Let S be a semiring given by the distributive lattice which has the following Hasse diagram:

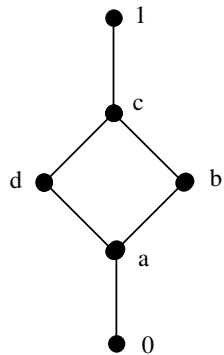

**Figure 2.5.3**

SS(3) be the semigroup semiring. Find

   a. The number of elements in SS(3)
   b. Can SS(3) have zero divisors?
   c. Find ideals of SS(3).
   d. Find a left ideal of SS(3).
   e. Find a subsemiring in SS(3) which is not an ideal.
   f. Can SS(3) have idempotents?
   g. Find invertible elements of SS(3).

## 2.6 Some Special semirings

In this section we introduce the new concepts of semirings given by Chris Monico, Zoltan Esik, Werner Kuich and others. We recall in this section the definition of ∗-semirings, Inductive ∗-semirings and c-simple semirings.

**DEFINITION (ZOLTAN ESIK AND WERNER KUICH)**: *A ∗-semiring is a semiring S equipped with a star operation ∗: S → S. Morphisms of ∗-semirings preserve the star operation.*

**DEFINITION (ZOLTAN ESIK AND WERNER KUICH)**: *An inductive ∗-semiring is a ∗-semiring, which is also an ordered semiring (An ordered semiring is a semiring S equipped with a partial order ≤ such that the operation are monotonic) and satisfies the fixed point inequation, $aa^o + 1 \leq a^*$ and the fixed point induction rule $ax + b \leq x \Rightarrow a * b < x$. A morphism of inductive ∗-semirings is an order preserving ∗-semiring morphism.*

*For more about these concepts please refer Inductive ∗-semirings by Zoltan Esik, Werner Kuich.*

**DEFINITION (ZOLTAN ESIK AND WERNER KUICH)**: *A continuous semiring is an ordered semiring S which is a complete partially ordered set with least element 0 and such that the sum and product operations are continuous. A semiring, which is both a ∗-semiring and a continuous semiring is called continuous ∗-semiring if the star operations on S is given by $a^* = \sum_{n \geq 0} a^n$ for all $a \in S$. It is proved a continuous ∗-semiring is an inductive ∗-semiring.*

They proceed on to define several equations to be satisfied by semiring. We ask the reader to refer this paper for the several definitions and results quoted in this book.

**DEFINITION (MONICO, CHRIS)**: *A congruence relation on a semiring S is an equivalence relation ~ that also satisfies*

$$x_1 \sim x_2 \Rightarrow \begin{cases} c + x_1 \sim c + x_2 \\ x_1 + c \sim x_2 + c \\ cx_1 \sim cx_2 \\ x_1 c \sim x_2 c \end{cases}$$

*for all $x_1, x_2 \in S$. A semiring that admits no congruence relation other than the trivial ones, identity S and $S \times S$ is said to be congruence simple or c-simple.*

He has proved if S is a finite field then S is a commutative c-simple finite semiring. We use his result and prove Smarandache analog. For more about these concepts please refer Chris Monico.

## Supplementary reading

1. Hebisch, Udo and Hanns Joachim Weinert. *Semirings and Semifield*, in "Handbook of Algebra", Vol. 1, Elsevier Science, Amsterdam, 1996.

2. Louis Dale. *Monic and Monic-free ideals in a polynomial semiring.* PAMS 56 45-50, 1976.

3. Monico, Chris *On finite congruence simple semiring*. http://arxiv.org/PS_cache/math/pdf/0205/0205083.pdf

# CHAPTER THREE
# SEMIFIELDS AND SEMIVECTOR SPACES

This chapter is completely devoted to the introduction of the concepts of semifields and semivector spaces. Here we give some examples and recall some interesting properties enjoyed by them. Neither this book nor this chapter does claim any complete recollection of all notions existing about semivector spaces and semifields, as the main motivation is only to study and introduce most of the basic properties about semivector spaces and semifields in the context of Smarandache notions.

This chapter has 3 sections, in the first section we introduce and study semifields, in the second section define and give examples of semivector spaces. In the final section recall some of the most important properties of semivector spaces.

## 3.1 Semifields

This section introduces the concept of semifields and give some examples and derive a few of the important properties about them.

**DEFINITION 3.1.1**: *Let S be a non-empty set. S is said to be a semifield if*

1. *S is a commutative semiring with 1.*
2. *S is a strict semiring. That is for a, b $\in$ S if a + b = 0 then a = 0 and b = 0.*
3. *If in S, a • b = 0 then either a = 0 or b = 0.*

*Example 3.1.1*: Let $Z^o$ be the semiring. $Z^o$ is a semifield.

*Example 3.1.2*: $Q^o$ is a semifield.

**DEFINITION 3.1.2**: *The semifield S is said to be of characteristic zero if 0 • x = 0 and for no integer n; n • x = 0 or equivalently if x $\in$ S \ {0}, nx = x + … + x, n times equal to zero is impossible for any positive integer n.*

*Example 3.1.3*: $R^o$ is a semifield of characteristic 0. The semifields given in examples 3.1.1 and 3.1.2 are semifields of characteristic 0.

**DEFINITION 3.1.3**: *Let S be a semifield; a subset N of S is a subsemifield if N itself is a semifield. If N is a subsemifield of S, we can call S as an extension semifield of N.*

*Example 3.1.4*: Let $Z^o$ and $Q^o$ be semifields. $Z^o$ is the subsemifield of $Q^o$ and $Q^o$ is an extension semifield $Z^o$.

*Example 3.1.5*: Let $Z^o$ and $R^o$ be semifields. $Z^o$ is the subsemifield of $R^o$ and $R^o$ is an extension semifield of $Z^o$.

***Example 3.1.6***: Let $Z^o[x]$ be the polynomial semiring which is a semifield. Clearly $Z^o[x]$ is an extension semifield of $Z^o$ and $Z^o$ is the subsemifield of $Z^o[x]$.

It is pertinent to note the following facts, which are put as remarks.

***Remark 1***: As in the case of fields of characteristic 0, semifield of characteristic zero also has infinite number of elements in them.

***Remark 2***: Unlike in fields where, R[x], polynomial ring is only an integral domain in case of semifields, the polynomial semiring is also a semifield. If S is a semifield, S[x] the polynomial semiring is an extension semifield of S that is $S \subseteq S[x]$. This is a special and a distinct property enjoyed solely by semifields and not by fields.

***Example 3.1.7***: Let $C_3$ be a chain lattice. $C_3$ is a semifield. $C_3[x]$ is a polynomial semiring, is an extension semifield of $C_3$ and $C_3$ is the subsemifield of $C_3[x]$.

Clearly $C_3[x]$ has no characteristic associated with it. In fact $C_3[x]$ is an infinite extension of the finite semifield $C_3$.

The following theorem is left as an exercise for the reader to prove.

**THEOREM 3.1.1**: *Every chain lattice is a semifield with no characteristic associated with it.*

We find the following two observations as appropriate corollaries.

**COROLLARY 3.1.2**: *Every distributive lattice, which is not a chain lattice, in general is not a semifield.*

*Proof*: By an example, consider the lattice S, whose Hasse diagram is given by the following figure.

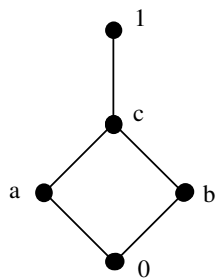

**Figure 3.1.1**

S is a semiring but S is not a semifield as $a \bullet b = 0$ for $a, b \in S$.

***Example 3.1.8***: Let S be the semiring given by the following Hasse diagram:

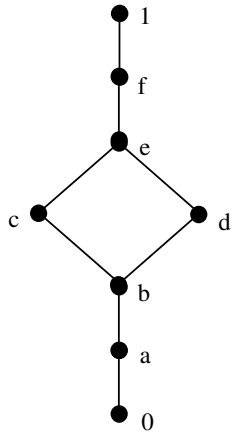

**Figure 3.1.2**

Clearly S is a semifield but S is not a chain lattice as S cannot be totally ordered, for c, d in S are not comparable. Thus we have semifields of this from also.

**COROLLARY 3.1.3**: *Let B be a Boolean algebra of order greater than two, B is never a semifield.*

*Proof*: Since B is a Boolean algebra which is not $C_2$ we see B has at least 2 atoms $a_1$, $a_2$ in B and if $a_1'$ is the complement of $a_1$ we have $a_1.a_1' = 0$ and $a_1' + a_1 = 1$. Thus B is not a semifield if $o(B) > 2$.

As in the case of fields we have the smallest prime field of characteristic 0 is Q, like wise in case of semifield we have the semifield $Z^o$ to be a prime semifield as it has no subsemifields.

**DEFINITION 3.1.4**: *Let S be a semifield, we say S is a prime semifield if S has no proper subsemifield.*

*Example 3.1.9*: Let $Z^o$ be the semifield. $Z^o$ has no proper subsemifields so $Z^o$ is the prime semifield of characteristic 0.

**THEOREM 3.1.4**: *Every semifield of characteristic zero contains $Z^o$ as a subsemifield.*

*Proof*: Let S be a semifield of characteristic 0. Since $1 \in S$ we see 1 generates $Z^o$ so $Z^o$ is a subsemifield of S.

**THEOREM 3.1.5**: *$Z^o$ is the smallest semifield of characteristic 0.*

*Proof*: $Z^o$ has no proper subsemifield. Since any subsemifield N of $Z^o$ must also contain 0 and 1 so $Z^o \subset N$ but $N \subset Z^o$ so $N = Z^o$.

From here onwards we will call the semifield $Z^o$ as the prime semifield of characteristic 0. $C_2$ the chain lattice is a prime semifield as $C_2$ has no subsemifields. This leads us to formulate the following important theorem.

**THEOREM 3.1.6**: *Let $C_n$ be a chain lattice with $n > 2$. $C_n$ is not a prime semifield.*

*Proof*: Now $C_n$ when $n > 2$ we have the subset $S = \{0, 1\}$ to be a subsemifield of $C_n$. So $C_n$ is not a prime semifield.

This gives the following:

**COROLLARY 3.1.7**: *$C_2$ is the prime semifield having no characteristic.*

*Proof*: Obvious as $C_2$ contains only 0 and 1.

**THEOREM 3.1.8**: *Let S be a distributive lattice, which is a semifield having more than 2 elements. S is not a prime semifield.*

*Proof*: If S is a semifield, S contains 0 and 1. So S has a subsemifield given by $\{0, 1\}$; thus S is not a prime semifield.

Just as in the case of fields direct product of semifields is not a semifield.

***Example 3.1.10***: Let $S = Z^o \times Z^o \times Z^o$ be the direct product of semifields. S is not a semifield as $\alpha = (0, 9, 0)$ and $\beta = (1, 0, 2) \in S$ is such that $\alpha.\beta = (0, 0, 0)$.

The major property which distinguishes the semifield and field is given by the following theorem. We know a field has no ideals. But semifields have ideals.

**THEOREM 3.1.9**: *Let $Z^o$ be a semifield. $Z^o$ has ideals.*

*Proof*: $nZ^o$ for any positive integer n is an ideal of $Z^o$.

**PROBLEMS**:

1. Find a subsemifield of $Q^o[x]$.

2. Find in $Z^o[x]$ a subsemifield.

3. Find ideals in $C_7[x]$.

4. Find subsemirings of $C_{11}[x]$ which are not ideals.

5. Can $C_{13}[x]$ have subsemirings which are not ideals?

6. In $Z^o$ can we say all subsemirings are ideals?

7. Find examples of semifields in which every subsemiring is an ideal.

## 3.2 Semivector Spaces and Examples

In this section we introduce the concept of semivector spaces and give examples of them. Since the notion of semivector spaces are not found in text and we have only one paper on semivector spaces published in 1993. So we have given a complete recaptivation of these concepts with examples.

**DEFINITION 3.2.1**: *A semivector space V over the semifield S of characteristic zero is the set of elements, called vectors with two laws of combination, called vector addition (or addition) and scalar multiplication, satisfying the following conditions:*

1. *To every pair of vectors $\alpha$, $\beta$ in V there is associated a vector in V called their sum, which we denote by $\alpha + \beta$.*
2. *Addition is associative $(\alpha + \beta) + \gamma = \alpha + (\beta + \gamma)$ for all $\alpha, \beta, \gamma \in V$.*
3. *There exists a vector, which we denote by zero such that $0 + \alpha = \alpha + 0 = \alpha$ for all $\alpha \in V$.*
4. *Addition is commutative $\alpha + \beta = \beta + \alpha$ for all $\alpha, \beta \in V$.*
5. *If $0 \in S$ and $\alpha \in V$ we have $0. \alpha = 0$.*
6. *To every scalar $s \in S$ and every vector $v \in V$ there is associated a unique vector called the product s.v which is denoted by sv.*
7. *Scalar multiplication is associative, $(ab) \alpha = a (b\alpha)$ for all $\alpha \in V$ and $a, b \in S$.*
8. *Scalar multiplication is distributive with respect to vector addition, $a (\alpha + \beta) = a\alpha + a\beta$ for all $a \in S$ and for all $\alpha, \beta \in V$.*
9. *Scalar multiplication is distributive with respect to scalar addition: $(a + b) \alpha = a\alpha + b\alpha$ for all $a, b \in S$ and for all $\alpha \in V$.*
10. *$1. \alpha = \alpha$ (where $1 \in S$) and $\alpha \in V$.*

*Example 3.2.1*: Let $Z^o$ be the semifield. $Z^o[x]$ is a semivector over $Z^o$.

*Example 3.2.2*: Let $Q^o$ be the semifield. $R^o$ is a semivector space over $Q^o$.

*Example 3.2.3*: $Q^o$ is a semivector space over the semifield $Z^o$.
It is important to note that $Z^o$ is not a semivector space over $Q^o$. Similarly $Z^o$ is not a semivector space over $R^o$.

*Example 3.2.4*: Let $M_{n \times n} = \{(a_{ij}) \mid a_{ij} \in Z^o\}$; the set of all $n \times n$ matrices with entries from $Z^o$. Clearly $M_{n \times n}$ is a semivector space over $Z^o$.

*Example 3.2.5*: Let $V = Z^o \times Z^o \times \ldots \times Z^o$ (n times), V is a semivector space over $Z^o$. It is left for the reader to verify.

*Example 3.2.6*: Let $C_n$ be a chain lattice. $C_n[x]$ is a semivector space over $C_n$.

*Example 3.2.7*: Let $V = C_n \times C_n \times C_n$, V is a semivector space over $C_n$. It is left for the reader to verify.

## 3.3 Properties about semivector spaces.

In this section we introduce some basic concepts like linear combination, linearly dependent, linearly independent, spanning set and subsemivector spaces and illustrate them with examples. Further we state several interesting results about semivector spaces and leave the proof of them as problems to the reader, as this text assumes the reader of a strong background in algebra. As the concept of semivector spaces was introduced in 1993 so we proceed to define related properties of semivector spaces.

Let S be a semifield and V be a semivector space over S. If $\beta = \sum \alpha_i v_i$ ($v_i \in V$, $\alpha_i \in$ S) which is in V, we use the terminology $\beta$ is a linear combination of $v_i$'s. We also say $\beta$ is linearly dependent on $v_i$'s if $\beta$ can be expressed as a linear combination of $v_i$'s We see the relation is a non trivial relation if at least one of the coefficients $\alpha_i$'s is non zero. This set $\{v_1, v_2, \ldots, v_k\}$ satisfies a non trivial relation if $v_j$ is a linear combination of $\{v_1, v_2, \ldots, v_{j-i}, v_{j+i}, \ldots, v_k\}$.

**DEFINITION 3.3.1**: *A set of vectors in V is said to be linearly dependent if there exists a non-trivial relation among them; otherwise the set is said to be linearly independent.*

*Example 3.3.1*: Let $Z^o[x]$ be the semivector space over $Z^o$. Now the set $\{1, x, x^2, x^3, \ldots, x^n\}$ is a linearly independent set. But if we consider the set $\{1, x, x^2, x^3, \ldots, x^n, x^3 + x^2 + 1\}$ it is a linearly dependent set.

**THEOREM 3.3.1**: *Let V be a semivector space over the semifield S. If $\alpha \in V$ is linearly dependent on $\{\beta_i\}$ and each $\beta_i$ is linearly dependent on $\{\gamma_i\}$ then $\alpha$ is linearly dependent on $\{\gamma_i\}$.*

*Proof*: Let $\alpha = \sum_i b_i \beta_i \in V$; $\beta_i \in V$ and $\beta_i = \sum c_{ij} \gamma_j$ for each i and $\gamma_j \in V$ and $c_j \in S$. Now $\alpha = \sum b_i \beta_i = \sum_i b_i \sum_j c_{ij} \gamma_j = \sum_j \left( \sum_i b_i c_{ij} \right) \gamma_j$. as $b_i c_{ij} \in S$ and $\gamma_j \in V$, $\alpha \in V$. Hence the claim.

The main difference between vector spaces and semivector spaces is that we do not have negative terms in semifields over which semivector spaces are built. So as in the case of vector spaces we will not be in a position to say if $\alpha_1 v_1 + \ldots + \alpha_n v_n = 0$ implies $v_1 = \frac{-1}{\alpha_1}(\alpha_2 v_2 + \ldots + \alpha_n v_n)$. To overcome this difficulty we have to seek other types of arguments. But this simple property has lead to several distortion in the nature of semivector spaces as we cannot define dimension, secondly many semivector spaces have only one basis and so on. Here also we do not go very deep into the study of semivector spaces as the chief aim of this book is only on the analogues study of Smarandache notions. So for the reader we have suggested at the end of this section the books or paper to be read as supplementary reading. This book at the end of each chapter has a set of books and papers enlisted for supplementary reading, which is a main feature of this book.

**DEFINITION 3.3.2**: *Let V be a semivector space over the semifield S. For any subset A of V the set of all linear combination of vectors in A, is called the set spanned by A and we shall denote it by $\langle A \rangle$. It is a part of this definition, $A \subset \langle A \rangle$.*

*Thus we have if $A \subset B$ then $\langle A \rangle \subset \langle B \rangle$.*

Consequent of this we can say the theorem 3.3.1 is equivalent to; if $A \subset \langle B \rangle$ and $B \subset \langle C \rangle$ then $A \subset \langle C \rangle$.

**THEOREM 3.3.2**: *Let V and S be as in the earlier theorem. If B and C be any two subsets of V such that $B \subset \langle C \rangle$ then $\langle B \rangle \subset \langle C \rangle$.*

*Proof*: Set $A = \langle B \rangle$ in theorem 3.3.1, then $B \subset \langle C \rangle$ implies $\langle B \rangle = A \subset \langle C \rangle$.

Now we still have many interesting observations to make if the set A is a linearly dependent set. We have the following theorem:

**THEOREM 3.3.3**: *Let V be a semivector space over S. $A = \{\alpha_1, \ldots, \alpha_k\}$ be a subset of V. If $\alpha_i \in A$ is dependent on the other vectors in A then $\langle A \rangle = \langle A \setminus \{\alpha_i\} \rangle$.*

*Proof*: The assumption is that $\alpha_i \in A$ is dependent on $A \setminus \{\alpha_i\}$, means that $A \subset \langle A \setminus \{\alpha_i\} \rangle$. It then follows from theorem 3.3.2 that $\langle A \rangle \subseteq \langle A \setminus \{\alpha_i\} \rangle$. Equality follows from the fact that the inclusion in other direction is evident.

**DEFINITION 3.3.3**: *A linearly independent set of a semivector space V over the semifield S is called a basis of V if that set can span the semivector space V.*

*Example 3.3.2*: Let $V = Z^o \times Z^o \times Z^o$ be a semivector space over $Z^o$. The only basis for V is $\{(1, 0, 0), (0, 1, 0), (0, 0, 1)\}$ no other set can span V.

This example is an evidence to show unlike vector spaces which can have any number of basis certain semivector spaces have only one basis.

*Example 3.3.3*: Let $Z^o$ be a semifield. $Z^o_n[x]$ denote the set of all polynomials of degree less than or equal to n. $Z^o_n[x]$ is a semivector space over $Z^o$. The only basis for $Z^o_n[x]$ is $\{1, x, x^2, \ldots, x^n\}$.

We have very distinguishing result about semivector spaces, which is never true in case of vector spaces.

**THEOREM 3.3.4**: *In a semivector space V, over the semifield S, the number of elements in a set which spans V need not be an upper bound for the number of vectors that can be linearly independent in V.*

*Proof*: The proof is given by an example. Consider the semivector space $V = Z^o \times Z^o$ over $Z^o$. We have $\{(0, 1), (1, 0)\}$ to be the only basis of V. In particular this set spans V. But we can find in $V = Z^o \times Z^o$ three vectors which are linearly independent. For example the set $U = \{(1, 1), (2, 1), (3, 0)\}$ is a linearly independent set in V for none of them is expressible in terms of the others. But note that this set is not a basis for V as this set does not span V. This can be found from the fact $(1, 3) \in V$ but it is not expressible as a linear combination of elements from U. Note that $U \cup \{(1, 3)\}$ is a linearly independent set.

**THEOREM 3.3.5**: *Let V be a semivector space over the semifield S. For any set $C \subset V$ we have $\langle\langle C \rangle\rangle = \langle C \rangle$.*

*Proof*: Clearly $\langle C \rangle \subseteq \langle\langle C \rangle\rangle$. Now replacing $B = \langle C \rangle$ in theorem 3.3.2 we have $\langle\langle C \rangle\rangle \subseteq \langle C \rangle$. Hence we have the equality $\langle C \rangle = \langle\langle C \rangle\rangle$.

**DEFINITION 3.3.4**: *Let V be a semivector space over the semifield S with the property that V has a unique basis. Then we define the number of elements in the basis of V to be the dimension of V.*

A semivector space is finite dimensional if the space has a unique basis and the number of elements in it is finite.

***Example 3.3.4***: Let $V = Z^o \times Z^o \times Z^o \times Z^o \times Z^o$ (5 times) be a semivector space over $Z^o$. Clearly dimension of V is 5.

***Example 3.3.5***: Let $V = Z_7^0[x] = \{$set of all polynomials of degree less than or equal to 7 with coefficient from $Z^o\}$ be a semivector space over $Z^o$. The dimension of V is 8.

**THEOREM 3.3.6**: *In a n- dimensional semivector space we need not have in general, every set of n + 1 vectors to be linearly independent.*

*Proof*: We have seen in case of the semivector space $V = Z^0 \times Z^0$ over $Z^0$, which is of dimension 2, has the 3 vectors $\{(1, 1), (2, 1), (3, 0)\}$ to be linearly independent.

This is a unique property enjoyed only by semivector spaces which can never be true in case of vector space for in case of vector spaces we know if dimension of a vector space is n then the vector space cannot contain more than n linearly independent vectors.

Now we proceed on to build semivector spaces using lattices.

**THEOREM 3.3.7**: *All finite lattices are semivector spaces over the two element Boolean algebra $C_2$.*

*Proof*: Let $(L, \cup, \cap)$ be any finite lattice. We see L is a semivector space over $C_2$. All axioms of the semivector space are true. To show scalar multiplication is distributive we have to prove $s \cap (a \cup b) = (s \cap a) \cup (s \cap b)$ for all $a, b \in L$ and $s \in C_2 = (0, 1)$. The two values s can take are either $s = 0$ or $s = 1$.

In case $s = 0$ we have $0 \cap (a \cup b) = 0$ (zero is the least element), $(0 \cap a) \cup (0 \cap b) = 0$.

When $s = 1$ we have $1 \cap (a \cup b) = a \cup b$ (1 is the greatest element)

$(1 \cap a) \cup (1 \cap b) = a \cup b$. Hence the claim.

Thus we can easily verify all lattices L whether L is a distributive lattice or otherwise L is a semivector space over the semifield $C_2$.

**DEFINITION 3.3.5**: *A subsemivector space W of a semivector space V over a semifield S is a non-empty subset of V, which is itself, a semivector space with respect to the operations of addition and scalar multiplication.*

<u>Remark</u>: The whole space V and the {0} element are trivially the subsemivector spaces of every semivector space V.

*Example 3.3.6*: $R^o$ is a semivector space over $Z^o$. $Q^o$ the subset of $R^o$ is a non-trivial subsemivector space of $R^o$.

*Example 3.3.7*: $Z^o[x]$ is a semivector space over $Z^o$. All polynomials of even degree in $Z^o[x]$ denoted by S is a subsemivector space over $Z^o$.

*Example 3.3.8*: Let $C_2$ be the semifield, the lattice L with the following Hasse diagram is a vector space over $C_2$.

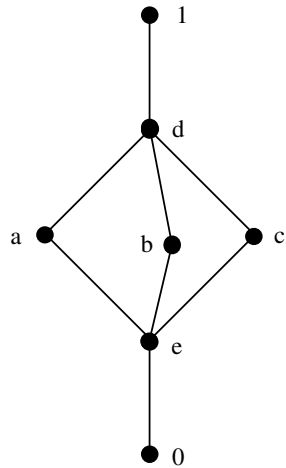

**Figure 3.3.1**

Let $S \subset L$ where $S = \{1, a, e, d, 0\}$ is subsemivector space over $C_2$.

*Example 3.3.9*: Let $M_{2 \times 2} = \left\{ \begin{pmatrix} a & b \\ c & d \end{pmatrix} \Big/ a,b,c,d \in C_8 \right\}$ be the semivector space over $C_8$.

Let $A = \left\{ \begin{pmatrix} a & 0 \\ 0 & 0 \end{pmatrix} \Big/ a \in C_8 \right\}$ be the subset of $M_{2 \times 2}$. A is a subsemivector space over $C_8$ ($C_8$ the chain lattice with 8 elements).

**THEOREM 3.3.8**: *Let V be a semivector space over the semifield S. The intersection of any collection of subsemivector spaces is a subsemivector space of V.*

*Proof*: The proof is as in the case of vector spaces.

Several other results from vector spaces can be adopted whenever possible, and at times it may not be possible to adopt them, So as the case may be, extension or introduction of the properties of vector spaces is done in case of semivector spaces. A property, which is different basically, is given below.

**THEOREM 3.3.9**: *In a semivector space an element need not in general have a unique representation in terms of its basis elements.*

*Proof*: This is proved by an example. Let $C_4$ be the chain lattice with elements say (0, b, a, 1) $0 < b < a < 1$. $C_4$ is a semivector space over $C_2$. Note the set $\{1, a, b\}$ is a linearly independent set as none of the elements in the given set is expressible as a linear combination of others. Further $\{1, a, b\}$ spans $C_4$. So $\{1, a, b\}$ is a basis, in fact a unique basis. It is interesting to note that the elements a and 1 do not have a unique representation in terms of the basis elements for we have

$$a = 1 \bullet a + 0 \bullet b + 0 \bullet 1$$
$$= 1 \bullet a + 1 \bullet b + 0 \bullet 1$$
$$1 = 0 \bullet a + 0 \bullet b + 1 \bullet 1$$
$$= 1 \bullet a + 1 \bullet b + 1 \bullet 1$$
$$= 1 \bullet a + 0 \bullet b + 1 \bullet 1$$
$$= 0 \bullet a + 1 \bullet b + 1 \bullet 1.$$

This is a unique feature enjoyed by semivector spaces built using lattices.
Lastly we proceed to define linear transformation and linear operators on semivector spaces.

**DEFINITION 3.3.6**: *Let $V_1$ and $V_2$ be any two semivector spaces over the semifield S. We say a map / function $T: V_1 \to V_2$ is a linear transformation of semivector spaces if $T(av + u) = aT(v) + T(u)$ for all $u, v \in V_1$ and $a \in S$.*

*Example 3.3.10*: Let $V = Z^o \times Z^o \times Z^o$ and $Z_6^o[x]$ be semivector spaces defined over $Z^o$. Define T: $V \to Z_6^o[x]$ by

$T(1, 0, 0) = x^6 + x^5$
$T(0, 1, 0) = x^3 + x^4$
$T(0, 0, 1) = x^2 + x + 1$
$T(9(3, 2, 1) + 6(1, 3, 0)) = 9[3(x^6 + x^5) + 2(x^3 + x^4) + x^2 + x + 1] + 6[x^6 + x^5 + 3(x^3 + x^4) + 0(x^2 + x + 1)] = 33x^6 + 33x^5 + 36x^4 + 36x^3 + 9x^2 + 9x + 9$.
$T(9(3, 2, 1) + 6(1, 3, 0)) = T((33, 36, 9)) = 33x^6 + 33x^5 + 36x^4 + 36x^3 + 9x^2 + 9x + 9$.

**DEFINITION 3.3.7**: *Let V be a semivector space over S. A map/ function T from V to V is called a linear operator of the semivector space V if $T(\alpha v + u) = \alpha T(v) + T(u)$ for all $\alpha \in S$ and $v, u \in V$.*

*Example 3.3.11*: Let $V = Z^o \times Z^o \times Z^o \times Z^o$ be a semivector space over $Z^o$. Define T: $V \to V$ by

$$T(1, 0, 0, 0) = (0, 1, 0, 0)$$
$$T(0, 1, 0, 0) = (0, 0, 1, 0)$$
$$T(0, 0, 1, 0) = (0, 0, 0, 1)$$
$$T(0, 0, 0, 1) = (1, 0, 0, 0)$$

It can be verified T is a linear operator on V.

**PROBLEMS**:

1. Is $Z^o \times Z^o$ an extension field of $Z^o$? Justify.

2. $C_7$ is a semifield, find an extension semifield of $C_7$.

3. $C_{20}$ is a semifield find a subsemifield of $C_{20}$.

4. Prove $C_9 \times C_{10} \times C_{11}$ is only a semiring and not a semifield.

5. Find a semiring of order 11, which is not a semifield.

6. Find an example of a non-commutative semiring with 16 elements.

7. Give an example of a lattice which not a semiring.

8. What is dimension of the semivector space $C_2[x]$ over $C_2$?

9. Can we find a dimension for the lattice L given by the following Hasse diagram which is a semivector space over $C_2$

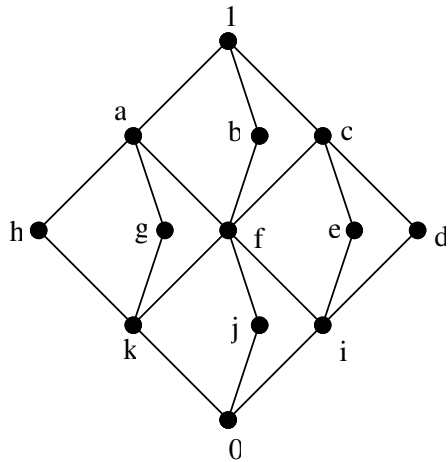

**Figure 3.1.2**

10. Can L be a semivector space over $C_8$? Substantiate your answer?

11. How many sets of basis does the semivector space L described in problem 9 have?

12. Let $M_{2 \times 2} = \{(a_{ij}) / a_{ij} \in Z^o\}$ be a semivector space over $Z^o$.

i. Find a basis for $M_{2\times 2}$.
  ii. Is it a unique basis?
  iii. Find the dimension of $M_{2\times 2}$ over $Z^o$.

13. Let $M_{2\times 2} = \{(a_{ij}) / a_{ij} \in Q^0\}$. $M_{2\times 2}$ is a semivector space over $Q^o$.

    i. Find a basis for $M_{2\times 2}$.
    ii. Does $M_{2\times 2}$ have dimension defined?

14. Let $M_{2\times 2} = \{(a_{ij}) / a_{ij} \in Q^o\}$ be a semivector space over $Z^o$.

    i. Find a basis for $M_{2\times 2}$.
    ii. Can we have dimension for $M_{2\times 2}$?

15. Let $M_{2\times 2} = \{(a_{ij}) / a_{ij} \in R^o\}$ be a vector space over $Z^o$.

    i. Find a basis for $M_{2\times 2}$.
    ii. Does a dimension exist for $M_{2\times 2}$ over $Z^o$?

16. Is $Q^o$ a finite dimensional semivector space over $Z^o$?

17. Is $R^o$ a finite dimensional semivector space over $Q^o$?

18. Is $R^o$ a finite dimensional semivector space over $Z^o$?

19. Is $R^o \times R^o$ a finite dimensional semivector space over $R^o$?

## Supplementary Reading

1. Vasantha Kandasamy, W. B. *Semivector spaces over semifield.* Mathematika 188, 43-50, 1993.

2. Vasantha Kandasamy, W.B. *On a new class of semivector spaces.* Varahmihir Jour of Math. Sci., Vol. 1, 23-30, 2001.

# CHAPTER FOUR
# SMARANDACHE SEMIRINGS

In this chapter we introduce and analyse the concept of Smarandache semiring. This chapter has six sections. The first section is devoted just to the introduction of the concept of Smarandache semiring and explaining them by illustrative examples. In the second section we introduce the substructures like Smarandache subsemiring, Smarandache ideals, Smarandache commutative semirings, Smarandache pseudo subsemirings, Smarandache dual ideals and so on and obtain some interesting results about them.

The next section is devoted to the study of special elements in Smarandache semirings (S-semirings). Using the recent literature of semirings corresponding S-semirings are defined and analysed in the chapter on special S-semirings. As in Smarandache notions there can be layers or levels of S-semigroups, this section defines S-semirings of level II so that it goes without saying that S-semirings which are studied till this point are S-semirings of level I. It is pertinent to mention here that except the concept of idempotent semiring all substructures constructed using S-semiring of level II are distinctly different from S-semiring of level I. The S-semiring of level II is achieved only by Smarandache mixed direct product. As all Smarandache notions are in a way non-classical mathematics, the final section is devoted to the study of Smarandache anti-semirings. Each section is filled with examples and problems for the student to take things in a serious way. As these solving process of problems will create in researchers and students a gripping interest in these concepts.

## 4. 1 Definition of S-semirings and examples

In this section we define Smarandache semirings and illustrate them by several examples. As this notion is very recent (2001) we explain it elaborately and substantiate it with examples. Further even the notion of the semirings is not found in literature in the form of specialized text books so we felt it essential to explain S-semirings.

**DEFINITION (SMARANDACHE, FLORENTIN)**: *The Smarandache semiring S which will be denoted from here onwards as S-semiring is defined to be a semiring S such that a proper subset B of S is a semifield (with respect to the same induced operation). This is $\phi \neq B \subset S$.*

*Example 4.1.1*: Let $Q^o$ be the semiring. $Q^o$ is a S-semiring, as $Z^o \subset Q^o$ is a proper subset which is a semifield.

*Example 4.1.2*: Let $Z^o$ be a semiring. $Z^o$ is not a S-semiring. $Z^o$ does not have any proper subset which is a semifield. In view of this we have the following.

**THEOREM 4.1.1:** *Every semiring in general need not be a S-semiring.*

*Proof*: By an example. $Z^o$ is a semiring which is a not a S-semiring. Also take the chain lattice $C_2$; $C_2$ is not a S-semiring for $C_2$ cannot have a proper subset which is a semifield.

***Example 4.1.3:*** Let $Z^o[x]$ be the polynomial semiring. $Z^o[x]$ is a S-semiring as $Z^o \subset Z^o[x]$ is a semifield.

***Example 4.1.4***: Let $C_2[x]$ be the polynomial semiring. $C_2[x]$ is a S-semiring, as $C_2 \subset C_2[x]$ is a semifield.

**THEOREM 4.1.2**: *All distributive lattices L with 0 and 1 and having more than 2 elements are S-semirings.*

*Proof*: Given L is a distributive lattice with 0 and 1. So L is a semiring. Since L contains 0 and 1, the set $A = (0, 1) \subset L$; A is a semifield so L is a S-semiring.

We have seen examples of commutative S-semirings. Now we will see the structure of both finite and infinite non-commutative S-semirings.

***Example 4.1.5***: Let $M_{n \times n} = \{(a_{ij})/a_{ij} \in Z^o\}$ be the set of all $n \times n$ matrices with entries from $Z^o$. $M_{n \times n}$ is a semiring under matrix addition and matrix multiplication. Clearly $M_{n \times n}$ is a non-commutative semiring. $M_{n \times n}$ is a S-semiring for if we take $A = \{(a_{ij})/a_{ij} = 0$ if $i \neq j$ and $a_{ii} \in Z^o \setminus \{0\}\} \cup \{$The zero matrix$\}$. Then $A \subset M_{n \times n}$ is a semifield so $M_{n \times n}$ is a S-semiring.

**DEFINITION 4.1.1**: *If the S-semiring has only finite number of elements we say the S-semiring is finite otherwise infinite.*

***Example 4.1.6***: Let $M_{3 \times 3}$ be the set of all $3 \times 3$ matrices with entries from the distributive lattice L having the following Hasse diagram $M_{3 \times 3}$ with '$\cup$' and '$\cap$' defined as in the earlier chapters in a semiring. $M_{3 \times 3}$ is a S-semiring. For take $A = \{(a_{ij})/a_{ij} = 0$ if $i \neq j$ and $a_{ii} \in L \setminus \{0\}\} \cup \left\{\begin{matrix} 0 & 0 & 0 \\ 0 & 0 & 0 \\ 0 & 0 & 0 \end{matrix}\right\} \subset M_{3 \times 3}$, A is a semifield. Hence $M_{3 \times 3}$ is a S-semiring.

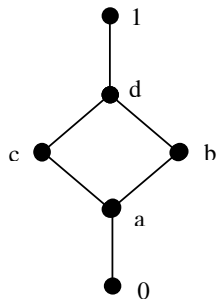

**Figure 4.1.1**

***Example 4.1.7***: Let $Z^o$ be the semiring. $S = Z^o \times Z^o \times Z^o \times Z^o$ is a semiring under component wise addition and multiplication. S is a S-semiring but $Z^o$ is not a S-semiring.

Thus the concept of direct product helps us to convert a non S-semiring into a S-semiring. Of course we will still see many more interesting properties about these direct products.

**PROBLEMS**:

1. Show the lattice given by the following Hasse diagram is a S-semiring

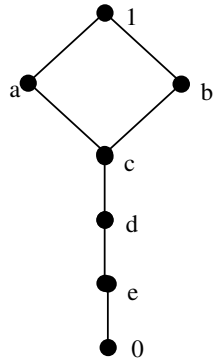

**Figure 4.1.2**

2. Is the semiring given in Problem 1 a semiring with unit?
3. Let $S = Z^o \times Z^o \times Z^o$; show S is a S-semiring.
4. Find the smallest semiring which is a S-semiring.
5. Prove the following: A commutative semiring $(S, +, \bullet)$ given by the following tables. Is it a S-semiring?

| + | 0 | 1 |
|---|---|---|
| 0 | 0 | 1 |
| 1 | 1 | 0 |

| $\bullet$ | 0 | 1 |
|---|---|---|
| 0 | 0 | 0 |
| 1 | 0 | 1 |

(This semiring is taken from the paper of Chris Monico.)

6. Prove all polynomial semirings R[x] where R is a semifield is a S-semiring.
7. Give an example of a semiring which is not a S-semiring.

## 4.2 Substructures in S-semirings

In this section we introduce many substructures in S-semirings like S-subsemirings, S-commutative semirings, S-non-commutative semirings, S-right(left) ideals, S-ideals in semirings and S-quotient semirings; study them and illustrate them with examples.

**DEFINITION 4.2.1**: *Let S be a semiring. A non-empty proper subset A of S is said to be a Smarandache subsemiring (S-subsemiring) if A is S-semiring i.e. A has a proper subset P such that P is a semifield under the operations of S.*

An immediate consequence of this is the following result.

**THEOREM 4.2.1**: *Let S be a semiring having a S-subsemiring then S is a S-semiring.*

*Proof*: Given S is a semiring and A ⊂ S is a S-subsemiring of S; so A has a proper subset P ⊂ A such that P is a semifield. Now P ⊂ A ⊂ S, so P ⊂ S is a semifield. Hence S is a S-semiring.

In view of this theorem it is not essential in the definition of S-subsemiring to mention that S is a S-semiring for the very definition of S-subsemiring forces S to be a S-semiring.

*Example 4.2.1*: Let $Z^o[x]$ be the semiring. $Z^o[x]$ has a proper S-subsemiring, for take A = {set of all polynomials of even degree with coefficients from $Z^o$}. Clearly A is a S-subsemiring as $Z^o \subset A$ is a semifield.

All S-semirings need not have S-subsemirings. Further all subsemirings of S need not be S-subsemirings of S. We illustrate these by the following examples.

*Example 4.2.2*: Let $Z^o$ be the semiring. $Z^o$ has subsemirings viz. pz = {0, p, 2p, …} for every integer p. But none of these subsemirings are S-subsemirings, also $Z^o$ is not a S-semiring.

Thus this example asserts if S is not a S-semiring. S has subsemirings but none of them are S-subsemirings.

*Example 4.2.3*: Let S be a semiring given by the following lattice whose Hasse diagram is:

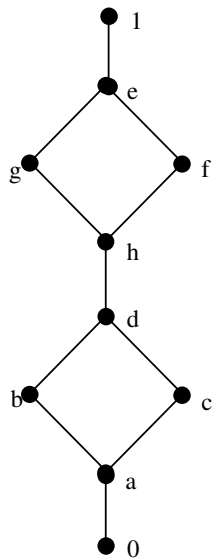

**Figure 4.2.1**

Now S is a S-semiring for $[0,1] = A \subset S$ is a semifield of S. Now S has several subsemiring but all of them are not S-subsemirings. For take the set A = {d, b, c, a, 0} in S. A is a subsemiring of S but A is not a S-subsemiring of S as A has no subfields contained in them.

In view of this we have the following theorem.

**THEOREM 4.2.2**: *Every subsemiring of a S-semiring need not in general be a S-subsemiring.*

*Proof*: To prove this consider the example 4.2.3 Clearly S is a S-semiring. Take A = {0, a}, A is a subsemiring of S but A is not a S-subsemiring of S. Hence the claim.

*Example 4.2.4*: Now consider the semiring $M_{3 \times 3} = \{(a_{ij})/a_{ij} \in Z^o\}$. $M_{3 \times 3}$ is a S-semiring. Let $A = \left\{ \begin{pmatrix} a_{11} & 0 & 0 \\ 0 & a_{22} & 0 \\ 0 & 0 & 0 \end{pmatrix} \middle/ a_{11}, a_{22} \in Z^o \right\}$. A is a subsemiring of $M_{3 \times 3}$. A is a S-subsemiring of $M_{3 \times 3}$ for take $P = \left\{ \begin{pmatrix} a_{11} & 0 & 0 \\ 0 & 0 & 0 \\ 0 & 0 & 0 \end{pmatrix} \middle/ a_{11} \in Z^o \right\}$, P is a semifield and $P \subset A$. So A is a S-semiring. The unit of the semifield is $\begin{pmatrix} 1 & 0 & 0 \\ 0 & 0 & 0 \\ 0 & 0 & 0 \end{pmatrix}$ which is not the unit of $M_{3 \times 3}$.

**DEFINITION 4.2.2**: *Let S be a S-semiring. We say S is a Smarandache commutative semiring (S-commutative semiring) if S has a S-subsemiring which is commutative. If the S-semiring has no commutative S-subsemiring then we say S is a Smarandache non-commutative semiring (S-non commutative semiring).*

The following facts are interesting about such S-semirings. If S is a semiring, which is commutative then trivially, S is a S-commutative semiring provided S is a S-semiring.

The important factor to be observed is even if S is not a commutative semiring still S can be a S-commutative semiring, which is evident from the following example.

*Example 4.2.5*: The semiring $M_{3 \times 3}$ described in example 4.2.4 is a non-commutative semiring which is a S-semiring. But $M_{3 \times 3}$ is a S-commutative semiring as the S-subsemiring $A = \left\{ \begin{pmatrix} a_{11} & 0 & 0 \\ 0 & a_{22} & 0 \\ 0 & 0 & 0 \end{pmatrix} \middle/ a_{11}, a_{22} \in Z^o \right\}$ is a commutative subsemiring. Hence the claim.

*Example 4.2.6*: A S-commutative semiring need not be a commutative semiring.

*Proof*: By an example. In the example 4.2.5 $M_{3\times 3}$ is a non-commutative semiring but is a S-commutative semiring.

**DEFINITION 4.2.3**: *Let S be a semiring. A non-empty subset P of S is said to be a Smarandache right (left) ideal (S-right (left) ideal) of S if the following conditions are satisfied.*

1. *P is a S-subsemiring.*
2. *For every $p \in P$ and $A \subset P$ where A is the semifield of P we have for all $a \in A$ and $p \in P$, ap (pa) is in A.*

*If P is simultaneously both a S-right ideal and a S-left ideal then we say P is a Smarandache ideal (S-ideal) of S.*

**Example 4.2.7**: Let $M_{3\times 3}$ be the semiring given in example 4.2.5. Clearly P = $\left\{ \begin{pmatrix} a_{11} & 0 & 0 \\ 0 & a_{22} & 0 \\ 0 & 0 & 0 \end{pmatrix} \middle/ a_{11}, a_{22} \in Z^o \right\}$ is a S-ideal of $M_{3\times 3}$. It is easily verified that P is both a S-left ideal and a right S-ideal of $M_{3\times 3}$.

**Example 4.2.8**: Let $M_{2\times 2} = \left\{ \begin{pmatrix} a & b \\ c & d \end{pmatrix} \middle/ a, b, c, d \in C_2 = (0,1) \right\}$ = set of all $2 \times 2$ matrices with entries from the chain lattice $C_2$.

$$M_{2\times 2} = \left\{ \begin{matrix} \begin{pmatrix} 0 & 0 \\ 0 & 0 \end{pmatrix}, \begin{pmatrix} 0 & 1 \\ 0 & 0 \end{pmatrix}, \begin{pmatrix} 1 & 0 \\ 0 & 0 \end{pmatrix}, \begin{pmatrix} 1 & 1 \\ 0 & 0 \end{pmatrix}, \\ \begin{pmatrix} 0 & 1 \\ 0 & 1 \end{pmatrix}, \begin{pmatrix} 0 & 0 \\ 1 & 1 \end{pmatrix}, \begin{pmatrix} 1 & 0 \\ 0 & 1 \end{pmatrix}, \begin{pmatrix} 0 & 1 \\ 1 & 0 \end{pmatrix}, \\ \begin{pmatrix} 1 & 1 \\ 1 & 0 \end{pmatrix}, \begin{pmatrix} 0 & 1 \\ 1 & 1 \end{pmatrix}, \begin{pmatrix} 1 & 0 \\ 1 & 1 \end{pmatrix}, \begin{pmatrix} 1 & 1 \\ 0 & 0 \end{pmatrix}, \\ \begin{pmatrix} 0 & 1 \\ 1 & 1 \end{pmatrix}, \begin{pmatrix} 1 & 0 \\ 1 & 1 \end{pmatrix}, \begin{pmatrix} 1 & 1 \\ 0 & 1 \end{pmatrix}, \begin{pmatrix} 1 & 1 \\ 1 & 1 \end{pmatrix} \end{matrix} \right\}$$

Clearly $M_{2\times 2}$ is a S-semiring with 16 elements in it.

For A = $\left\{ \begin{pmatrix} 1 & 0 \\ 0 & 1 \end{pmatrix}, \begin{pmatrix} 0 & 0 \\ 0 & 0 \end{pmatrix} \right\}$ is a semifield. To find ideals in $M_{2\times 2}$. The set S = $\left\{ \begin{pmatrix} 0 & 0 \\ 0 & 0 \end{pmatrix}, \begin{pmatrix} 1 & 0 \\ 0 & 0 \end{pmatrix} \right\}$ is a subsemiring of $M_{2\times 2}$. Clearly S is not a S-subsemiring so S cannot be a S-ideal of $M_{2\times 2}$. $A = \left\{ \begin{pmatrix} 1 & 0 \\ 0 & 0 \end{pmatrix}, \begin{pmatrix} 1 & 0 \\ 0 & 1 \end{pmatrix}, \begin{pmatrix} 0 & 0 \\ 0 & 0 \end{pmatrix} \right\}$ is a S-subsemiring of $M_{2\times 2}$, which is not a S-ideal of $M_{2\times 2}$.

From this example we derive the following result.

**THEOREM 4.2.3**: *Let S be a S-semiring. Every S-ideal of S is a S-subsemiring of S but every S-subsemiring of S in general need not be a S-ideal of S.*

*Proof*: Clearly by the very definition of the S-ideal we see every S-ideal is a S-subsemiring. Conversely to show every S-subsemiring in general need not be an S-ideal of S. We give the following example. Clearly the subset $A = \left\{ \begin{pmatrix} 1 & 0 \\ 0 & 0 \end{pmatrix}, \begin{pmatrix} 1 & 0 \\ 0 & 1 \end{pmatrix}, \begin{pmatrix} 0 & 0 \\ 0 & 0 \end{pmatrix} \right\}$ in $M_{2\times 2}$ given in example 4.2.8 is a S-subsemiring which is not an S-ideal of $M_{2\times 2}$. Hence the claim.

**DEFINITION 4.2.4**: *Let S be a semiring. A non-empty proper subset A of S is said to be Smarandache pseudo-subsemiring (S-pseudo subsemiring) if the following condition is true.*

*If there exists a subset P of S such that $A \subset P$; where P is a S-subsemiring i.e. P has a subset B such that B is a semifield under the operations of S or P itself is a semifield under the operations of S.*

*Example 4.2.9*: Let S be the semiring which is a lattice having the following Hasse diagram:

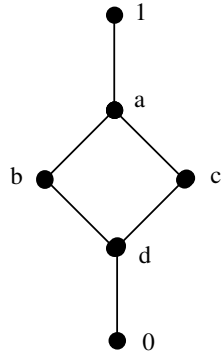

**Figure 4.2.2**

S is S-semiring. Let $A = \{0, a, b, d\}$. Now A is a S-pseudo subsemiring for $A \subset B = \{1, a, b, d, 0\}$. Hence the claim. We see all subsets need not be S-pseudo subsemirings. Let $A = \{0, a, d, b, c\}$. Clearly A cannot be contained in a proper semifield. So the set $A = \{0, a, b, c, d\}$ is not a S-pseudo subsemiring.

**THEOREM 4.2.4**: *Let S be a semiring every proper subset of S need not in general be a S-pseudo subsemiring.*

*Proof*: By an example. Consider the semiring given in example 4.2.9. The set {0, a, b, c, d} is not a S-pseudo subsemiring.

**THEOREM 4.2.5**: *Let S be a semiring if S contains a S-pseudo subsemiring then S is a S-semiring.*

*Proof*: S is a semiring. S contains a S-pseudo subsemiring A; i.e. A is contained in a semifield P, P contained in S. So S is a S-semiring.

The concept of S-pseudo subsemiring leads to the definition of S-pseudo ideals in the semiring.

**DEFINITION 4.2.5**: *Let S be a semiring. A non-empty subset P of S is said to be a Smarandache pseudo right(left) ideal (S-pseudo right (left) ideal) of the semiring S if the following conditions are true.*

1. *P is a S-pseudo subsemiring i.e. P ⊂ A, A a semifield in S.*
2. *For every p ∈ P and every a ∈ A, ap ∈ P (pa ∈ P).*

*If P is simultaneously both a S-pseudo right ideal and S-pseudo left ideal we say P is a Smarandache pseudo ideal (S-pseudo ideal).*

Now we define two new notions about S-semiring viz. Smarandache dual ideal and Smarandache pseudo dual ideal of a semiring S.

**DEFINITION 4.2.6**: *Let S be a semiring. A non-empty subset P of S is said to be a Smarandache dual ideal of S (S-dual ideal) if the following conditions hold good.*

1. *P is a S-subsemiring*
2. *For every p ∈ P and a ∈ A\{0} a + p is in A, where A ⊂ P.*

*Example 4.2.10*: Let S be the semiring given by the Hasse diagram:

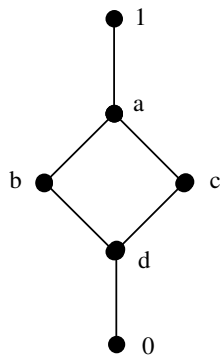

**Figure 4.2.3**

S is a S-semiring. Let P = {0, a, b, d, 1}. P is a S-dual ideal of S. For A = {0, 1} is a semifield.

***Example 4.2.11***: Let $M_{2\times 2} = \left\{ \begin{pmatrix} a & b \\ c & d \end{pmatrix} \middle/ a,b,c,d \in C_2, \text{the two element chain lattice} \right\}$, $M_{2\times 2}$ is a semiring. Take $P = \left\{ \begin{pmatrix} a & 0 \\ 0 & b \end{pmatrix} \middle/ a,b \in C_2 \right\}$. P is a S-dual ideal of $M_{2\times 2}$. For A $= \left\{ \begin{pmatrix} 1 & 0 \\ 0 & 1 \end{pmatrix}, \begin{pmatrix} 0 & 0 \\ 0 & 0 \end{pmatrix} \right\}$ is a semifield.

***Example 4.2.12***: Let $Z^o[x]$ be the S-semiring. $P = \{Z_2^o[x] =$ all polynomials of even degree$\}$.

P is a S-subsemiring for $Z^o \subset Z_2^o[x]$ so $Z^o$ is a semifield. Thus P is not a S-dual ideal. If we take $P = Z_2^o[x]$. Clearly P is not a S-ideal but P is a S-subsemiring, for $Z^o \subset Z_2^o[x]$. If $p(x) \in Z_2^o[x]$, $p(x)Z^o \not\subset Z^o$. Hence the claim. This leads to the following.

**THEOREM 4.2.6**: *Let S be a S-semiring. If P is a S-dual ideal then P need not be a S-ideal.*

*Proof*: By example 4.2.11 the result in the theorem is true.

Now we define S-pseudo dual ideal.

**DEFINITION 4.2.7**: *Let S be a semiring. A non-empty subset P of S is said to be a Smarandache pseudo dual ideal (S-pseudo dual ideal) of S if the following conditions are true.*

1. *P is a S-pseudo subsemiring i.e. $P \subset A$, A is a semifield in S or A contains a semifield.*

2. *For every $p \in P$ and $a \in A$, $p + a \in P$. Clearly P is simultaneously left and right S-pseudo dual ideal of S as S is additively commutative.*

The reader is assigned the work of finding examples. Of course, several examples exist. Here we introduce yet another nice substructure in a semiring called a Smarandache semidivision ring.

**DEFINITION 4.2.8**: *Let S be a S-semiring. S is said to be a Smarandache semidivision ring (S-semidivision ring) if the proper subset $A \subset S$ is such that*

1. *A is a Smarandache subsemiring which is non-commutative.*
2. *A contains a subset P such that P is a semidivision ring that is P has no zero divisors and P is a non-commutative semiring.*

<u>*Remark*</u>: The concept of S-semidivision ring can be defined only when the semiring under consideration is non-commutative.

***Example 4.2.13***: Let $M_{2\times 2} = \{(a_{ij})/a_{ij} \in Z^o\}$. Let $A = \left\{ \begin{pmatrix} x & y \\ 0 & z \end{pmatrix} \bigg/ x, y, z \in Z^o \right\}$. A is a S-subsemiring as the set $B = \left\{ \begin{pmatrix} x & 0 \\ 0 & 0 \end{pmatrix} \bigg/ x \in Z^o \right\}$ is a semifield. If we take $P = \left\{ \begin{pmatrix} x & y \\ 0 & z \end{pmatrix} \bigg/ x, y, z \in Z^o \setminus \{0\} \right\} \cup \left\{ \begin{pmatrix} 0 & 0 \\ 0 & 0 \end{pmatrix} \right\}$ then P is a semidivision ring. Hence $M_{2\times 2}$ is a S-semidivision ring.

**THEOREM 4.2.8**: *Let S be a S-semidivision ring. Then we have two properties to be true.*

1. *S is a S-semiring*
2. *S has a S-subsemiring*

*Proof*: By the very definition of S-semidivision ring, the results are true.

The natural question would be if a semiring S has S-subsemiring A, then can we say A is a S-semidivision ring. The answer is, all S-subsemiring are not semidivision rings so A need not in general be a S-semidivision ring. This is illustrated by the following:

***Example 4.2.14***: Let $M_{2\times 2}$ be as in the above example we see $S_{2\times 2} = \left\{ \begin{pmatrix} a & b \\ 0 & 0 \end{pmatrix} \bigg/ a, b \in Z^o \right\}$ is a S-subsemiring but $S_{2\times 2}$ does not contain a semidivision ring.

**PROBLEMS**:

1. Find S-subsemirings and S-ideals in $S = Z^o \times Z^o \times Z^o$?

2. Find S-subsemirings in the semiring S which are not ideals (S given in Problem 1).

3. Can $Z^o$ have S–subsemirings which are not S-ideals?

4. Does $Q^o$ have S-ideals?

5. Find S-subsemirings in $R^o$. Does $R^o$ have S-ideals?

6. Find S-ideals and S-subsemirings in $M_{2\times 2} = \left\{ \begin{pmatrix} a & b \\ c & s \end{pmatrix} \bigg/ a, b, c, d \in Z^o \right\}$ where $M_{2\times 2}$ is a semiring.

7. Find S-subsemirings in $M_{2\times 2}$ given in Problem 6 which are not S-ideals.

8. Can $M_{2\times 2}$ given in Problem 6 have S-right ideals? If so find them.

9. Find the order of the S-semiring $M_{4\times4} = \{(a_{ij}) / a_{ij} \in C_2; C_2$ is the chain lattice with two elements 0 and 1$\}$.

   Find
       a. S-ideals of $M_{4\times4}$
       b. S-subsemirings which are not S-ideals
       c. S-right ideals which are not S-left ideals
       d. Subsemirings which are not S-subsemirings.

10. Find a semiring of finite order n; n > 2 which is not a S-semiring.

11. Find S-pseudo ideal of $M_{2\times2}$ given in problem 6.

12. Find S-dual ideal in the semiring $M_{4\times4}$ given in problem 9.

13. Find S-pseudo dual ideals in the semiring $M_{2\times2}$ of problem 6 and $M_{4\times4}$ of problem 9.

14. Give an example of a S-pseudo dual ideal which is not a S-pseudo ideal.

15. Give an example of a S-pseudo ideal which is not an ideal.

16. Give an example of S-pseudo dual ideal which is not an ideal.

17. Find a semiring which has i) S-ideals, ii) S-pseudo ideals, iii) S-dual ideals, and iv) S-pseudo dual ideals.

18. Find in the semiring $M_{4\times4}$ given in example 9 a S-semidivision ring.

19. Find in $M_{4\times4}$ given in example 9 a S-subsemiring which has no S-semidivision rings.

## 4.3 Smarandache Special Elements in S-Semirings

In this section we introduce the concept of Smarandache idempotents, Smarandache zero divisors, Smarandache units and Smarandache inverses in a semiring and illustrate them by examples and study them.

**DEFINITION 4.3.1**: *Let S be any semiring. We say a and b $\in$ S is a Smarandache zero divisor (S-zero divisor) if a $\bullet$ b = 0 and there exists x, y $\in$ S\\{a, b, 0\}, x $\neq$ y with*

    *1. ax = 0 or xa = 0*
    *2. by = 0 or yb = 0 and*
    *3. xy $\neq$ 0 or yx $\neq$ 0*

*Clearly if S is a semifield we will not have S-zero divisors.*

**Example 4.3.1**: Let $S = Z^o \times Z^o \times Z^o \times Z^o$ be a semiring. Clearly S is a S-semiring as $P = \{(x, 0, 0, 0) / x \in Z_o\}$ is a semifield, $P \subset Z^o \times Z^o \times Z^o \times Z^o$. Now a = (0, 0, 4, 2) and

b = (5, 0, 0, 0) in S is such that a • b = (0, 0, 0, 0) i.e. they are zero divisors. Take x = (2, 8, 0, 0) and y = (0, 1, 0, 0) ∈ S\ {a, b, (0, 0, 0, 0)}. Clearly xa = 0 and ax = 0, yb = 0 and by = 0, but xy = (0, 8, 0, 0) ≠ (0, 0, 0, 0). So a = (0, 0, 4, 2) is a S-zero divisor in S.

***Example 4.3.2***: Let $M_{2\times 2} = \left\{ \begin{pmatrix} a & b \\ c & d \end{pmatrix} \Big/ a,b,c,d \in Z^o \right\}$ be the set of all $2 \times 2$ matrices with entries from $Z^o$. $M_{2\times 2}$ is a semiring. Now consider $A = \begin{pmatrix} 1 & 0 \\ 0 & 0 \end{pmatrix}$ and $B = \begin{pmatrix} 0 & 0 \\ 1 & 0 \end{pmatrix}$. Clearly

$$A \bullet B = \begin{pmatrix} 1 & 0 \\ 0 & 0 \end{pmatrix} \begin{pmatrix} 0 & 0 \\ 1 & 0 \end{pmatrix} = \begin{pmatrix} 0 & 0 \\ 0 & 0 \end{pmatrix}$$

Choose $X = \begin{pmatrix} 0 & 0 \\ 0 & 1 \end{pmatrix}$ and $Y = \begin{pmatrix} 0 & 0 \\ 1 & 1 \end{pmatrix}$ from $M_{2\times 2}$.

$$AX = \begin{pmatrix} 1 & 0 \\ 0 & 0 \end{pmatrix} \begin{pmatrix} 0 & 0 \\ 0 & 1 \end{pmatrix} = \begin{pmatrix} 0 & 0 \\ 0 & 0 \end{pmatrix}$$

$$XA = \begin{pmatrix} 0 & 0 \\ 0 & 1 \end{pmatrix} \begin{pmatrix} 1 & 0 \\ 0 & 0 \end{pmatrix} = \begin{pmatrix} 0 & 0 \\ 0 & 0 \end{pmatrix}$$

$$Y \bullet B = \begin{pmatrix} 0 & 0 \\ 1 & 1 \end{pmatrix} \begin{pmatrix} 0 & 0 \\ 1 & 0 \end{pmatrix} = \begin{pmatrix} 0 & 0 \\ 1 & 0 \end{pmatrix}$$

$$B \bullet Y = \begin{pmatrix} 0 & 0 \\ 1 & 0 \end{pmatrix} \begin{pmatrix} 0 & 0 \\ 1 & 1 \end{pmatrix} = \begin{pmatrix} 0 & 0 \\ 0 & 0 \end{pmatrix}$$

$$Y \bullet B \neq \begin{pmatrix} 0 & 0 \\ 1 & 0 \end{pmatrix}$$

But $XY = \begin{pmatrix} 0 & 0 \\ 0 & 1 \end{pmatrix} \begin{pmatrix} 0 & 0 \\ 1 & 1 \end{pmatrix} \neq \begin{pmatrix} 0 & 0 \\ 1 & 1 \end{pmatrix}$. $YX = \begin{pmatrix} 0 & 0 \\ 1 & 1 \end{pmatrix} \begin{pmatrix} 0 & 0 \\ 0 & 1 \end{pmatrix} \neq \begin{pmatrix} 0 & 0 \\ 0 & 1 \end{pmatrix}$.

Thus A, B ∈ $M_{2\times 2}$ are S-zero divisors of $M_{2\times 2}$. Now consider, if we take A = B = $\begin{pmatrix} 0 & 0 \\ 1 & 0 \end{pmatrix}$. Clearly $A^2 = \begin{pmatrix} 0 & 0 \\ 1 & 0 \end{pmatrix} \begin{pmatrix} 0 & 0 \\ 1 & 0 \end{pmatrix} = \begin{pmatrix} 0 & 0 \\ 0 & 0 \end{pmatrix}$. We will find out whether A is a S-zero divisor of the semiring $M_{2\times 2}$. Take $X = \begin{pmatrix} 0 & 0 \\ 1 & 1 \end{pmatrix}$ then we have $AX = \begin{pmatrix} 0 & 0 \\ 0 & 0 \end{pmatrix}$, XA

$= \begin{pmatrix} 0 & 0 \\ 1 & 0 \end{pmatrix} \neq \begin{pmatrix} 0 & 0 \\ 0 & 0 \end{pmatrix}$ and $X^2 = \begin{pmatrix} 0 & 0 \\ 1 & 1 \end{pmatrix}\begin{pmatrix} 0 & 0 \\ 1 & 1 \end{pmatrix} = \begin{pmatrix} 0 & 0 \\ 1 & 1 \end{pmatrix} \neq \begin{pmatrix} 0 & 0 \\ 0 & 0 \end{pmatrix}$. Thus A is also a S-zero divisor of $M_{2 \times 2}$.

**THEOREM 4.3.1**: *Let S be a semiring if $x \in S$ is a S-zero divisor then x is a zero divisor.*

*Proof*: By the very definition of the S-zero divisor we see every S-zero divisor is a zero divisor.

The converse of this theorem remains as an open problem. We define a concept called Smarandache anti zero divisor in a semiring S.

**DEFINITION 4.3.2**: *Let S be a semiring. An element $x \in S$ is said to be a Smarandache anti-zero divisor (S-anti zero divisor) if we can find a y such that $xy \neq 0$ and a, b $\in$ S\\{0, x, y} such that*

1. *$ax \neq 0$ or $xa \neq 0$*
2. *$by \neq 0$ or $yb \neq 0$*
3. *$ab = 0$ or $ba = 0$.*

*Example 4.3.3*: Let $M_{2 \times 2} = \left\{ \begin{pmatrix} a & b \\ c & d \end{pmatrix} \middle/ a, b, c, d \in Z^o \right\}$ be the semiring. $A = \begin{pmatrix} 1 & 0 \\ 0 & 1 \end{pmatrix} \in M_{2 \times 2}$ is an anti-zero divisor. For take $B = \begin{pmatrix} 0 & 1 \\ 1 & 0 \end{pmatrix}$ is such that $AB \neq (0)$. Choose

$X = \begin{pmatrix} 0 & 0 \\ 1 & 0 \end{pmatrix}$ and $Y = \begin{pmatrix} 0 & 0 \\ 0 & 1 \end{pmatrix}$ in $M_{2 \times 2}$. $AB \neq \begin{pmatrix} 0 & 0 \\ 0 & 0 \end{pmatrix}$, $AX \neq \begin{pmatrix} 0 & 0 \\ 0 & 0 \end{pmatrix}$, $BY \neq \begin{pmatrix} 0 & 0 \\ 0 & 0 \end{pmatrix}$.

But $XY = \begin{pmatrix} 0 & 0 \\ 1 & 0 \end{pmatrix}\begin{pmatrix} 0 & 0 \\ 0 & 1 \end{pmatrix} = \begin{pmatrix} 0 & 0 \\ 0 & 0 \end{pmatrix}$. Thus A is a S-anti-zero divisor.

This leads us to the following theorem which is very unique in its own way.

**THEOREM 4.3.2**: *Let S be a semiring. If $x \in S$ is a S-anti zero divisor then x need not in general be a zero divisor.*

*Proof*: We see from the above example 4.3.3. The element $\begin{pmatrix} 1 & 0 \\ 0 & 1 \end{pmatrix}$ is the unit element of $M_{2 \times 2}$ so for no non-zero element in $M_{2 \times 2}$ $\begin{pmatrix} 1 & 0 \\ 0 & 1 \end{pmatrix}$ is a zero divisor. But $\begin{pmatrix} 1 & 0 \\ 0 & 1 \end{pmatrix}$ is a S-anti zero divisor of $M_{2 \times 2}$. Hence the claim.

**THEOREM 4.3.3**: *Let S be a semiring. If S has a S-anti zero divisor then S has a non-trivial zero divisor.*

*Proof*: From the very definition of the S-anti zero divisor in a semiring we are guaranteed of a non-trivial zero divisor in S.

*Example 4.3.4*: Let $S = Z^o \times Z^o \times Z^o \times Z^o \times Z^o \times Z^o \times Z^o$ (7 times). Clearly S is a semiring with (1, 1, 1, 1, 1, 1, 1) as a unit and S has non-trivial zero divisors. To show (1, 1, 1, 1, 1, 1, 1) is an S-anti zero divisor. Consider any x in S\{0}. x(1,1,1,1,1,1,1) = x. Let x = (0, 0, 2, 3, 1, 0, 0). Now let a = (2, 1, 0, 0, 6, 0, 0) and b = (0, 0, 1, 2, 0, 3, 4). Clearly a(1, 1, 1, 1, 1, 1, 1) ≠ 0, bx ≠ 0 but ab = 0. So (1, 1, 1, 1, 1, 1, 1) is an S-anti zero divisor of S.

We have the following result in view of these examples:

**THEOREM 4.3.4**: *If S is a semiring with unit and if S has zero divisors then the unit is an S-anti zero divisor.*

*Proof*: Given S has zero divisors say $\alpha, \beta \in S\setminus\{0\}$ with $\alpha\beta = 0$. Let 1 be the identity of S. Choose any x, we have $1 \bullet x \neq 0$ such that $\alpha \bullet 1 \neq 0$, $\beta x \neq 0$ but $\alpha\beta = 0$. Hence the claim.

One may come to think only the unit of a semiring is a S-anti zero divisor. It is not so we have other elements to be also S-anti zero divisors, which is evident from the following example:

*Example 4.3.5*: Let $S = Z_o \times Z_o \times Z_o \times Z_o \times Z_o$ be the semiring. Clearly x = (1, 1, 1, 1, 0) ∈ S. $x \bullet y_1 = 0$ for $y_1 = (0, 0, 0, 0, 4)$. But for y = (0, 0, 6, 7, 0) we have xy = (0, 0, 6, 7, 0). Choose a = (3, 2, 0, 0, 0) and b = (0, 0, 0, 9, 2), ax = (3, 2, 0, 0, 0), by = (0, 0, 0, 6, 3, 0). But $a \bullet b = (0, 0, 0, 0, 0)$. Thus (1, 1, 1, 1, 0) is a S-anti zero divisor in this case the two things are to be observed:

        a.  (1, 1, 1, 1, 0) is a zero divisor.
        b.  (1, 1, 1, 1, 0) is not a unit.

Now we proceed on to define Smarandache idempotents in a semiring S.

**DEFINITION 4.3.3**: *Let S be a semiring. An element $0 \neq a \in S$ is a Smarandache idempotent (S-idempotent) of S if*

    i)    $a^2 = a$
    ii)   *There exists $b \in S \setminus \{x\}$ such that i) $b^2 = a$ and ii) ab = b (ba = b) or ba = a (ab = a)*

*'or' in condition 2(ii) is mutually exclusive.*

*Example 4.3.6*: Let $C_2$ be the chain lattice. $S_3$ be the symmetric group of degree 3. $C_2S_3$ be the group semiring of the group $S_3$ over the semiring $C_2$. Now $1 + p_4 + p_5 \in S_3$ where

$$1 = \begin{pmatrix} 1 & 2 & 3 \\ 1 & 2 & 3 \end{pmatrix}, \quad p_1 = \begin{pmatrix} 1 & 2 & 3 \\ 1 & 3 & 2 \end{pmatrix}, \quad p_2 = \begin{pmatrix} 1 & 2 & 3 \\ 3 & 2 & 1 \end{pmatrix},$$

$$p_3 = \begin{pmatrix} 1 & 2 & 3 \\ 2 & 1 & 3 \end{pmatrix}, \quad p_4 = \begin{pmatrix} 1 & 2 & 3 \\ 2 & 3 & 1 \end{pmatrix}, \quad p_5 = \begin{pmatrix} 1 & 2 & 3 \\ 3 & 1 & 2 \end{pmatrix}$$

$(1 + p_4 + p_5)^2 = 1 + p_4 + p_5$. Take $p_1 + p_2 + p_3 \in S_3$ is such that $(p_1 + p_2 + p_3)^2 = 1 + p_4 + p_5$. $(1 + p_4 + p_5)(p_1 + p_2 + p_3) = p_1 + p_2 + p_3$. Thus $1 + p_4 + p_5$ is a S-idempotent of $C_2S_3$.

*Example 4.3.7*: Let $S = Z^o \times Z^o \times Z^o \times Z^o \times Z^o$ (5 times). S is a semiring. S has non-trivial idempotents which are other than (1, 1, 1, 1, 1) viz. (1, 0, 0, 0, 0), (0, 1, 0, 0, 0), (0, 0, 1, 0, 0) …. (1, 1, 1, 1, 0), (1, 1, 0, 1, 1) and (0, 1, 1, 1, 1). But none of them are S-idempotents.

Thus we have the following theorem in view of this example.

**THEOREM 4.3.5**: *Let S be a semiring. Every idempotent in general need not be a S-idempotent.*

*Proof*: In the semiring S in example 4.3.7 we see the idempotents in $S = Z^o \times Z^o \times Z^o \times Z^o \times Z^o$ are not S-idempotents.

Next we proceed on to define S-units in a semiring.

**DEFINITION 4.3.4**: *Let S be a semiring with identity 1. We say $x \in S\setminus\{1\}$ to be a Smarandache unit (S-unit) if there exists a $y \in S$ such that*

1. *xy = 1*
2. *There exists a, b $\in S \setminus \{x, y, 1\}$ such that (i) xa = y or ax = y, or (ii) yb = x or by = x and ab = 1.*

We leave it as an exercise for a reader to construct an example of a S-unit. We see as in the case of rings, group rings or semigroup rings which imply the existence of zero divisors if the ring has idempotents and alike results, can never be proved in case of semirings. For in a ring if e is a non-trivial idempotent (e ≠ 0 and e ≠ 1) such that $e^2 = e$ implies $e^2 - e = 0$ so $e(e - 1) = 0$ is a zero divisor. But we cannot arrive at a zero divisor using the existence of an idempotent in case of semirings.

Thus semirings have very unique and special properties enjoyed by them.

**PROBLEMS**:

1. Find S-units if any in the S-semiring, $M_{5 \times 5} = \{(a_{ij})/a_{ij} \in C_2\}$.

2. Does $M_{5 \times 5}$ in Problem 1 have S-idempotents? S – zero divisors?

3. Can P(X) where X = (1, 2, 3, 4, 5, 6, 7) have

        i.    S-idempotents?
        ii.   S-units?
        iii.  S-zero divisors?

Justify your answer for the existence or non-existence.

4. Find S – zero divisors if it exists in the S-semiring, $M_{2\times 2} = \{(a_{ij})/ a_{ij} \in Z^o\}$.. Can $M_{2\times 2}$ have S-idempotents?

5. Give an example of a S-semiring which has no S-units but S-zero divisor and S-idempotents.

6. Can a S-semiring have only S-zero divisors and no S-units?

7. Find any relation between the existence or non-existence of S- units, S-idempotents or S-zero divisors.

8. Find all S-idempotents in $S = Z^o \times Z^o \times Z^o$.

9. How many S-idempotents does $S = Z^o \times Z^o \times Z^o \times Z^o \times Z^o$ contain?

10. Find S-idempotents in $C_2S(5)$.

## 4.4 Special S-semirings

In this section we introduce some special types of Smarandache semirings like Smarandache compact semiring, Smarandache e-semiring, Smarandache congruence simple semirings, Smarandache ∗- semirings, Smarandache inductive ∗ semirings and Smarandache chain semirings. We define and give examples of them wherever possible. Apart from this the reader is left to develop these concepts. Further we make a special assumption that when we say semirings we mean only semirings, which are not rings, as rings or fields are trivially semirings.

**DEFINITION 4.4.1**: *Let S be a semiring, we say S has a Smarandache congruence relation (S-congruence relation) ~ if we have a S-subsemiring A of S such that '~' is an equivalence relation that also satisfies*

$$x_1 \sim x_2 \Rightarrow \begin{cases} c + x_1 \sim c + x_2 \\ x_1 + c \sim x_2 + c \\ cx_1 \sim cx_2 \\ x_1 c \sim x_2 c \end{cases}$$

*for all $x_1, x_2, c \in P \subset A$ where P is a semifield under the operations of S. A semiring S that admits no Smarandache congruence relation other than the trivial ones, identity and $P \times P$, is that semiring which has no proper S- subsemiring on which the S-congruence relation can be defined, is called the Smarandache congruence simple semiring ( S-c-simple semiring).*

**DEFINITION 4.4.2**: *Let S be Smarandache semigroup that is S is a semigroup which has a proper subset A $\subseteq$ S and A is a group under the operations of S. We define for any S-semigroup V(S) = S $\cup$ {$\infty$}. Extend the multiplication in S to V(S) by the rule x$\infty$ = $\infty$x = $\infty$ for all x $\in$ V(S). Define an addition on V(S) by x + x = x and x + y = $\infty$ for all x, y $\in$ V(S) with x $\neq$ y.*

We define this analogous to V(G) for any abelian group G defined by Chris Monico to define S-c-simple semiring.

**THEOREM 4.4.1**: *Let S be a finite semigroup. If S is S-commutative semigroup then V(S) is a S-c-simple semiring.*

Proof: V(S) is S-c-simple whenever S is a Smarandache commutative semiring, that is S is a S-semiring and has at least a proper subset G $\subset$ S where G is a commutative group under the operations of S. Hence the claim, as V(G) when G is a finite commutative group is a S-c-simple semiring.

*Example 4.4.1*: Let S(4) be the semigroup of mappings of the set ($x_1$, $x_2$, $x_3$, $x_4$) to itself. V(S(4)) is a S-c-semiring. For the group G generated by $\begin{pmatrix} x_1 & x_2 & x_3 & x_4 \\ x_2 & x_3 & x_4 & x_1 \end{pmatrix}$ is cyclic of degree 4. Hence the claim.

*Example 4.4.2*: Let S(n) be a multiplicative semigroup of maps from a set of n elements to itself. V(S(n)) will be a S-c-simple semiring. For V(S(n)) contains a group G generated by $x = \begin{pmatrix} 1 & 2 & 3 & \ldots & n-1 & n \\ 2 & 3 & 4 & \ldots & n & 1 \end{pmatrix}$ which is a cyclic group of degree n which is abelian. Now V(G) $\subset$ V(S(n)) so V(S(n)) is a S-c-simple semiring. In view of this we have the following theorem, which gives a class of S-commutative c-simple finite semirings.

**THEOREM 4.4.2**: *V(S(n)) is a S-commutative congruence simple semiring for all finite positive integer n.*

*Proof*: Since every S(n) is a S-semigroup and every S(n) has proper subsets which are abelian groups we have the result to be true.

**THEOREM 4.4.3**: *V(S(n)) has at least (n-1) proper subsets which are abelian groups. So V(S(n)) is a S-commutative congruence simple finite semiring.*

*Proof*: Obvious from the fact for every m, 1 < m < n we have cyclic group of order m, which we denote by $G_m$, so that V($G_m$) is a commutative, congruence simple finite semiring, hence V(S(n)) is a S-commutative congruence simple finite semiring.

Using the classical theorem of Cayley which states, "Every group is isomorphic to a subgroup of A(S) for some appropriate S". We can modify it and restate it as " Every finite abelian group is isomorphic to a subgroup of $S_n$ for some appropriate n". Now we use the analogue of this classical theorem for S-semigroups which is given by us as "Every S-semigroup is isomorphic to a S-subsemigroup in S(n)".

Using this we have the following very interesting result.

**THEOREM 4.4.4**: *Let G be a finite abelian group V(G) be isomorphic to S for some commutative congruence simple finite semiring. Then there exists a suitable n so that every V(G) has a isomorphic image in the S-c-simple semiring V(S(n)) for a suitable n.*

*Proof:* We know by classical Cayley's theorem. G is isomorphic a subgroup in $S_n$ for a suitable n. Now each $S_n \subset S(n)$ and we have S(n) is a S-semigroup.

Now $V(G) \subset V(S_n)$ s$\subset$ V(S(n)) so we have the theorem to be true.

*Example 4.4.3*: Let $C_n$ be a chain lattice with n-elements, $M_{t \times t}$ be the collection of all $t \times t$ matrices with entries from $C_n$. $M_{t \times t}$ is a S-finite congruence simple semiring. Thus we have a class of Smarandache finite congruence simple semiring for varying t.

**DEFINITION 4.4.3**: *Let S be a semiring we say S is a Smarandache right chain semiring (S-right chain semiring) if the S-right ideals of S are totally ordered by inclusion. Similarly we define Smarandache left chain semiring (S-left chain semiring). If all the S-ideals are totally ordered by inclusion we say S is a Smarandache chain semiring (S-chain semiring).*

*Example 4.4.4*: Let $C_9$ be a chain lattice, which is a semiring of order 9. Clearly $C_9$ is a S-chain semiring.

**THEOREM 4.4.5**: *All chain lattices are S-chain semirings.*

*Proof*: Left for the reader to verify.

**DEFINITION 4.4.4**: *Let S be a semiring. If $S_1 \subset S_2 \subset ...$ is a monotonic ascending chain of S-ideals $S_i$, and there exists a positive integer r such that $S_r = S_s$ for all $s \geq r$. then we say the semiring S satisfies the Smarandache ascending chain conditions (S-acc) for S-ideals in the semiring S.*

*Example 4.4.5*: If we take the semiring S to be any chain lattice S, then it satisfies S-acc.

**DEFINITION 4.4.5**: *Let S be any semiring. We say S satisfies Smarandache descending chain condition (S-dcc) on S-ideals $S_i$ if every strictly decreasing sequence of S-ideals $N_1 \supset N_2 \supset N_3 \supset ...$ in S is of finite length. The Smarandache min condition (S – mc) for S-ideals holds in S if given any set P of S-ideals in S, there is an ideal of P that does not properly contain any other ideal in the set P.*

**DEFINITION 4.4.6**: *Let S be a chain ring such that S-acc for ideals holds in S. The Smarandache maximum condition (S-MC) for S-ideals holds in S if for every non-empty set P of S-ideals in S contains a S-ideal not properly contained in any other S-ideal of the set P.*

**DEFINITION 4.4.7**: *Let S be a semiring we say S is a Smarandache compact semiring (S-compact semiring) if A ⊂ S where A is a S subsemiring of S is a compact semiring under the operations of S. It is important to note if S is to be a S-compact semiring it is not necessary that S is a compact semiring. If S has a S-subsemiring which is a compact subsemiring then it is sufficient.*

So we have the following.

**THEOREM 4.4.6**: *If S is a compact semiring and S is a S-semiring then S is a S-compact semiring provided S has a S-subsemiring.*

*Proof*: Obvious by the very definition of S-compact semiring.

*Example 4.4.6*: All chain semirings $C_n[x]$ are S-compact semirings. For $C_n \subset C_n[x]$ is a S-subsemiring of $C_n[x]$ which is a compact semiring as a . b = b and a + b = a.

*Example 4.4.7*: Let X be a set with n elements, P(X) the power set of X which is the Boolean algebra. Hence P(X) is a semiring. P(X) is a S-compact semiring for every chain in P(X) connecting X and $\phi$ is a S-subsemiring which is a compact semiring.

**DEFINITION 4.4.8**: *Let S be a semiring, S is said to be Smarandache ∗-semiring (S-∗-semiring) if S contains a proper subset A satisfying the following conditions:*

1. *A is a subsemiring*
2. *A is a S-subsemiring*
3. *A is a ∗-semiring*

*So if S is a ∗- semiring and if S has a S-subsemiring then obviously S is a S-∗-semiring.*

**DEFINITION 4.4.9**: *Let S be any semiring. We say S is a Smarandache inductive ∗-semiring (S-inductive ∗-semiring) if S contains a proper subset A such that the following conditions are true*

1. *A is a subsemiring of S.*
2. *A is a S-subsemiring of S*
3. *A is inductive ∗-semiring.*

**THEOREM 4.4.7**: *If S is a inductive ∗-semiring and S has a S-subsemiring then S is a Smarandache inductive ∗-semiring.*

*Proof*: By the very definition of S-inductive ∗-semirings

**DEFINITION 4.4.10**: *A semiring S is said to be a Smarandache continuous semiring (S-continuous) if a proper subsemiring A of S satisfies the following 2 conditions*

1. *A is a S- subsemiring.*
2. *A is a continuous semiring.*

**DEFINITION 4.4.11**: *Let S be a semiring. S is said to be a Smarandache idempotent semiring (S-idempotent semiring) if A proper subset P of S, which is a subsemiring of S satisfies the following conditions:*

1. *P is a S-subsemiring.*
2. *P is an idempotent semiring.*

*Example 4.4.8*: Let $C_7[x]$ be the polynomial semiring. $C_7[x]$ is a S-idempotent semiring for $C_7 \subset C_7[x]$ is

1. Subsemiring of $C_7[x]$.
2. $C_7$ is a S-subsemiring.
3. In $C_7$ we have $a + a = a$ for all $a \in C_7$.

**DEFINITION 4.4.12**: *Let S be a semiring. S is said to be a Smarandache e-semiring (S-e-semiring) if S contains a proper subset A satisfying the following conditions:*

1. *A is a subsemiring*
2. *A is a S-subsemiring*
3. *A is a e-semiring*

*Example 4.4.9*: Let $C_n[x]$ be a polynomial semiring over the semiring $C_n$ (the chain lattice with n elements). $C_n[x]$ is a S-e-semiring for $C_n \subset C_n[x]$ satisfies all conditions. $C_n$ is clearly a S-e-semiring.

**DEFINITION 4.4.13**: *Let S be any semiring. G be a Smarandache semigroup. Consider the semigroup semiring SG. We call SG the Smarandache group semiring.(S-group semiring)*

*Example 4.4.10*: Let $Z^o S(n)$ be the semigroup semiring. $Z^o S(n)$ is a Smarandache group semiring.

Then the natural question would be what is the definition of Smarandache semigroup semiring.

**DEFINITION 4.4.14**: *Let G be a group, G is a Smarandache anti group, that is G contains a proper subset which is a semigroup. The group semiring FG where F is any semiring is called the Smarandache semigroup semiring.(S-semigroup semiring)*

*Example 4.4.11*: Let $Q^+$ be the group under multiplication. Clearly $Z^+$ is a semigroup in $Q^+$. Consider $C_3$ a chain lattice, that is a semiring, $C_3 Q^+$ is a Smarandache semigroup semiring.

We are not able to construct any finite Smarandache semigroup semiring, as we do not have examples of finite group, which contain proper subsets which are semigroups. This is left as an open problem.

**PROBLEMS**:

1. Give an example of a finite S-c- simple semiring?

2. Find an example of a S-commutative c-simple semiring.

3. Can you make a cyclic group of prime order a S-c-simple semiring? Justify your answer.

4. Is the group semiring $Z^o S_n$ a S-chain semiring?

5. Give an example of a semigroup semiring which is a S-chain semiring.

6. Give an example of S-∗-semiring.

7. Give an example S- inductive ∗-semiring.

8. Will $Z^o S_5$ group semiring satisfy S-acc condition?

9. Can $Z^o G$ the group semiring where G is an infinite cyclic group satisfy S-dcc condition?

10. Give an example of a semiring, which satisfies both S-dcc and S-acc.

11. Give an example of a semiring, which satisfies S-dcc and S-mc. What is the relation between S-dcc and S-mc?

12. Give an example of a S-compact semiring.

13. Is $M_{3\times 3} = \{(a_{ij}) / a_{ij} \in Z^o\}$ a S-compact semiring? The operation on $M_{3\times 3}$ is the usual matrix addition and matrix multiplication.

14. Can $M_{5\times 5} = \{(a_{ij}) / a_{ij} \in C_9\}$, $C_9$ the chain lattice be a S-compact semiring?

15. Is the semiring given in example 14 a S-∗ semiring?

16. Can the semiring given in example 13 be a S-inductive ∗ semiring?

17. Give an example of a S-group semiring.

18. Give an example of a Smarandache semigroup semiring.

19. Give an example of an infinite Smarandache group semiring.

## 4.5 S-semiring of second level

As inspired and suggested by Minh Perez in this section, I venture to define what are called Smarandache semirings (S-semirings) of second level. All semirings studied in the sections 4.1 to 4.4 are S-semirings of first level. Even to define the very concept of second level S-semiring we define a new Smarandache mixed direct product.

**DEFINITION 4.5.1**: *Let $S_1$ and $S_2$ be two different algebraic structures. We define Smarandache mixed direct product (S-mixed direct product) as $S_1 \times S_2$ where $S_1 \times S_2$ = $\{(s_1, s_2) / s_1 \in S_1$ and $s_2 \in S_2\}$ so the S-mixed direct product will enjoy simultaneously or separately the properties of the algebraic structures $S_1$ and $S_2$.*

*Example 4.5.1*: Let $Z^o$ be the semiring. $Z_7$ be the field, $S = Z^o \times Z_7$ is a semiring, is the Smarandache mixed direct product. The semiring S contains $\{0\} \times Z_7 \subseteq Z^o \times Z_7$ as a field paving way for us to define the second level Smarandache semirings.

**DEFINITION 4.5.2**: *A semiring S is said to be a Smarandache semiring of II level or Smarandache semiring of level II (S-semiring of level II) if S contains a proper subset P, which is a field.*

*Example 4.5.2*: Let $S = C_{10} \times Q$ be the S-mixed direct product of the semiring $C_{10}$ and the field Q. Clearly S is a semiring and S is in fact a S-semiring of II level.

*Examp4le 4.5.3*: Let $S = Z^o \times Q$. This is also a semiring, which is a S-semiring of level II.

We see that if S is just a semiring not got as a S-mixed direct product we see S cannot contain a proper subset which is a field. Clearly S-semiring of level II are richer structures than the S-semiring of level I.

*Example 4.5.4*: Let $S = Z_5 \times C_{10}$ is a S-semiring of level II, S has S-ideals and S-subsemiring and S is also an S-idempotent semiring.

*Example 4.5.5*: $S = C_2 \times R$ is a S-semiring of level II, $S_1 = C_2 \times R^o$ is a S-subsemiring. $S_2 = C_2 \times Q^o$ is also a S-subsemiring but $S_3 = \{0\} \times Z^o$ is not a S-subsemiring. So we see as in the case of S-semirings of level I we can have subsemirings which may not in general be a S-subsemiring of level II.

**THEOREM 4.5.1**: *All S-semirings of level II are non-strict semiring.*

*Proof*: Given S is a S-semiring of level II we know this means S contains a proper subset P which is a field. So in P we have $x + y = 0$ with $x \neq 0$ and $y \neq 0$ as P is a field. So all S-semirings of level II are non strict semirings.

In view of this we see that only S-semirings of level II alone can give way to non strict semirings, as we do not know any S-semiring of level I which is not a strict semiring. Thus we see that there exist a clear distinction between S-semirings of level I and level II.

***Example 4.5.6***: S = $Z_7 \times C_3$ is a S-semiring of level II. Clearly (5, 0) + (2, 0) = (0, 0), (3, 0) + (4, 0) = (0, 0). In fact S has S-ideals of the from S = $Z_7 \times \{0, a\}$.

These examples forces us to define Smarandache substructures in a different way. We find no change between two levels in case of S-zero divisors, S-units and S-idempotents. But we see we have to redefine S-subsemirings and S-ideals.

**DEFINITION 4.5.3**: *Let S be a S-semiring of level II. we say a proper subset A of S will be called a Smarandache subsemiring of level II (S-subsemiring II) if A contains a proper subset P which is a field.*

***Example 4.5.7***: Let S = $L \times Z_{11}$, where L is the lattice with Hasse diagram and $Z_{11}$ the prime field of characteristic 11. S is a S-semiring of level II. A = $\{1, 0\} \times Z_{11}$ is a S-subsemiring of level II. A S-semiring of level II can contain S-subsemiring of level I also.

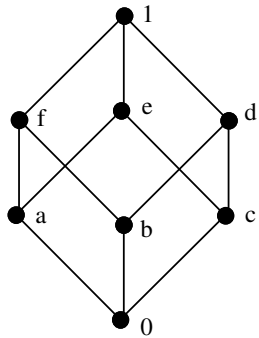

**Figure 4.5.1**

***Example 4.5.8***: Let S = $C_8 \times Z$, S is not a S-semiring of level II. The set A = $\{0, 1\} \times 2Z$ is a S- subsemiring for $\{0, 1\} \times \{0\}$ is a semifield.

In view of this we have the following theorem:

**THEOREM 4.5.2**: *If S is a semiring of level II, in general S can have both S-subsemirings of level I and II.*

*Proof*: Obvious by examples.

**THEOREM 4.5.3**: *If S is a S-semiring of level I which is not of level II. S cannot have a S-subsemiring of level II.*

*Proof*: If the S-semiring has a S-subsemiring of Level II then S contains a proper subset A such that A has a proper subset which is a field so this in turn will make S a S-semiring of level II a contrary to our assumption S is a S-semiring of level I only and not of level II. Hence the claim.

**DEFINITION 4.5.4**: *Let S be a semiring, which is a S-semiring of level II. A proper subset A of S is said to be a Smarandache right ideal of level II (S-right (Left) ideal II) if A is a S-subsemiring of level II and for all a ∈ A and for all s ∈ P ⊂ A (where P is a field) we have sa ∈ P (or as ∈ P).*

*If S is a S-semiring and A is simultaneously a S-right ideal II and a S-left ideal II we say A is a Smarandache ideal of level II (S-ideal II).*

***Example 4.5.9***: Let $S = Z_{10} \times Z^o$, S is a S-semiring of level II for A = {0, 5} × {0} is a field and (5, 0) is its identity. Also $A_1$ = {0, 2, 4, 6, 8} × {0} is a field with {6, 0} as the identity. We can get S-ideal II from these semirings as {0, 5} × $2Z^o$ is an S-ideal II of S, Now P = {0, 6} × $3Z^o$ is another S-ideal II of S.

***Example 4.5.10***: $S = C_6 \times Z_{12}$ is S-semiring of level II for {0} ×{0, 4, 8} is a subfield of S.

**THEOREM 4.5.4**: *Let S be S-semiring of level II given by $S = C_n \times Z_p$ where $Z_p$ is a prime field of characteristic p. Then S is a Smarandache commutative congruence simple finite semiring.*

*Proof*: $S = C_n \times Z_p$ has a S-subsemiring which is a Smarandache-commutative congruence simple finite semiring. (For more details refer Chris Monico).

<u>*Remark*</u>: Only S-semiring of II level can be S commutative c-simple finite semiring.

**DEFINITION 4.5.5**: *Let S be a S-semiring of level II. S is said to satisfy Smarandache assending chain condition II (S-accII) if $S_1 \subset S_2 \subset …$ is a monic ascending chain of S-ideals $S_i$ (of level II) and there exists a positive integer r such that $S_r = S_s$ for all $s \geq r$.*

Thus we see S-acc II is different from S-acc. The reader to requested to construct examples in S-semirings of level II which have both sets of S-acc's that is S-acc and S-acc II are satisfied. Clearly a S-semiring of level I can never satisfy S-acc II.

Similarly we define S-dcc II for S-semirings of level II. Here also a S-semiring of level I can never satisfy the S-dcc II condition. Now we proceed on to define compact semiring of level II.

**DEFINITION 4.5.6**: *Let S be a S-semiring of Level II. S is said to be a Smarandache compact semiring of level II (S-compact II semiring), if A ⊂ S where A is a S-subsemiring II of S is a compact subsemiring under the operations of S. Here also if S is a S-semiring of level I then S cannot be a S-compact II semiring.*

**DEFINITION 4.5.7**: *Let S be a semiring, which is a S-semiring of II level. S is said to be Smarandache * semiring of level II (S-* semiring II) if S contains a proper subset A satisfying the following conditions*

1. *A is a subsemiring of S*
2. *A is a S-subsemiring of Level II.*
3. *A is a *-semiring.*

*So if S is a S-∗semiring II then S need not be a S ∗-semiring.*

**DEFINITION 4.5.8**: *Let S be a S-semiring of level II. S is said to be Smarandache inductive ∗ semiring II (S inductive ∗ -semiring II) if S contains a proper subset A satisfying the following conditions:*

1. *A is a subsemiring of S.*
2. *A is a S-subsemiring II*
3. *A is a inductive ∗- semiring*

*It is easily verified a S-semiring of level I can never be a S-∗ semiring of level II or S-inductive ∗ -semiring of level II.*

**DEFINITION 4.5.9**: *Let S be a S-semiring of level II. S is said to be Smarandache continuous semiring of level II (S-continuous semiring II) if a proper subset A of S satisfies the following condition:*

1. *A is a S-subsemiring of level II*
2. *A is a continuous semiring*

*In case of S semiring of level II also, the definition remains the same as that of the definition of S-idempotent semiring given in S-semiring of level I. While defining for S-semiring of level II the concept of Smarandache e-semiring we replace the subsemiring A which is a S-subsemiring by S-subsemiring of level II and A is a e-semiring. Thus with these definitions about substructures in S-semiring of level II we propose the following problems for the reader to solve.*

**Notation**: A semiring of all types in level two will shortly be denoted by S-semiring II for example S-continuous semiring II etc.

**PROBLEMS**:

1. Give an example of a S-semiring of level II of
    i. Finite order.
    ii. Infinite order.

2. Find for the S-semiring, $S = Z_{15} \times C_5$
    i. S-subsemiring II.
    ii. S-ideal II.

3. Can the S-semiring II where $S = Z_{17} \times C_8$ be a S-continuous semiring II?

4. Is the S-semiring II where $S = Q \times C_2$ be a S-continuous semiring II?

5. Is the semiring $S = Z_{12} \times G$ be a S-e semiring II?

6. Find a S-semiring II which is a S-compact semiring II.

7. Give an example of a S-∗ semiring II.

8. Find an example of S-semiring II which is a S-inductive ∗ semiring II.

9. Give an example of S-semiring II which is a S-acc II semiring.

10. Give an example of a S-semiring II which is a
    i. S-dcc II semiring
    ii. S-MC II semiring
    iii. S-mc II semiring.

11. Give an example of a S-semiring II which is a S-c-semiring II.

12. Can a S-semiring II be not a S-dcc II semiring?

## 4.6 Smarandache Anti Semiring

Here we introduce yet an interesting property of semiring viz. Smarandache anti semiring. Florentin Smarandache has introduced a new concept called Smarandache anti structures. A set that is a strong structure contain a proper subset that has a weaker structure, for example if G is a group, we consider a subset S of G that is a semigroup, for the stronger structure groups contain subsets which are semigroups which are known as Smarandache anti semigroups. Suppose Z denotes the group under + we see $Z^+$ the set of integers without zero is a semigroup. Like wise we introduce here the concept of Smarandache anti semiring.

**DEFINITION 4.6.1**: *Let R be a ring. R is said to be Smarandache anti semiring (S-anti semiring) if R contains a subset S such that S is just a semiring.*

*Example 4.6.1*: Let Z be the ring. Z is an S-anti semiring for $Z^+$ the set of positive integers is a semiring.

*Example 4.6.2*: Let Q be the field of rationals, $Q^+$ is the semiring. So Q is S-anti semiring.

*Example 4.6.3*: Let R be the field of reals, $R^+$ is a semiring so R is a S-anti semiring.

*Example 4.6.4*: C be the field of complex numbers. This has subsets $Z^+$, $Q^+$ and $R^+$ to be semirings. Hence C is a S-anti semiring.

All these are examples of commutative rings of infinite order, that is of characteristic 0. Now we proceed on to study non-commutative rings with characteristic 0.

*Example 4.6.5*: Let $M_{3\times3} = \{(a_{ij}) / a_{ij} \in Z;$ the ring of integers$\}$ be the set of all $3 \times 3$ matrices under matrix addition and matrix multiplication. $M_{3\times3}$ is a non-commutative ring of characteristic 0. Clearly, $M_{3\times3}$ is a S-anti-semiring as $M_{3\times3} = \{(a_{ij}) / a_{ij} \in Z^o\}$ is a semiring.

***Example 4.6.6***: Let $Q[x]$ be the polynomial ring. The subset $P = \{Q^+[x] / Q^+$ is the positive rationals is a semiring$\}$, So $Q[x]$ is a S-anti semiring.

***Example 4.6.7***: Let M be any modular lattice having the following Hasse diagram

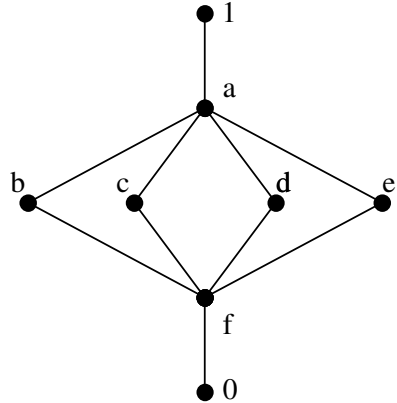

**Figure 4.6.1**

$S = \{1, a, b, f, 0\}$ is a semiring, so can we say M is a S-anti semiring of finite order? Note: A modular lattice can never be a ring.

***Example 4.6.8***: Let $P_{3\times 3} = \{(a_{ij}) / a_{ij} \in M\}$ where M is a lattice having the following Hasse diagram:

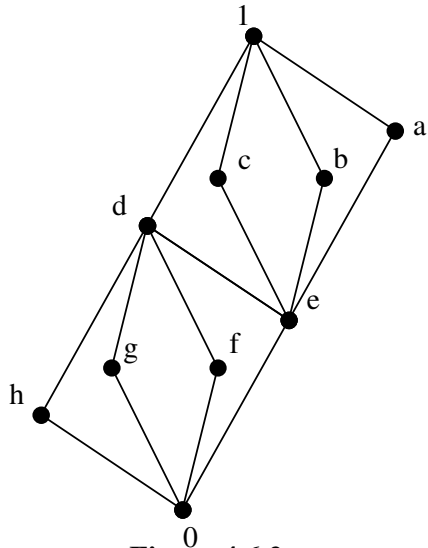

**Figure 4.6.2**

Let $H = \{0, e, a, 1\}$ is a distributive lattice. Let $M_{3\times 3} = \{(a_{ij}) / a_{ij} \in H\}$. Can we say $M_{3\times 3}$ is a S-anti semiring? From these examples we see a distributive lattice is a semiring and the class of distributive lattices is contained in the class modular lattices are never rings. So the concept of S-anti semirings using distributive lattices cannot be defined. Here $M_{3\times 3}$ is a non-commutative algebraic structure. Now we have the following open problem. Can we say all rings are S-anti semirings? The answer to this question is no, for when we take the rings $Z_n = \{0, 1, 2, \ldots, n-1\}$, n any positive integer we see $Z_n$ has no subset which is a semiring. So $Z_n$ for no n is a S anti semiring.

*Example 4.6.9*: Let Q be the field of rationals. $G = <g \mid g^3 = 1>$. QG be the group ring of G over K. Take $Q^+G = \{\sum \alpha_i g_i / \alpha_i \in Q^+ \text{ and } g_i \in G\}$, $Q^+G$ is a semiring. So QG is a S-anti semiring.

**THEOREM 4.6.1**: *Let F be a field of characteristic zero and G any group, the group ring FG is a S-anti semiring.*

*Proof*: Since F is a field of characteristic zero so $Q \subset F$ or $Q = F$. In both cases we have $Q^+$ is a semiring. Thus $Q^+G$ is also semiring. Hence FG is a S-anti semiring.

*Example 4.6.10*: $Z_nG$ be the group ring of the group G over the ring $Z_n$. Is the group ring $Z_nG$ a S-anti semiring? Left for the reader to prove.

**PROBLEMS**:

1. Give an example of a finite S-anti semiring.

2. Prove if F is a field of characteristic 0. $F \times \ldots \times F$ is a S-anti semiring.

3. Let F be a field of charcteristic p, p a prime. Can $S = Z_p \times \ldots \times Z_p$ (n times) be a S-anti semiring?

4. Find whether the following lattice with the Hasse diagram is a S-anti semiring

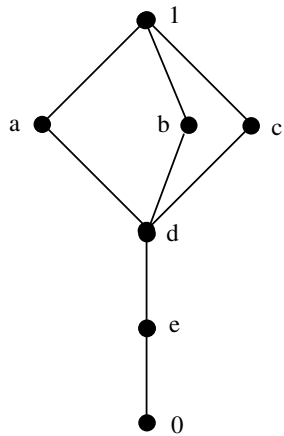

**Figure 4.6.3**

5. Can $M_{2\times 2} = \{(a_{ij}) / a_{ij} \in Z_3\}$ the ring of 2×2 matrices be a S-anti semiring? Justify your answer.

6. Is $Z_5S_3$, a S-anti semiring?

7. Can $Z_{12}S(4)$ be a S-anti semiring?

8. Can the ring $R = Z_5 \times Z_{12} \times Z_{10}$ be a S-anti semiring?

9. Prove $ZS(3)$ is a S-anti semiring.

10. Will ever $S = Z_5 \times Z_5 \times Z_5$ be a S-anti semiring? Justify your answer?

## Supplementary Reading

# CHAPTER FIVE
# SMARANDACHE SEMIFIELDS

In this chapter we introduce the concept of Smarandache semifields (S-semifield) and obtain some interesting results about them. Semifields can be of characteristic 0 or have no characteristic associated with them. We define S-substructure in S-semifields. We also define in this chapter S-semifields of Level II and finally introduce the concept of S-anti semifields.

## 5. 1 Definition and examples of Smarandache semifields

This section is devoted to the definition of Smarandache semifields and contains illustrative examples.

**DEFINITION 5.1.1**: *Let S be a semifield. S is said to be a Smarandache semifield (S-semifield) if a proper subset of S is a k-semi algebra, with respect to the same induced operations and an external operator.*

*Example 5.1.1*: Let $Z^o$ be a semifield. $Z^o$ is a S-semifield for A = {0, p, 2p, ..} is a proper subset of $Z^o$ which is a k-semi algebra.

*Example 5.1.2*: Let $Z^o[x]$ be a semifield. $Z^o[x]$ is a S-semifield as $pZ^o[x]$ is a proper subset which is a k-semi algebra.

It is important to note that all semifields need not be S-semifields.

*Example 5.1.3*: Let $Q^o$ be the semifield, $Q^o$ is not a S-semifield.

**THEOREM 5.1.1.** *Let S be the semifield. Every semifield need not be a S-semifield.*

*Proof:* It is true by the above example, as $Q^o$ is a semifield which is not a S-semifield.

All the while we have introduced only S-semifield of characteristic zero. Now we will proceed onto define S-semifields which has no characteristic associated with it.

*Example 5.1.4*: Let $C_n$ be a chain lattice. $C_n$ is a semifield. Any set of the form {0, $a_1$, …, $a_r$} such that $a_1 < a_2 < … < a_r$ and $a_r \neq 1$ is a k-semi algebra so $C_n$ is a semifield.

**THEOREM 5.1.2**: *All chain lattices $C_n$ are S-semifields.*

*Proof*: Obvious from the fact $C_n$ forms a semiring and has no zero divisors and has k-semi algebras in them.

The next natural question would be: Are all distributive lattices S-semifields?

In view of this we have an example.

*Example 5.1.5*: Let L be the lattice given by the following Hasse diagram.

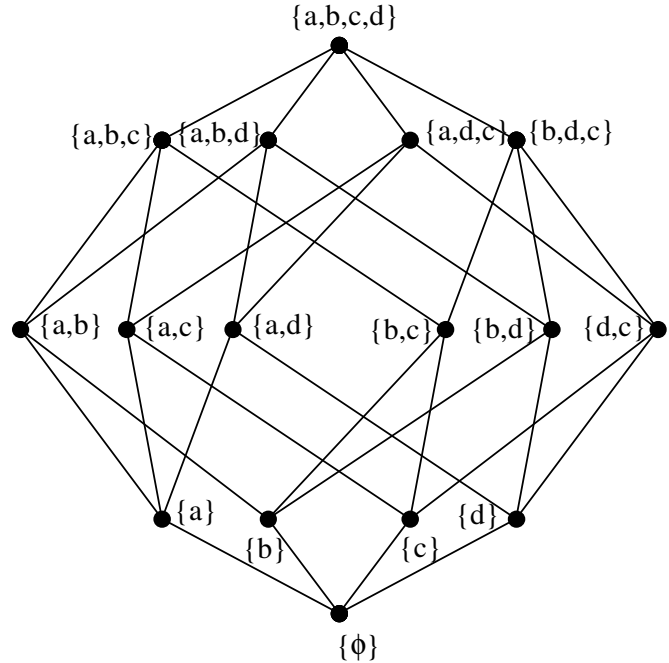

**Figure 5.1.1**

This lattice is distributive but not a semifield.

**THEOREM 5.1.3**: *All distributive lattices in general are not S-semifields.*

*Proof*: In view of the example 5.1.4. we see in general all distributive lattices are not S-semifields.

All distributive lattices are not S-semifields is untrue for we can have distributive lattices that are not chain lattices can be S-semifields.

*Example 5.1.6*: The following lattice with the Hasse diagram is a S-semifield.

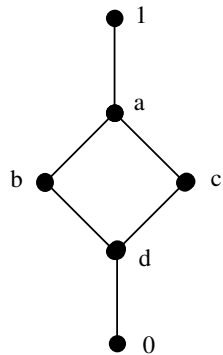

**Figure 5.1.2**

Hence the claim.

**PROBLEMS**:

1. Is $C_3$ a S-semifield?
2. Can a Boolean algebra of order greater than two be a S-semifield?
3. Is $R^o$ a S-semifield?
4. Can $Q^o[x]$ be a S-semifield?
5. Prove $C_2[x]$ is a S-semifield.
6. Prove $C_n[x]$ is a S-semifield.
7. Is $R^o[x]$ a S-semifield?

## 5.2 S-weak semifields

In this section we define a generalized concept of S-semifields viz. S-weak semifield and illustrate them with examples.

**DEFINITION 5.2.1**: *Let S be a semiring. S is said to be a Smarandache weak semifield (S-weak semifield) if S contains a proper subset P which is a semifield and P is a Smarandache semifield.*

Thus we see the following example is a S-weak semifield.

*Example 5.2.1*: Let S be a semiring given by the following Hasse diagram:

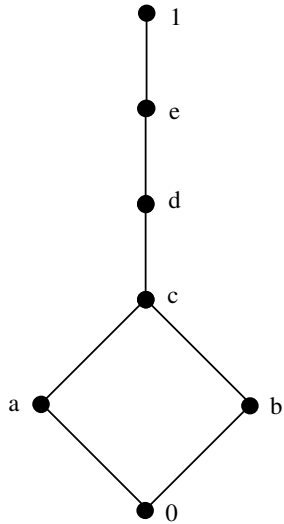

**Figure 5.2.1**

Clearly S is not a semifield as $a \bullet b = 0$ ($a \neq 0$, $b \neq 0$). S is only a semiring. Take P = {1, e, d, c, a, 0}. P is a semifield and P is in fact a S-semifield as T = {0, a, c, d, e} is a k-algebra over S. S is a S-weak semifield.
Thus we can say the following:

**THEOREM 5.2.1**: *Let S be a S-semifield. Then S is a S-weak semifield. But in general a S-weak semifield cannot be a S-semifield.*

*Proof*: Clearly by the very definitions of S-semifield and S-weak semifield we see every S-semifield is a S-weak semifield. But a S-weak semifield is not a S-semifield for if we take S to be just a semiring with zero divisor or a semiring which is non-commutative we see S cannot be a S-semifield. Example 5.2.1 is a S-weak semifield which is not a S-semifield. Thus we see we can also define S-weak semifield using non-commutative semirings.

*Example 5.2.2*: Let $S = Z^o \times Z^o \times Z^o$ be a semiring. S is a S-weak semifield. It is left for the reader to verify.

*Example 5.2.3*: Let $M_{2 \times 2} = \{(a_{ij}) / a_{ij} \in Z^o\}$, $M_{2 \times 2}$ is a semiring which is not a semifield. Take $P = \left\{ \begin{pmatrix} a & 0 \\ 0 & b \end{pmatrix} \bigg/ a, b \in Z^o \setminus \{0\} \right\} \cup \left\{ \begin{pmatrix} 0 & 0 \\ 0 & 0 \end{pmatrix} \right\}$. Clearly P is a semifield. Consider the set $A = \left\{ \begin{pmatrix} a & 0 \\ 0 & 0 \end{pmatrix} \bigg/ a \in Z^o \setminus \{0\} \right\} \cup \left\{ \begin{pmatrix} 0 & 0 \\ 0 & 0 \end{pmatrix} \right\}$. A is a k algebra over P. Hence the claim. So $M_{2 \times 2}$ is a S-weak semifield.

**PROBLEMS**:

1. Give an example of a S-weak semifield of finite order.
2. Give an example of a S-semifield of order 7.
3. Is $C_2 \times C_2 \times C_2 \times C_2 = S$, a S-semifield?
4. Can $S = C_2 \times Z^o$ be a S-weak semifield?
5. Show $S = Z^o \times Q^o$ is a S-weak semifield.
6. Give an example of a S-weak semifield of order 12.
7. Is S[x] where S is the lattice given by the Hasse diagram, a S-weak semifield?

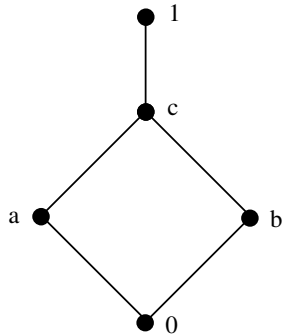

**Figure 5.2.2**

## 5.3 Special types of S-semifields

Most of the results given about the special semirings can be easily adopted to semifields.

**THEOREM 5.3.1**: *Let $C_n$ be the semifield. $C_n$ is a finite additively commutative S-c-simple semiring.*

*Proof*: From the fact that $C_n$ has 1 to be only additive absorbing element 1. (Results of Chris Monico).

**THEOREM 5.3.2**: *Let $C_n^t[x]$ be a semifield of all polynomials of degree $\leq t$. $C_n^t[x]$ is a finite additively S-commutative c-simple semiring.*

*Proof*: True from the fact $C_n^t[x]$ is additively idempotent.

*Example 5.3.1*: Let $C_n$ be a chain lattice with n elements. $C_n^m[x]$ be the set of all polynomials of degree less than or equal to m with coefficients from $C_n$. $C_n^m[x]$ is a finite Smarandache c-simple ring. Left as an exercise for the reader to verify.

## 5.4 Smarandache semifields of Level II

We have defined S-semifields in the earlier chapters and they will be called as S-semifield of level I and we proceed on to define the concept of Smarandache semifields of level II. We will make use of the Smarandache mixed direct product defined in Chapter 4 to define S-semifields of level II.

**DEFINITION 5.4.1**: *Let $S = C_n \times Z_p$ be the S-mixed direct product of the field $Z_p$ and the semifield $C_n$. Clearly $S = C_n \times Z_p$ contains subfields and subsemifields.*

**DEFINITION 5.4.2**: *Let S be a semifield. S is said to be a Smarandache semifield of level II (S-semifield II) if S contains a proper subset which is a field.*

Just as in the case of Smarandache semirings of level I and level II we have in the case of S-semifields of level I is different and disjoint from that of the S-semifield of level II. For this has motivated us to define in the next chapter the concept of Smarandache semivector spaces.

*Example 5.4.1*: Let $S = C_7 \times Z_5$. S is a S-semifield of level II.

**THEOREM 5.4.1**: *If S is a finite S-semifield of level II then S is a S-finite c-simple semiring.*

*Proof*: By the very definition of S to be a S-semifield of level II, S has a proper subset which is a field, since S is finite so is the field contained in it. Hence is view of Chris Monico, S is a S-c-simple semiring.

*Example 5.4.2*: Let $S = Z^o[x] \times Q$, S is a S-semifield of level II.

In case of fields we cannot define idempotents but in case of S-semifields we can have non-trivial idempotents also.

*Example 5.4.3*: Let $S = C_5 \times Q$. S is S-semifield of level II. $C_5$ has the following Hasse diagram. All elements of the form $(a_1, 1), (a_2, 1), (a_3, 1), (a_1, 0), (a_2, 0)$ and $(a_3, 0)$ are some of the idempotents in S.

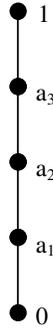

**Figure 5.4.1**

One more interesting property about S-semifields of order II is that S-semifields can have ideals.

**DEFINITION 5.4.3**: *Let S be a semifield. A proper subset P of S is said to be Smarandache- subsemifield of level I (S-subsemifield I) of S if P is a S-semifield of level I.*

In view of this we have the following theorem:

**THEOREM 5.4.2**: *Let S be a semifield. If S has a S-subsemifield of level I then S is a S-semifield of level I.*

*Proof*: Obvious by the very definition of S-subsemifields of level I, so S is a S-semifield of level I.

**DEFINITION 5.4.4**: *Let S be a semifield. A proper subset P of S is said to be a Smarandache subsemifield of level II if P is a S-semifield of level II.*

In view of this we have the following theorem.

**THEOREM 5.4.3**: *Let S be a semifield if S has subset P which is a S-subsemifield of level II then S is a S-semifield of level II.*

*Proof:* If S has a S-subsemifield of level II then S is a S-semifield of level II.

*Example 5.4.4*: Let $S = Z^o \times R$ is a S-semifield of level II of characteristic 0.

*Example 5.4.5:* Let $S = L \times R$ where L is a finite distributive lattice and R a field of characteristic 0 is a S-semifield of level II and this semiring has no characteristic associated with it.

## 5.5 Smarandache Anti-semifield

In this section we define Smarandache anti-semifields and study their properties.

**DEFINITION 5.5.1**: *Let S be a field or a ring. S is said to be a Smarandache anti semifield (S-anti semifield) is S has a proper subset A which is a semifield.*

*Example 5.5.1*: Q is field. $A = Q^o \subset Q$ is a semifield so Q is a S-anti semifield.

*Example 5.5.2*: Let Z be the ring of integers. Z is a S-anti semifield for $Z^o \subset Z$ is a semifield.

*Example 5.5.3*: $M_{3\times 3} = \{(a_{ij})/ a_{ij} \in Q\}$ be the ring of n × n matrices. $M_{3\times 3}$ is a S-anti semifield. For $S = \left\{ \begin{pmatrix} a_{11} & 0 & 0 \\ 0 & 0 & 0 \\ 0 & 0 & 0 \end{pmatrix} \middle/ a_{11} \in Q \right\}$ is a semifield. So $M_{3\times 3}$ is a S-anti semifield.

**DEFINITION 5.5.2**: *Let S be a ring or a field. A proper subset P in S is said to be a Smarandache anti-subsemifield (S-anti subsemifields) of S if P is itself a S-anti semifield.*

**THEOREM 5.5.1**: *If a ring or a field S has a S-anti subsemifield then S is a S-anti semifield.*

*Proof*: Obvious by the very definition of S-anti semifields and S-anti subsemifields.

*Example 5.5.4*: Let $Z_7$ be a field. $Z_7$ is not a S-anti semifield.

In view of this we get the following.

**THEOREM 5.5.2**: *All fields/rings are not in general S-anti semifields.*

*Proof*: By an example. Consider the collections of prime fields of characteristic p. p a prime. None of them are S-anti semifields.

**THEOREM 5.5.3**: *All fields of characteristic zero are S-anti semifields.*

*Proof:* Since F is a field of characteristic 0, we see Q the prime field of characteristic zero is either contained in F or F = Q. In both the cases we see F is a S-anti semifield as $Z^o \subset F$ or $Q^o \subset F$ are semifields; so F is a S-anti semifield.

**THEOREM 5.5.4**: *All rings S, commutative or non-commutative with unit 1 and characteristic 0 is a S-anti semifield.*

*Proof*: Since $1 \in S$ we see $Z \subset S$, as S is a ring of characteristic 0. Now $Z \subset S$ so $Z^o \subset Z$, is a semifield hence S is a S-anti semifield.

*Example 5.5.5*: Let F[x] be the polynomial ring. F is a ring or field of characteristic 0. Clearly F[x] is a S-anti semifield as $Z^o \subset F[x]$, is a semifield of P[x].

We now proceed on to define S-anti ideals in S-anti semifields.

**DEFINITION 5.5.3**: *Let S be a field/ ring which is a S-anti semifield. If we can find a subset P in the subsemifield T in S such that*

1. *P is a semiring.*
2. *for all $p \in P$ and $t \in T$, $pt \in P$.*

*Then P is called the Smarandache anti ideal (S-anti ideal) of the S-anti semifield. Note we cannot have the concept of right or left ideal as the subsemifield is commutative.*

*Example 5.5.6*: Let Q be the field. Q is a S-anti semifield. Clearly $pZ^o = \{0, p, 2p, \ldots\}$ is a S-anti ideal of Q.

*Example 5.5.7*: Let Q[x] be the polynomial ring. Q[x] is a S-anti semifield and $(pZ^o)[x]$ is a S-anti ideal of Q[x].

Thus we see even fields can have S-anti ideals in them

**PROBLEMS**:

1. Find a S-anti semifield of finite order.
2. Is $Q^o[x]$ a S-anti semifield? Justify your answer?
3. Give an example of a S-anti semifield using direct product of $Z \times Z_3$.
4. What is the characteristic of the S-anti semifield $Q \times Z_{12} \times Z_7$?
5. Show that S-anti semifields may or may not have characteristic associated with them.
6. Does there exist a S-anti semiring of characteristic n. n a positive finite integer?
7. Find an S-anti ideal in $S = Q \times Z_3 \times Z_7$.
8. Can $M_{3 \times 3} = \{(a_{ij})/ a_{ij} \in Q\}$ have S-anti ideals?
    i)   if so find them.
    ii)  find S-anti subsemifields which are not S-anti ideals.

## Supplementary Reading

1. Vasantha Kandasamy, W. B. *Semivector spaces over semifield.* Mathematika 188, 43-50, 1993.

2. Vasantha Kandasamy, W.B. *On a new class of semivector spaces.* Varahmihir Jour of Math. Sci., Vol. 1, 23-30, 2001.

3. Vasantha Kandasamy, W. B. *Smarandache Semirings and Semifields.* Smarandache Notions Journal, Vol. 7, 1-2-3, 88-91, 2001.
   http://www.gallup.unm.edu/~smarandache/SemiRings.pdf

# CHAPTER SIX

# SMARANDACHE SEMIVECTOR SPACES AND ITS PROPERTIES

In this chapter we introduce the concept of Smarandache semivector spaces and study some of its properties. We know every Smarandache semigroup (S-semigroup) is a semigroup but every semigroup in general is not a S-semigroup. Now we make use of Smarandache semigroup to construct Smarandache semivector spaces. Throughout this paper we assume all semigroups, which we construct, are going to be only semigroups under addition.

## 6.1 Definition of Smarandache semivector spaces with examples

Here we introduce the concept of Smarandache semivector spaces and show all Smarandache semivector spaces are semivector spaces but not conversely.

*Example 6.1.1*: Let $S = Z \times Z^o = \{(a, b) / a \in Z \text{ and } b \in Z^o\}$. Clearly $S = Z \times Z^o$ is a semigroup under component wise addition. In fact this semigroup is a Smarandache semigroup.

*Example 6.1.2*: $S = Z^o \times Z^o$ is not a S-semigroup.

**DEFINITION 6.1.1**: *Let G be a semigroup under the operation +, S any semifield. Let G be a semivector space (S-semivector space) over S. G is said to be a Smarandache semivector space (S-semivector space) over S if G is a Smarandache semigroup (S-semigroup).*

*Example 6.1.3*: Let $S = Q \times Z^o$ be a semigroup under component wise addition. S is a semivector space over $Z^o$ the semifield. Now we see S is a S-semivector space over $Z^o$. It is important to note $S = Q \times Z^o$ is not a semivector space over the semifield $Q^o$.

*Example 6.1.4*: Let $Q^o \times Q^o \times Q^o = S$ be a semigroup under component wise addition. Clearly S is a semivector space over $Q^o$ but S is not a S-semivector space as $S = Q^o \times Q^o \times Q^o$ is not a S-semigroup.

**THEOREM 6.1.1**: *All S-semivector spaces over a semifield S are semivector spaces but all semivector spaces need not be S-semivector spaces.*

*Proof*: By the very definition of S-semivector spaces we see all S-semivector spaces are semivector spaces. We note that all semivector spaces need not in general be S-semivector spaces as seen from example 6.1.4.

*Example 6.1.5*: Let $S = R^o \times Q^o \times Z$ be a S-semigroup. Clearly S is a S-semivector space over $Z^o$.

**Note**: $S = R^o \times Q^o \times Z$ is not even a semivector space over $Q^o$ or $R^o$.

**PROBLEMS**:

1. Give some semivector over $Z^o$.
2. Is $S = Z^o \times Z^o \times Q^o \times Q^o$ a S-semivector space over $Z^o$?
3. Can $S = Z^o \times Q^o \times R$ be a S- semivector space over $R^o$? Justify your answer.
4. Prove $S = Z^o \times Q^o \times R$ is a S-semivector over $Z^o$.
5. Why is $S = Z^o \times Q^o \times R$ a S-semivector over $Z^o$?

## 6.2 S- subsemivector spaces

Here we define the concept of Smarandache subsemivector spaces and give some examples. Further we define the notion of linear combination and Smarandache linearly independent vectors in the case of S-semivector spaces.

**DEFINITION 6.2.1**: *Let V be a Smarandache semigroup which is a S-semivector space over a semifield S. A proper subset W of V is said to be Smarandache subsemivector space (S-subsemivector space) of V if W is a Smarandache subsemigroup or W itself is a S-semigroup.*

*Example 6.2.1*: Let $V = Q^o \times Z^o \times Z$, V is a S-semivector space over $Z^o$. $W = Q^o \times Z^o \times 2Z$ is a S-subsemivector space of V. In fact $W_1 = Q^o \times \{0\} \times Z \subseteq V$ is also a S-subsemivector space of V. But $W_2 = Q^o \times Z^o \times Z^o \subset V$ is not a S- subsemivector space of V over $Z^o$. But $W_2$ is a subsemivector space of V over $Z^o$.

**THEOREM 6.2.1**: *Let V be a S semivector space over the semifield F. Every S-subsemivector space over S is a subsemivector space over S. But all subsemivector spaces of a S- semivector space need not be S-subsemivector space over S.*

*Proof*: By the very definition of S-subsemivector spaces $W \subset V$ we see W is a subsemivector space of V. But every subsemivector space W of V in general is not a S-subsemivector space as is evidenced from example 6.2.1 the subsemivector space $W_2 = Q^o \times Z^o \times Z^o \subset V$ is not a S-subsemivector space of V. Hence the claim.

*Example 6.2.2*: Consider $V = Z \times Z^o$, V is a S-semigroup. V is a S-semivector space over $Z^o$. We see the set $\{(-1, 1), (1, 1)\}$ will not span V completely $\{(-1, 0) (1, 0), (0, 1)\}$ will span V. It is left for the reader to find out sets, which can span V completely. Can we find a smaller set, which can span V than the set, $\{(-1, 0), (1, 0), (0, 1)\}$?

Let V be any S-semivector space over the semifield S. Suppose $v_1, \ldots, v_n$ be n set of elements in V then we say $\alpha = \sum_{i=1}^{n} \alpha_i v_i$ in V to be a linear combination of the $v_i$'s. We see when V is just a semivector space given in chapter III we could find semivector spaces using finite lattices but when we have made the definition of S-semivector spaces we see the class of those semivector spaces built using lattices can never be S-semivector spaces as we cannot make even semilattices into S-semigroups as x . x = x for all x in a semilattice. So we are left only with those semivector spaces built using $Q^o$, $Z^o$ and $R^o$ as semifields.

*Example 6.2.3*: Let $V = Q \times Z^o$ be a semivector space over $Z^o$. Clearly V is a S semivector space. In fact V has to be spanned only by a set which has infinitely many elements.

*Example 6.2.4*: Let $V = Q \times Z^o \times R$ be a S-semigroup. We know V cannot be a S-semivector space over $Q^o$ or $R^o$. V can only be a S-semivector space over $Z^o$. We see the set, which can span V, is bigger than the one given in example 6.2.3.

*Example 6.2.5*: Let $V = Z \times Z^o$ be a S-semigroup. V is a S-semivector space over $Z^o$. Clearly $\{(-1, 0), (0, 1), (1, 0)\} = \beta$ spans V. Our main concern is that will $\beta$ be the only set that spans V or can V have any other set which span it. Most of these questions remain open.

*Example 6.2.6*: Let $V = Z^o[x] \times Z$ be a S-semigroup. V is a S-semivector space over $Z^o$. The set, which can span V, has infinite cardinality.

**DEFINITION 6.2.2**: *Let V be a S-semigroup which is a S-semivector space over a semifield S. Let $P = \{v_1, \ldots, v_n\}$ be a finite set which spans V and the $v_i$ in the set P are such that no $v_i$'s in P can be expressed as the linear combination of the other elements in $P \setminus \{v_i\}$. In this case we say P is a linearly independent set, which span V.*

**DEFINITION 6.2.3**: *Let V be a S-semigroup which is a S-semivector space over a semifield S. If only one finite set P spans V and no other set can span V and if the elements of that set is linearly independent, that is no one of them can be expressed in terms of others then we say V is a finite dimensional S-semivector space and the cardinality of P is the dimension of V.*

*We see in the case of S-semivector spaces V the number of elements which spans V are always larger than the number of elements which spans the semivector spaces, which are not S-semivector spaces.*

**DEFINITION 6.2.4**: *Let V be a semigroup which is a S-semivector space over a semifield S. A Smarandache basis for V (S-basis for V) is a set of linearly independent elements, which span a S-subsemivector space P of V, that is P, is a S-subsemivector space of V, so P is also a S-semigroup. Thus for any semivector space V we can have several S-basis for V.*

*Example 6.2.7*: Let $V = Z^o \times Z$ be a S-semivector space over Z. Let $P = \{0\} \times \{pZ\}$ be a S-subsemivector space of V. Now the S-basis for P is $\{(0, p), (0, -p)\}$. We see for each prime p we can have S-subsemivector space which have different S-basis.

<u>PROBLEMS</u>:

1. Find a S-subsemivector space of $V = Q^o \times Z$ over $Z^o$.

2. Is $V = R^o \times Z$ a S-semivector space over the semifield $Z^o$?

3. Find a S-basis for $V = Z^o \times Z \times Z$, V the S-semivector space over $Z^o$.

4. Can the space V given in problem 2 be spanned by a finite set?

5. Let $V = Q \times Q^o$ be a S-semivector space over $Q^o$. Does V have a S-subsemivector space?

6. Let $V = Q \times Q^o$ be a S-semivector space over $Z^o$. Find a set that spans V. Can V have S-basis?

7. Give an example of S-semivector space which has infinitely many S-basis.

8. Does there exists a S-semivector space with no S-basis?

9. Give an example of a S-semivector space which has only one S-basis.

10. Does there exist a S-semivector space for which both the basis and the S-basis coincide?

## 6.3 Smarandache linear transformations

In this section we define Smarandache linear transformation of S-semivector spaces and study them. Here we consider only S-semivector spaces over semifields.

**DEFINITION 6.3.1**: *Let V and W be any two S-semigroups. We assume $P \subset V$ and $C \subset W$ are two proper subsets which are groups in V and W respectively. V and W be S-semivector spaces over the same semifield F. A map $T: V \to W$ is said to be a Smarandache linear transformation (S-linear transformation) if $T(cp_1 + p_2) = cTp_1 + Tp_2$ for all $p_1, p_2 \in P$ and $Tp_1, Tp_2 \in C$ i.e. T restricted to the subsets which are subgroups acts as linear transformation.*

*Example 6.3.1*: Let $V = Z^o \times Q$ and $W = Z^o \times R$ be S –semivector spaces over the semifield $Z^o$. We have $P = \{0\} \times Q$ and $C = \{0\} \times R$ are subsets of V and W respectively which are groups under +. Define $T: V \to W$, a map such that $T(0, p) \to (0, 2p)$ for every $p \in P$. Clearly T is a S-linear transformation. We see the maps T need not even be well defined on the remaining parts of V and W. What we need is $T: P \to C$ is a linear transformation of vector spaces.

*Example 6.3.2*: Let $V = Q^o \times R$ and $W = Z^o \times Z$ be S-semigroups which are S-semivector spaces over $Z^o$. $T: V \to W$ such that $T:\{0\} \times R \to \{0\} \times Z$ defined by $T(0, r) = (0, 0)$ if $r \notin Z$ and $T(0, r) = (0, r)$ if $r \in Z$. It is easily verified T is a S-linear transformation.

*Example 6.3.3*: Let $V = Z^o \times Z^o \times Z^o$ be a semigroup under addition. Clearly V is a semivector space over $Z^o$ but V is never a S-semivector space.

In view of this we have got a distinct behaviour of S-semivector space. We know if F is a field $V = F \times F \times \ldots \times F$ (n times) is a vector space over F. If S is a semifield then $W = S \times S \times \ldots S =$ (n times) is a semivector over S. But for a S- semivector space we cannot have this for we see none of the standard semifields defined using $Z^o$, $Q^o$ and $R^o$ are S-semigroups. They are only semigroups under addition and they do not contain a proper subset which is a group under addition.

*Example 6.3.4*: Let $V = Z^o \times Q \times Z^o$ be a S-semivector space over $Z^o$. Clearly $Z^o \times Z^o \times Z^o = W \subset V$, W is a semivector space which is not a S-semivector space. We are forced to state this theorem.

**THEOREM 6.3.1**: *Let V be a S-semivector space over $Q^o$ or $Z^o$ or $R^o$, then we can always find a subspace in V which is not a S-semivector space.*

*Proof*: If V is to be a S-semivector space the only possibility is that we should take care to see that V is only a semigroup having a subset which is a group i.e. our basic assumption is V is not a group but V is a S-semigroup. Keeping this in view, if V is to be a S-semivector space over $Z^o$ (or $Q^o$ or $R^o$) we can have $V = Z^o \times Z^o \times Z^o \times Q \times R \times \ldots \times Z^o$ i.e. V has at least once $Z^o$ occurring or $Q^o$ occurring or $R^o$ occurring and equally it is a must that in the product V, either Z or Q or R must occur for V to be a S-semigroup. Thus we see if V is a S-semivector space over $Z^o$. Certainly $W = Z^o \times \ldots \times Z^o \subset V$ is a semivector space over $Z^o$ and is not a S-semivector space over $Z^o$. Hence the claim.

*Example 6.3.5*: Let $V = Z^o \times Q^o \times R$ be a S-semigroup. V is a S-semivector space over $Z^o$. Clearly $W = Z^o \times Z^o \times Z^o$ is a subsemivector space of V which is not a S-semivector space.

**THEOREM 6.3.2**: *Let $V = S_1 \times \ldots \times S_n$ is a S-semivector spaces over $Z^o$ or $R^o$ or $Q^o$ where $S_i \in \{Z^o, Z, Q^o, Q, R^o, R\}$.*

1. *If one of the $S_i$ is Z or $Z^o$ then V can be a S-semivector space only over $Z^o$.*

2. *If none of the $S_i$ is Z or $Z^o$ and one of the $S_i$ is Q or $Q^o$, V can be a S-semivector space only over $Z^o$ or $Q^o$.*

3. *If none of the $S_i$ is Z or $Z^o$ or Q or $Q^o$ only R or $R^o$ then V can be a S-semivector space over $Z^o$ or $Q^o$ or $R^o$.*

*Proof*: It is left for the reader to verify all the three possibilities.

**THEOREM 6.3.3**: *Let $V = S_1 \times \ldots \times S_n$ where $S_i \in \{Z^o, Z, Q^o, Q, R$ or $R^o\}$ be a S-semigroup.*

1. *If V is a S-semivector space over $Z^o$ then $W = Z^o \times \ldots \times Z^o$ (n times) is a subsemivector space of V which is not a S-subsemivector space of V.*

2. *If V is a S-semivector space over $Q^o$ then $W = Q^o \times \ldots \times Q^o$ (n times) is a subsemivector space of V which is not a S-subsemivector space of V.*

3. *If V is a S-semivector space over $R^o$ then $W = R^o \times \ldots \times R^o$ (n times) is a subsemivector space of V and is not a S-subsemivector space of V.*

*Proof*: Left for the reader to do the proofs as an exercise.

**THEOREM 6.3.4**: *Let $V = S_1 \times \ldots \times S_n$ where $S_i \in \{Z^o, Z, R^o, R, Q^o, Q\}$ if V is a S-semivector space over $Q^o$. Then $W = Z^o \times \ldots \times Z^o$ ( n times) $\subset V$ is only a subset of V but never a subspace of V.*

*Proof*: Use the fact V is defined over $Q^o$ and not over $Z^o$.

We define a new concept called Smarandache pseudo subsemivector space.

**DEFINITION 6.3.2**: *Let V be a vector space over S. Let W be a proper subset of V. If W is not a subsemivector space over S but W is a subsemivector space over a proper subset $P \subset S$, then we say W is a Smarandache pseudo semivector space (S- pseudo semivector space) over $P \subset S$.*

*Example 6.3.6*: Let $V = Q \times R^o$ be a S-semivector space over $Q^o$. Clearly $W = Z^o \times R^o$ is not a subsemivector space over $Q^o$ but $W = Z^o \times R^o$ is a S- pseudo semivector space over $Z^o$.

*Example 6.3.7*: Let $V = Q^o \times R^o \times Q$ be a S-semivector space over $Q^o$. Now $W = Z^o \times Z^o \times Z^o$ and $W_1 = Q^o \times Q^o \times Q^o$ and S-pseudo semivector spaces over $Z^o \subset Q^o$.

Thus only these revolutionary Smarandache notions can pave way for definitions like Smarandache pseudo subsemivector spaces certainly which have good relevance.

**PROBLEMS**:

1. Find for $V = Z^o \times Q \times R^o$ a S-semivector space over $Z^o$
   a. S-subsemivector space
   b. S-pseudo subsemivector space

2. Can $V = Z^o \times Z$ over $Z^o$ have S-pseudo subsemivector space?

3. Find a S-pseudo subsemivector space $V = Q^o \times R$ over $Q^o$.

4. Can $V = Z^o \times Q^o \times R^o \times Z$ over $Z^o$ have S–subsemivector space? Find a S-pseudo subsemivector space.

5. Given $V = Z^o \times Q^o \times R^o \times Z$ over $Z^o$. Is V a semivector space? Is V a S-pseudo subsemivector space?

6. Define a S-linear transformation from V to W where $V = Z^o \times Q$ and $W = Q^o \times R$ are semivector spaces over $Z^o$.

7. Let $V = Z^o \times Q^o \times R$ and $W = Z \times Z^o \times Q$ be S-semivector spaces over $Z^o$. Define a S-linear transformation from V to W.

8. If we define a map ( V and W given in problem 7) T: V $\to$ W where $V_1 = \{0\} \times \{0\} \times R$ and $W_1 = \{0\} \times \{0\} \times Q$. T(0, 0, r) = (0, 0, 0) if r $\in$ R\Q and (0 0 r) if r $\in$ Q. Is T a S-linear transformation?

9. Let $V = Z^o \times Q$ and $W = Z \times Q^o$ be S-semivector spaces defined over $Z^o$. Define a S-linear transformation from V to W.

10. Let $V = Z \times Q^o \times R^o$ be a semivector space over $Z^o$. Find a S-linear operator on V.

## 6.4 S-anti semivector spaces

In this section we define Smarandache anti semivector spaces and obtain some interesting results about them.

**DEFINITION 6.4.1**: *Let V be a vector space over the field F. We say V is a Smarandache anti semivector space (S-anti semivector space) over F if there exists a subspace W $\subset$ V such that W is a semivector space over the semifield S $\subset$ F. Here W is just a semigroup under '+' and S is a semifield in F.*

*Example 6.4.1*: Let R be the field of reals. R is a vector space over Q. Clearly R is a S-anti semivector space as $Q^o \subset R$ is a S-semivector space over $Z^o$.

*Example 6.4.2*: Let $V = Q \times R \times Q$ be a vector space over Q. We see $W = Q^o \times R^o \times Q$ is a S-semivector space over $Q^o$. $W_1 = Z \times Z^o \times Z^o$ is not a S-semivector space over $Q^o$. But V is a S-anti semivector space over Q as $P = Z^o \times Z^o \times Z^o$ is a semivector space over $Z^o$.

*Example 6.4.3*: Let $V = Q \times Q \times \ldots \times Q$ (n-times), V is a vector space over Q. Clearly V is a S- anti semivector space for $Z^o \times Z^o \times \ldots \times Z^o \times Z$ is a S-semivector space over $Z^o$.

Many important questions are under study. The first is if V is a vector space over F and has a finite basis then it does not in general imply the S-anti semivector space has a finite basis. We have several facts in this regard, which are illustrated by the following examples.

*Example 6.4.4:* Let $V = Q \times Q \times Q \times Q \times Q$, (5 times) is a vector space over Q. Now $W = Z \times Z^o \times Z^o \times Z^o \times Z^o$ is a S-semivector space over $Z^o$. So V is a S-anti semivector space. The basis for $V = Q \times Q \times Q \times Q \times Q$ is {(1, 0, 0, 0, 0) (0, 1, 0, 0, 0), (0, 0, 1, 0, 0), (0, 0, 0, 0, 1), (0, 0, 0, 1, 0)} as a vector space over Q.

Now what is the basis or the set which spans $W = Z \times Z^o \times Z^o \times Z^o \times Z^o$ over $Z^o$. Clearly the set of 5 elements cannot span W. So we see in case of S-anti semivector spaces the basis of V cannot generate W. If we take $W_1 = Q^o \times Q^o \times Q^o \times Q^o \times Z$ as a S-semivector space over $Z^o$. Clearly $W_1$ cannot be finitely generated by a set. Thus a vector space, which has dimension 5, becomes infinite dimensional when it is a S-anti semivector space.

**DEFINITION 6.4.2**: *Let V and W be any two vector spaces over the field F. Suppose U ⊂ V and X ⊂ W be two subsets of V and W respectively which are S-semigroups and so are S-semivector spaces over S ⊂ F that is V and W are S-anti semivector spaces. A map T: V → W is said to be a Smarandache T linear transformation of the S-anti semivector spaces if T: U → X is a S-linear transformation.*

*Example 6.4.5:* Let $V = Q \times Q \times Q$ and $W = R \times R \times R \times R$ be two vector spaces over Q. Clearly $U = Z \times Z^o \times Z^o \subset V$ and $X = Q \times Z \times Z^o \times Z^o \subset W$ are S-semigroups and U and X are S-semivector spaces so V and W are S-anti semivector spaces. T: V → W be defined by T(x, y, z) = (x, x, z, z) for (x, y, z) ∈ $Z \times Z^o \times Z^o$ and (x, x, z, z) ∈ X is a Smarandache T linear operator.

Such nice results can be obtained using Smarandache anti semivector spaces.

**PROBLEMS**:

1. Let $V = Z^o \times Q$ and $W = Q^o \times R$ be S-semivector spaces over $Z^o$. Find S-linear transformation from V to W.

2. Let $V = Q \times Q \times R$ be vector space over Q. Is V a S-anti semivector space?

3. Can V and W in problem 1 be made into S-anti semivector space?

4. Find a basis for V given in problem 2.
    i. As a vector space
    ii. As a S-anti semivector space.

5. Let V = Q[x] be a vector space over Q can V be a S-anti semivector space?

    i. What is dimension of V as a vector space?
    ii. What is the dimension of V as a S-anti semivector space if we take $W = Q^o[x]$ and W is a semivector space over $Q^o$.

## Supplementary Reading

1. Vasantha Kandasamy, W. B. *Semivector spaces over semifield.* Mathematika 188, 43-50, 1993.

2. Vasantha Kandasamy, W.B. *On a new class of semivector spaces.* Varahmihir Jour of Math. Sci., Vol. 1, 23-30, 2001.

3. Vasantha Kandasamy, W. B. *Smarandache Semirings and Semifields.* Smarandache Notions Journal, Vol. 7, 1-2-3, 88-91, 2001.
   http://www.gallup.unm.edu/~smarandache/SemiRings.pdf

# CHAPTER SEVEN
# RESEARCH PROBLEMS

The main attraction of any textbook for any researcher in the list of open research problems proposed by it. In this book we have enlisted 25 open research problems for a student/ a researcher. Certainly these problems will throw open many more interesting results and this will certainly lead to a lot of researchers studying Smarandache notions. Smarandache notions are the only tools to study mathematics in unconventional ways and at the same time they can be used to analyse all types of relations between any two algebraic structures, this study will emerge as an attractive one among researchers. Except for Smarandache notions such rich type of mathematical analysis would be completely absent in the realm of mathematics.

1. Does there exist a semiring S of characteristic n, where n is a finite positive integer? (To show this one has to prove in a semiring S that for every $s \in S$; ns $= s + s + \ldots + s$ (n times) is zero)

2. Find a non-strict semiring S. (Hint: To prove this we have to show in a semiring S. $a, b \in S\setminus\{0\}$, $a + b = 0$ is possible.)

3. Let V be a semivector space over $S = Z^o$ where $V = Z^o \times Z^o \times \ldots \times Z^o$ (n times). Is it possible to find an upper bound for the number of linearly independent vectors in the semivector space V?

4. A simplified problem of problem 3 is if $V = Z^o \times Z^o$ is a semivector space over $Z^o$. Find an upper bound for the number of linearly independent elements in the semivector space $V = Z^o \times Z^o$ over $Z^o$.

5. Does there exist a semivector space V of dimension n over a semifield F such that every set of (n + 1) vectors is linearly dependent?

6. Let V be a semivector space over a semifield S. Let us assume V is not endowed with a unique basis. Let the number of elements in a basis be n (n > 1). Suppose V has another basis, then is it true in general that the number of elements in these two basis will be the same?

7. Characterize those semivector spaces which do not have a unique basis.

8. Characterize those semirings in which every S-subsemiring is an S-ideal.

9. Characterize those semirings in which all S-pseudo ideals are S-pseudo dual ideals.

10. Does there exist semirings in which all the four concepts i) S-ideals, ii) S-dual ideals, iii) S-pseudo ideals, iv) S-pseudo dual ideals coincide on every S-subsemiring?

11. Prove if $S = Z^o \times Z^o \times \ldots \times Z^n$ (n times) is a semiring. Then every zero divisor is a S-zero divisor.

12. Can we prove in any semiring every zero divisor is a S-zero divisor?

13. Let $S = Z^o \times Z^o \times \ldots \times Z^o$ (n times) be a semiring. S has idempotents but S has no S-idempotents. Is this true in case of all semifields of characteristic 0 ($Z^o$, $R^o$ and $Q^o$)?

14. Characterize those non-commutative semirings in which every S-subsemiring is a S-semidivision ring. (At least give an example of a non-commutative semiring in which every S-subsemiring is a S-semidivision ring).

15. Can one prove in any semiring the condition S-dcc and S-mc are equivalent? (Study this in case of i) Group semirings, ii) semigroup semirings, iii) S-group semirings, iv) S-semigroup semirings.)

16. Develop using the concept of Smarandache e-semiring

    a. Smarandache classical CSP (S-classical CSP)
    b. Smarandache fuzzy CSP (S-fuzzy CSP)
    c. Smarandache probabilistic CSP (S-probabilistic CSP)
    d. Smarandache weighted CSP (S-weighted CSP)

17. Can group semirings FG where F is a semiring of characteristic 0 (i.e. $Z^o$ or $Q^o$ or $R^o$) and G is a finite group be a S-chain semiring?

18. Find an example of a finite Smarandache semigroup semiring. (Hint: This is equivalent to finding a finite Smarandache anti-group. That is does there exist a finite group which is a Smarandache anti-group)

19. Does there exist a finite semiring S (semirings not constructed using S-mixed direct product) which has a proper subset which is a field?

20. Obtain some interesting results on S-inductive ∗ semirings.

21. Give interesting results on S-semivector spaces and S-basis.

22. Characterize those S-semivector spaces for which S-basis and basis coincide.

23. Characterize these S-anti semivector spaces V such that both V as a vector space as well as a S-anti semivector space have same dimension.

24. Can we have S-semivector space V having the same number of basis as the semivector space V?

25. Define Smarandache characteristic equation (S-characteristic equation), Smarandache eigen values (S-eigen values), Smarandache eigen vectors (S-eigen vectors) for S-linear operators of S-semivector spaces and obtain some interesting results about them.

# INDEX

*About the Author*

Dr. W. B. Vasantha is an Associate Professor in the Department of Mathematics, Indian Institute of Technology Madras, Chennai, where she lives with her husband Dr. K. Kandasamy and daughters Meena and Kama. Her current interests include Smarandache algebraic structures, fuzzy theory, coding/ communication theory. In the past decade she has guided seven Ph.D. scholars in the different fields of non-associative algebras, algebraic coding theory, transportation theory, fuzzy groups, and applications of fuzzy theory of the problems faced in chemical industries and cement industries. Currently, six Ph.D. scholars are working under her guidance. She has to her credit 241 research papers of which 200 are individually authored.  Apart from this, she and her students have presented around 262 papers in national and international conferences. She teaches both undergraduate and post-graduate students and has guided over 41 M.Sc. and M.Tech projects. She has worked in collaboration projects with the Indian Space Research Organization and with the Tamil Nadu State AIDS Control Society.

She can be contacted at vasantha@iitm.ac.in
You can visit her work on the web at: http://mat.iitm.ac.in/~wbv